\documentclass[amscd, leqno]{amsart}
\usepackage{amscd}
\theoremstyle{plain}
\newtheorem{thm}{Theorem}[section]
\newtheorem{Mthm}[thm]{Main Theorem}
\newtheorem{prop}[thm]{Proposition}
\newtheorem{lem}[thm]{Lemma}
\newtheorem{claim}[thm]{Claim}

\newtheorem{conj}[thm]{Conjecture}
\theoremstyle{definition}
\newtheorem{defn}[thm]{Definition}

\theoremstyle{remark}
\newtheorem{rem}[thm]{Remark}
\newenvironment{pf}{\begin{proof}}{\end{proof}}

\newcommand{\cal}{\mathcal}
\makeatletter

\@addtoreset{equation}{section}
\makeatother
\begin{document}
\title{Analytic torsion for Calabi-Yau threefolds}
\author{Hao Fang}
\address[Hao Fang]{
Department of Mathematics, 
University of Iowa, 
Iowa City, IA 52245, USA}
\email[Hao Fang]{haofang@math.uiowa.edu}
\author{Zhiqin Lu}
\address[Zhiqin Lu]{
Department of Mathematics, 
University of California Irvine,
Irvine, CA 92697, USA}
\email[Zhiqin Lu]{zlu@math.uci.edu}
\author{Ken-Ichi Yoshikawa}
\address[Ken-Ichi Yoshikawa]
{Graduate School of Mathematical Sciences,
University of Tokyo,
3-8-1 Komaba, Tokyo 153-8914, JAPAN}
\email[Ken-Ichi Yoshikawa]{yosikawa@ms.u-tokyo.ac.jp}
\thanks{The first-named author is partially supported 
by a grant from the New York University Research 
Challenge Fund Program and by NSF through
Institute for Advanced Study; 
the second-named author is partially supported by NSF 
Career Award DMS-0347033 and the Alfred P. Sloan
Research Fellowship; 
the third-named author is partially supported by
the Grants-in-Aid for Scientific Research 
for Encouragement of Young Scientists (B) 16740030, JSPS}

\begin{abstract}
After Bershadsky-Cecotti-Ooguri-Vafa, we introduce
an invariant of Calabi-Yau threefolds, which we call
the BCOV invariant and which we obtain 
using analytic torsion.
We give an explicit formula for the BCOV invariant
as a function on the compactified moduli space, 
when it is isomorphic to a projective line. 
As a corollary, we prove the formula for 
the BCOV invariant of quintic mirror threefolds 
conjectured by Bershadsky-Cecotti-Ooguri-Vafa.
\end{abstract}

\maketitle

{\bf Contents}
\begin{itemize}
\item[1.] 
Introduction
\item[2.] 
Calabi-Yau varieties with at most one ordinary
double point
\item[3.]
Quillen metrics
\item[4.] 
The BCOV invariant of Calabi-Yau manifolds
\item[5.] 
The singularity of the Quillen metric
on the BCOV bundle
\item[6.] 
The cotangent sheaf of the Kuranishi space
\item[7.] 
Behaviors of the Weil-Petersson metric
and the Hodge metric
\item[8.] 
The singularity of the BCOV invariant I
-- the case of ODP
\item[9.] 
The singularity of the BCOV invariant II
-- general degenerations
\item[10.]
The curvature current of the BCOV invariant
\item[11.]
The BCOV invariant of Calabi-Yau threefolds
with $h^{1,2}=1$
\item[12.]
The BCOV invariant of quintic mirror threefolds
\item[13.]
The BCOV invariant of FHSV threefolds
\newline
\end{itemize}

%%%%%%%%%%%%%%%%%%%%%%%%%%%%%%%%%%%%%%%%%%%%%%%%%%%%%%%%%%%
%
%                 Section 1
%
%%%%%%%%%%%%%%%%%%%%%%%%%%%%%%%%%%%%%%%%%%%%%%%%%%%%%%%%%%%

\section
{\bf Introduction}
\par
%%%%%%%%%%%%%%%%%%%%%%%%%%%%%%%%%%%%%%%%%%%%%%%%%%%%%%%%%%%
%
%                 (The BCOV invariant)
%
%%%%%%%%%%%%%%%%%%%%%%%%%%%%%%%%%%%%%%%%%%%%%%%%%%%%%%%%%%%
In the outstanding papers \cite{BCOV93}, \cite{BCOV94},
Bershadsky-Cecotti-Ooguri-Vafa made a deep study
of the generating function $F_{g}$ of genus-$g$ 
Gromov-Witten invariants for Calabi-Yau threefolds.
One mathematical surprise, which they obtained from 
physical arguments, is a system of holomorphic
anomaly equations satisfied by the functions 
$F_{g}$, $g\geq1$. From the holomorphic anomaly equations, 
they obtained a conjectural explicit formula 
for $F_{g}$ of a general quintic threefolds 
in ${\Bbb P}^{4}$ and thus they extended 
the mirror symmetry conjecture of 
Candelas-de la Ossa-Green-Parkes \cite{CdlOGP93}. 
\par
By focusing on
the genus-$1$ holomorphic anomaly equation, 
they conjectured that $F_{1}$ of a Calabi-Yau
threefold is expressed as a certain linear 
combination of the Ray-Singer analytic torsions 
(cf. \cite{BGS88}, \cite{RaySinger73}) 
of its mirror Calabi-Yau threefolds. 
After Bershadsky-Cecotti-Ooguri-Vafa,
we call the linear combination of Ray-Singer analytic 
torsions in \cite{BCOV94} the BCOV torsion,
which is the main subject of this paper.
\par
By making use of the curvature formula for
Quillen metrics \cite{BGS88}, 
Bershadsky-Cecotti-Ooguri-Vafa obtained
a variational formula for the BCOV torsion of
Ricci-flat Calabi-Yau manifolds \cite{BCOV94}. 
Fang-Lu \cite{FangLu03} expressed the variation
of the BCOV torsion of Ricci-flat Calabi-Yau 
manifolds as a linear combination of 
the Weil-Petersson metric \cite{Tian87}
and the generalized Hodge metrics \cite{Lu01}, 
which led them to some new results 
on the moduli space of 
polarized Calabi-Yau manifolds.
\par
On the other hand,
as a consequence of the duality in string theory,
Harvey-Moore \cite{HarveyMoore98} conjectured 
that the BCOV torsion of 
certain Ricci-flat Calabi-Yau threefolds 
is expressed as the product of the norms of 
the Borcherds $\Phi$-function \cite{Borcherds96}
and the Dedekind $\eta$-function.
Their conjecture was proved by
Yoshikawa \cite{Yoshikawa04}. In his proof,
an invariant of $K3$ surfaces with involution, 
which he obtained using
equivariant analytic torsion \cite{Bismut95} and 
a Bott-Chern class \cite{BGS88},
played a crucial role. 
\par
In this paper,
we extend the constructions of 
Bershadsky-Cecotti-Ooguri-Vafa and Yoshikawa to 
introduce a new invariant of Calabi-Yau threefolds,
which we call the BCOV invariant,
and we get an explicit formula for the BCOV invariant
as a function on the compactified moduli space 
when it is isomorphic to ${\Bbb P}^{1}$.
As a corollary of our formula, 
we prove one part of the conjecture
of Bershadsky-Cecotti-Ooguri-Vafa concerning
the BCOV torsion of quintic mirror threefolds.
Let us explain our results in more details.
\par
Let $X$ be a Calabi-Yau threefold.
Let $g$ be a K\"ahler metric on $X$ 
with K\"ahler form $\gamma$.  
We set $\overline{X}=(X,\gamma)$. 
Let $\tau(\overline{X},\overline{\Omega}^{p}_{X})$
be the Ray-Singer analytic torsion
of $\Omega_{X}^{p}=\wedge^{p}T^{*}X$ 
with respect to $g$.
We define the BCOV torsion of $\overline{X}$ as
$$
{\cal T}_{\rm BCOV}(\overline{X})
=
\prod_{p\geq0}
\tau(\overline{X},\overline{\Omega}^{p}_{X})^{(-1)^{p}p}.
$$
\par
Let
$\{{\bf e}_{1},\ldots,{\bf e}_{b_{2}(X)}\}$ 
be an integral basis of
$H^{2}(X,{\Bbb Z})/{\rm Torsion}$.
By Hodge theory and the Lefschetz decomposition
theorem, $H^{2}(X,{\Bbb R})$ is equipped
with the $L^{2}$-metric
$\langle\cdot,\cdot\rangle_{L^{2},[\gamma]}$,
which depends only on the K\"ahler class $[\gamma]$.
We define
$$
{\rm Vol}_{L^{2}}(H^{2}(X,{\Bbb Z}),[\gamma])
=
\det\left(
\langle{\bf e}_{i},{\bf e}_{j}\rangle_{L^{2},[\gamma]}
\right)_{1\leq i,j\leq b_{2}(X)},
$$
which is independent of the choice of 
an integral basis of 
$H^{2}(X,{\Bbb Z})/{\rm Torsion}$.
\par
Let $\eta$ be a nowhere vanishing holomorphic 
$3$-form on $X$. 
Let $c_{3}(X,\gamma)$ be the top Chern form 
of $(TX,g)$. We set 
${\rm Vol}(X,\gamma)=
(2\pi)^{-3}\int_{X}\gamma^{3}$ and
$\chi(X)=\int_{X}c_{3}(X,\gamma)$. We define
$$
{\cal A}(\overline{X})
=
{\rm Vol}(X,\gamma)^{\frac{\chi(X)}{12}}\,
\exp\left[-\frac{1}{12}\int_{X}
\log\left(
\frac{\sqrt{-1}\eta\wedge\bar{\eta}}
{\gamma^{3}/3!}
\cdot
\frac{{\rm Vol}(X,\gamma)}
{\|\eta\|_{L^{2}}^{2}}
\right)\,c_{3}(X,\gamma)\right],
$$
which is independent of the choice of $\eta$.
We define the real number $\tau_{\rm BCOV}(X)$ as
$$
\tau_{\rm BCOV}(X)
=
{\rm Vol}(X,\gamma)^{-3}\,
{\rm Vol}_{L^{2}}(H^{2}(X,{\Bbb Z}),[\gamma])^{-1}\,
{\cal A}(\overline{X})\,
{\cal T}_{\rm BCOV}(\overline{X}).
$$
In Sect.\,4.4, we show that 
$\tau_{\rm BCOV}(X)$ is independent of the choice
of $\gamma$. Hence $\tau_{\rm BCOV}(X)$ is 
an invariant of $X$, which we call 
the BCOV invariant. 
The purpose of this paper is to 
study $\tau_{\rm BCOV}$ 
as a function on the moduli of 
Calabi-Yau threefolds.

%%%%%%%%%%%%%%%%%%%%%%%%%%%%%%%%%%%%%%%%%%%%%%%%%%%%%%%%%%%
%
%                (Main Results)
%
%%%%%%%%%%%%%%%%%%%%%%%%%%%%%%%%%%%%%%%%%%%%%%%%%%%%%%%%%%%

\par
Let $\cal X$ be a (possibly singular) 
irreducible projective fourfold.
Let $\pi\colon{\cal X}\to{\Bbb P}^{1}$ be
a surjective flat morphism with discriminant locus 
${\cal D}$. 
Let $\psi$ be the inhomogeneous coordinate 
of ${\Bbb P}^{1}$, and set $X_{\psi}:=\pi^{-1}(\psi)$ 
for $\psi\in{\Bbb P}^{1}$.
We assume the following:
\newline{(i)}
${\cal D}$ is a finite subset of ${\Bbb P}^{1}$ 
such that $\infty\in{\cal D}$ and 
${\cal D}\setminus\{\infty\}\not=\emptyset$;
\newline{(ii)}
$X_{\psi}$ is a Calabi-Yau threefold with
$h^{2}(\Omega^{1}_{X_{\psi}})=1$ 
for $\psi\in{\Bbb P}^{1}\setminus{\cal D}$;
\newline{(iii)}
${\rm Sing}\,X_{\psi}$ consists of a unique 
ordinary double point (ODP)
for $\psi\in{\cal D}\setminus\{\infty\}$;
\newline{(iv)}
${\rm Sing}({\cal X})\cap X_{\infty}=\emptyset$
and $X_{\infty}$ is a divisor of normal crossing.
\par
Under these assumptions, the relative dualizing
sheaf $K_{{\cal X}/{\Bbb P}^{1}}$ is locally free
on ${\cal X}$, and its direct image sheaf
$\pi_{*}K_{{\cal X}/{\Bbb P}^{1}}$ is locally free
on ${\Bbb P}^{1}$.
\par
For $\psi\in{\Bbb P}^{1}\setminus\{\infty\}$, 
let $({\rm Def}(X_{\psi}),[X_{\psi}])$ be
the Kuranishi space of $X_{\psi}$.
Since $X_{\psi}$ is Calabi-Yau,
$\dim{\rm Def}(X_{\psi})=1$. We identify
$({\rm Def}(X_{\psi}),[X_{\psi}])$ with 
$({\Bbb C},0)$ by the smoothness of 
the Kuranishi space 
(cf. \cite{Tian87}, \cite{Tian92}, \cite{Todorov89}). 
Let $\mu_{\psi}\colon
({\Bbb P}^{1},\psi)\to({\rm Def}(X_{\psi}),[X_{\psi}])$
be the map of germs that induces the family
$\pi\colon{\cal X}\to{\Bbb P}^{1}$ near $\psi$.
The ramification index $r(\psi)$ of
$\pi\colon{\cal X}\to{\Bbb P}^{1}$ at 
$\psi\in{\Bbb P}^{1}$ is defined as the vanishing
order of $\mu_{\psi}$ at $\psi$.
Let $\{R_{j}\}_{j\in J}$ be the set of points 
of ${\Bbb P}^{1}$ with ramification index $>1$,
and write
${\cal D}\setminus\{\infty\}=\{D_{k}\}_{k\in K}$.
We set
$r_{j}=r(R_{j})$ for $j\in J$ 
and $r_{k}=r(D_{k})$ for $k\in K$.
\par
Outside ${\cal D}\cup\{R_{j}\}_{j\in J}$,
$T{\Bbb P}^{1}$ is equipped
with the Weil-Petersson metric.
Let $\|\cdot\|$ be the singular Hermitian metric on
$(\pi_{*}K_{{\cal X}/{\Bbb P}^{1}})^{\otimes(48+\chi)}
\otimes(T{\Bbb P}^{1})^{\otimes12}$ 
induced from the $L^{2}$-metric on 
$\pi_{*}K_{{\cal X}/{\Bbb P}^{1}}$
and from the Weil-Petersson metric on $T{\Bbb P}^{1}$.

\begin{Mthm}
Let $\varXi$ be a meromorphic section of
$\pi_{*}K_{{\cal X}/{\Bbb P}^{1}}$ with
$$
{\rm div}(\varXi)=
\sum_{i\in I}m_{i}\,P_{i}+m_{\infty}\,P_{\infty},
\qquad
P_{i}\not=P_{\infty}\,(i\in I).
$$ 
Identify the points $P_{i},R_{j},D_{k}$ 
with their coordinates 
$\psi(P_{i}),\psi(R_{j}),\psi(D_{k})\in{\Bbb C}$,
respectively.
Set $\chi=\chi(X_{\psi})$,
$\psi\in{\Bbb P}^{1}\setminus{\cal D}$.
Then there exists $C\in{\Bbb R}_{>0}$ 
such that
$$
\tau_{\rm BCOV}(X_{\psi})
=
C\,\left\|
\prod_{i\in I,j\in J,k\in K}
\frac{(\psi-D_{k})^{2r_{k}}}
{(\psi-P_{i})^{(48+\chi)m_{i}}
(\psi-R_{j})^{12(r_{j}-1)}}
\,\varXi_{\psi}^{48+\chi}
\otimes
\left(\frac{\partial}{\partial\psi}\right)^{12}
\right\|^{\frac{1}{6}}.
$$
\end{Mthm}

%%%%%%%%%%%%%%%%%%%%%%%%%%%%%%%%%%%%%%%%%%%%%%%%%%%%%%%%%%%
%
%                 (BCOV's conjecture)
%
%%%%%%%%%%%%%%%%%%%%%%%%%%%%%%%%%%%%%%%%%%%%%%%%%%%%%%%%%%%
\par
As a corollary of the Main Theorem 1.1,
we give a partial answer to the conjecture of 
Bershadsky-Cecotti-Ooguri-Vafa,
which we explain briefly (cf. Sect.\,12).
\par
Let $p\colon{\cal X}\to{\Bbb P}^{1}$ be the pencil of 
quintic threefolds in ${\Bbb P}^{4}$ defined by
$$
{\cal X}:=
\{([z],\psi)\in{\Bbb P}^{4}\times{\Bbb P}^{1};\,
z_{0}^{5}+z_{1}^{5}+z_{2}^{5}+z_{3}^{5}+z_{4}^{5}
-5\psi\,z_{0}z_{1}z_{2}z_{3}z_{4}=0
\},
\quad
p={\rm pr}_{2}.
$$
Let ${\Bbb Z}_{5}$ be the set of fifth roots
of unity and define
$$
G:=
\{(a_{0},a_{1},a_{2},a_{3},a_{4})\in({\Bbb Z}_{5})^{5};\,
a_{0}a_{1}a_{2}a_{3}a_{4}=1\}/{\Bbb Z}_{5}(1,1,1,1,1)
\cong
{\Bbb Z}_{5}^{3}.
$$
We regard $G$ as a group of projective 
transformations of ${\Bbb P}^{4}$.
Since $G$ preserves the fibers of $p$, 
we have the induced family
$p\colon{\cal X}/G\to{\Bbb P}^{1}$.
Let ${\cal D}$ be the discriminant locus of
the family $p\colon{\cal X}\to{\Bbb P}^{1}$.
By \cite{Batyrev94}, \cite{Morrison93},
there exists a resolution 
$q\colon{\cal W}\to{\cal X}/G$
such that $W_{\psi}=q^{1}(X_{\psi})$ 
is a smooth Calabi-Yau threefold for 
$\psi\in{\Bbb P}^{1}\setminus{\cal D}$
and such that
${\rm Sing}\,W_{\psi}$ consists of a unique ODP
if $\psi^{5}=1$.
The family of Calabi-Yau threefolds
$\pi\colon{\cal W}\to{\Bbb P}^{1}$ 
is called a family of quintic mirror threefolds.
\par
After 
Candelas-de la Ossa-Green-Parkes \cite{CdlOGP93}, 
$\pi_{*}K_{{\cal W}/{\Bbb P}^{1}}$ and
$T{\Bbb P}^{1}$ are trivialized as follows
near $\psi=\infty$.
For $\psi\in{\Bbb P}^{1}\setminus{\cal D}$,
we define a holomorphic $3$-form on $X_{\psi}$ by
$$
\varOmega_{\psi}
=
\left(\frac{2\pi\sqrt{-1}}{5}\right)^{-3}5\psi\,
\frac{z_{4}\,dz_{0}\wedge dz_{1}\wedge dz_{2}}
{\partial F_{\psi}(z)/\partial z_{3}}.
$$
Since $\varOmega_{\psi}$ is $G$-invariant,
$\varOmega_{\psi}$ induces a holomorphic $3$-form
on $X_{\psi}/G$ in the sense of orbifolds.
We identify $\varOmega_{\psi}$ with the corresponding
holomorphic $3$-form on $X_{\psi}/G$, and
we define a holomorphic $3$-form $\varXi_{\psi}$
on $W_{\psi}$ as
$\varXi_{\psi}=q_{\psi}^{*}\Omega_{\psi}$.
We define
$$
y_{0}(\psi)
=
\sum_{n=1}^{\infty}\frac{(5n)!}{(n!)^{5}(5\psi)^{5n}},
\qquad
|\psi|>1.
$$
Then 
$\pi_{*}K_{{\cal W}/{\Bbb P}^{1}}$ 
is trivialized by the local section 
$\varXi_{\psi}/y_{0}(\psi)$ near $\psi=\infty$.
\par
Let $q$ be the coordinate of the unit disc 
in ${\Bbb C}$.
We identify the parameters $\psi^{5}$ and $q$
via the mirror map \cite{CdlOGP93}.
Then $T{\Bbb P}^{1}$ is trivialized 
by the local section
$q\,d/dq=q\,(d\psi/dq)\,d/d\psi$ near $\psi=\infty$.
(See Sect.\,12.)
\par
We define a multi-valued analytic function 
$F_{1,B}^{\rm top}(\psi)$ near 
$\infty\in{\Bbb P}^{1}$ as
$$
F_{1,B}^{\rm top}(\psi)
=
\left(\frac{\psi}{y_{0}(\psi)}\right)^{\frac{62}{3}}
(\psi^{5}-1)^{-\frac{1}{6}}\,
q\,\frac{d\psi}{dq}
$$
and a power series in $q$ as
$F_{1,A}^{\rm top}(q)=F_{1,B}^{\rm top}(\psi(q))$.
The conjectures of Bershadsky-Cecotti-Ooguri-Vafa
\cite{BCOV93}, \cite{BCOV94} can be formulated
as follows:

\begin{conj}
{\bf (A)}
Let $N_{g}(d)$ be the genus-$g$
Gromov-Witten invariant of degree $d$ 
of a general quintic threefold in ${\Bbb P}^{4}$.
Then 
$$
q\frac{d}{dq}\log F_{1,A}^{\rm top}(q)
=
\frac{50}{12}
-
\sum_{n,d=1}^{\infty}N_{1}(d)\,
\frac{2nd\,q^{nd}}{1-q^{nd}}
-
\sum_{d=1}^{\infty}N_{0}(d)\,
\frac{2d\,q^{d}}{12(1-q^{d})}.
$$
\newline{\bf (B)}
The following identity holds near $\psi=\infty$:
$$
\tau_{\rm BCOV}(W_{\psi})
=
{\rm Const.}\,
\left\|\frac{1}{F_{1,B}^{\rm top}(\psi)^{3}}\,
\left(
\frac{\varXi_{\psi}}{y_{0}(\psi)}
\right)^{62}
\otimes\left(q\frac{d}{dq}\right)^{3}
\right\|^{\frac{2}{3}}.
$$
\end{conj}

In Sect.\,12, we prove the following:

\begin{thm}
The Conjecture 1.2 (B) holds.
\end{thm}

For the remaining Conjecture 1.2 (A),
see Li-Zinger \cite{LiZinger04}.
In \cite{FLY05},
we shall study the BCOV invariant of Calabi-Yau 
threefolds with higher dimensional moduli and
the BCOV torsion of Calabi-Yau manifolds of 
dimension greater than $3$.

%%%%%%%%%%%%%%%%%%%%%%%%%%%%%%%%%%%%%%%%%%%%%%%%%%%%%%%%%%%
%
%                 (Strategy of the proof)
%
%%%%%%%%%%%%%%%%%%%%%%%%%%%%%%%%%%%%%%%%%%%%%%%%%%%%%%%%%%%

\par
Let us briefly explain our approach to prove
the Main Theorem 1.1. We follow the approach
in \cite{Yoshikawa04}.
Let $\Omega_{\rm WP}$ be the Weil-Petersson form 
on ${\Bbb P}^{1}\setminus{\cal D}$, 
and let ${\rm Ric}\,\Omega_{\rm WP}$ be
the Ricci-form of $\Omega_{\rm WP}$. 
By \cite{Lu01}, \cite{LuSun04},
the $(1,1)$-forms $\Omega_{\rm WP}$ and
${\rm Ric}\,\Omega_{\rm WP}$ have Poincar\'e growth
on ${\Bbb P}^{1}\setminus{\cal D}$, 
so that they extend trivially to closed positive
$(1,1)$-currents on ${\Bbb P}^{1}$ (cf. Sect.\,7.3).
We identify $\Omega_{\rm WP}$ and
${\rm Ric}\,\Omega_{\rm WP}$ with their trivial
extensions. For a divisor $D$ on ${\Bbb P}^{1}$,
let $\delta_{D}$ denote the Dirac $\delta$-current
on ${\Bbb P}^{1}$ associated to $D$.
Regard $\tau_{\rm BCOV}$ as a function
on ${\Bbb P}^{1}\setminus{\cal D}$.
By making use of the Poincar\'e-Lelong formula, 
the Main Theorem 1.1 is deduced from 
the following:

\begin{claim}
Set ${\cal D}^{*}=\sum_{k\in K}r_{k}\,D_{k}$.
Then there exists $a\in{\Bbb R}$ 
such that 
\begin{equation}
dd^{c}\log\tau_{\rm BCOV}
=
-\left(\frac{\chi}{12}+4\right)\,
\Omega_{\rm WP}
-{\rm Ric}\,\Omega_{\rm WP}
+\frac{1}{6}\,\delta_{{\cal D}^{*}}
+a\,\delta_{\infty}.
\end{equation}
\end{claim}

We shall establish Claim 1.4 as follows:
\par
{\bf (a)}
By making use of the curvature formula for 
Quillen metrics of Bismut-Gillet-Soul\'e \cite{BGS88},
we prove the variational formula like (1.1)
for an arbitrary family of Calabi-Yau threefolds.
As a result, we get Eq.\,$(1.1)$
on the open part ${\Bbb P}^{1}\setminus{\cal D}$.
More precisely,
we introduce a Hermitian line, 
called the BCOV Hermitian line,
for an arbitrary Calabi-Yau manifold
of arbitrary dimension, which we obtain 
using determinants of cohomologies
\cite{KnudsenMumford76},
Quillen metrics \cite{BGS88}, \cite{Quillen85}, 
and a Bott-Chern class like ${\cal A}(\cdot)$. 
Then the BCOV Hermitian line of a Calabi-Yau 
manifold depends only on the complex structure 
of the manifold. 
The Hodge diamond of Calabi-Yau threefolds are
so simple that the BCOV Hermitian line
reduces to the scalar invariant 
$\tau_{\rm BCOV}$ in the case of threefolds. 
Hence Eq.\,$(1.1)$ 
on ${\Bbb P}^{1}\setminus{\cal D}$ is deduced 
from the curvature formula for the BCOV Hermitian 
line bundles. (See Sect.\,4).
\par
{\bf (b)}
To establish the formula for $\log\tau_{\rm BCOV}$
near ${\cal D}$, we fix a specific holomorphic 
extension of the BCOV bundle from 
${\Bbb P}^{1}\setminus{\cal D}$
to ${\Bbb P}^{1}$, which we call the K\"ahler
extension. (See Sect.\,5.)
Since $\tau_{\rm BCOV}$ is the ratio of 
the Quillen metric and the $L^{2}$-metric
on the BCOV bundle,  
it suffices to determine the singularities of
the Quillen metric and the $L^{2}$-metric
on the extended BCOV bundle.
We determine the singularity
of the Quillen metric on the extended BCOV bundle
with respect to the metric on 
$T{\cal X}/{\Bbb P}^{1}$ induced from 
a K\"ahler metric on ${\cal X}$.
The anomaly formula for Quillen metrics
of Bismut-Gillet-Soul\'e \cite{BGS88} 
and a formula for the singularity
of Quillen metrics \cite{Bismut97},
\cite{Yoshikawa05} play the central role. 
(See Sect.\,5.).
\par
{\bf (c)}
By the smoothness of ${\rm Def}(X_{\psi})$ at
$\psi\in{\cal D}^{*}$
\cite{Kawamata92}, \cite{Ran92}, \cite{Tian92},
the behavior of the $L^{2}$-metric
on the extended BCOV bundle near ${\cal D}^{*}$
is determined by the singularity of 
$\Omega_{\rm WP}$ near ${\cal D}^{*}$,
which was computed by Tian \cite{Tian92}.
(See Sects.\,6, 7, 8.)
To determine the behavior of the $L^{2}$ metric
on the extended BCOV bundle at $\psi=\infty$, 
one may assume that 
$\pi\colon{\cal X}\to{\Bbb P}^{1}$
is semi-stable at $\psi=\infty$ 
by Mumford \cite{Mumford73}. 
We consider another
holomorphic extension of the BCOV bundle, i.e.,
the canonical extension in Hodge theory 
\cite{Schmid73}. 
With respect to the canonical extension,
the $L^{2}$-metric has at most an algebraic 
singularity at $\psi=\infty$ 
by Schmid \cite{Schmid73}. 
Comparing the two extensions, we show that
the $L^{2}$-metric has at most an algebraic 
singularity at $\psi=\infty$
with respect to the K\"ahler extension.
(See Sect.\,9.)
By the residue theorem and assumption (ii),
the number $a$ in Eq.\,$(1.1)$ is determined
by the degrees of the divisors ${\cal D}^{*}$, 
${\rm div}(\varXi)$, $\sum_{j\in J}(r_{j}-1)\,R_{j}$.
(See Sect.\,11.)

%%%%%%%%%%%%%%%%%%%%%%%%%%%%%%%%%%%%%%%%%%%%%%%%%%%%%%%%%%%
%
%            (The organization of this paper)
%
%%%%%%%%%%%%%%%%%%%%%%%%%%%%%%%%%%%%%%%%%%%%%%%%%%%%%%%%%%%
This paper is organized as follows.
In Sect.\,2, we recall the deformation theory
of Calabi-Yau threefolds.
In Sect.\,3, we recall the definition of Quillen 
metrics and the corresponding curvature formula.
In Sect.\,4, we introduce the BCOV invariant and
prove its variational formula.
In Sect.\,5, we study the boundary behavior of 
Quillen metrics. 
In Sect.\,6, we study the boundary behavior
of Kodaira-Spencer map.
In Sect.\,7, we study the boundary behavior
of the Weil-Petersson metric and the Hodge metric.
In Sects.\,8 and 9, we study the boundary behavior
of the BCOV invariant.
In Sect.\,10, we extend the variational formula
for the BCOV invariant to the boundary of 
moduli space.
In Sect.\,11, we prove the Main Theorem.
In Sect.\,12, we study a conjecture of 
Bershadsky-Cecotti-Ooguri-Vafa.
In Sect.\,13, we study a conjecture of Harvey-Moore.

%%%%%%%%%%%%%%%%%%%%%%%%%%%%%%%%%%%%%%%%%%%%%%%%%%%%%%%%%%%
%
%            (Acknowledgements)
%
%%%%%%%%%%%%%%%%%%%%%%%%%%%%%%%%%%%%%%%%%%%%%%%%%%%%%%%%%%%
{\bf Acknowledgements}
The first-named author thanks Professors Jeff Cheeger
and Gang Tian for helpful discussions.
The second-named author thanks Professors Gang Tian
and Duong H. Phong for helpful discussions.
The third-named author thanks Professors Shinobu Hosono, 
Shu Kawaguchi, Yoshinori Namikawa and Gang Tian 
for helpful discussions, and his special thanks are 
due to Professor Jean-Michel Bismut, 
who suggested him, together with many other ideas,
one of the most crucial constructions in this paper, 
the Bott-Chern term ${\cal A}(\overline{X})$.

%%%%%%%%%%%%%%%%%%%%%%%%%%%%%%%%%%%%%%%%%%%%%%%%%%%%%%%%%%%
%
%                 Section 2
%
%%%%%%%%%%%%%%%%%%%%%%%%%%%%%%%%%%%%%%%%%%%%%%%%%%%%%%%%%%%
\section
{\bf 
Calabi-Yau varieties with at most one 
ordinary double point}
\par

\subsection
{}
{\bf Calabi-Yau varieties with at most one ODP and 
their deformations}

\subsubsection
{Calabi-Yau varieties with at most one ODP}

Recall that an $n$-dimensional singularity is an 
{\it ordinary double point} ({\it ODP} for short) 
if it is isomorphic to the hypersurface singularity 
at $0\in{\Bbb C}^{n}$ defined by the equation
$z_{0}^{2}+\cdots+z_{n}^{2}=0$.

\begin{defn}
A complex projective variety $X$ of dimension $n\geq3$ 
satisfying the following conditions is called a 
{\it Calabi-Yau $n$-fold with at most one ODP}: 
\newline
{(i) }
There exists a nowhere vanishing holomorphic $n$-form
on $X_{\rm reg}=X\setminus{\rm Sing}(X)$;
\newline
{(ii) }
$X$ is connected and
$H^{q}(X,{\cal O}_{X})=0$ for $0<q<n$;
\newline
{(iii) }
The singular locus ${\rm Sing}(X)$ consists of 
empty or at most one ODP.
\end{defn}

Throughout this paper, we use the following notation:
For a complex space $Y$,
let $\Theta_{Y}$ be the tangent sheaf of $Y$,
let $\Omega^{1}_{Y}$ be the sheaf of 
K\"ahler differentials on $Y$,
and let $K_{Y}$ be the dualizing sheaf of $Y$.
The sheaf $\Omega^{p}_{Y}$ is defined as
$\bigwedge^{p}\Omega^{1}_{Y}$. On the regular part
of $Y$, the sheaves $\Theta_{Y}$, $\Omega^{p}_{Y}$,
$K_{Y}$ are often identified with the corresponding
holomorphic vector bundles
$TY$, $\bigwedge^{p}T^{*}Y$, $\det T^{*}Y$,
respectively.
\par
We set $\varDelta(r):=\{t\in{\Bbb C};\,|t|<r\}$
and $\varDelta(r)^{*}:=\varDelta(r)\setminus\{0\}$
for $r>0$. We write $\varDelta$ (resp. $\varDelta^{*}$)
for $\varDelta(1)$ (resp. $\varDelta(1)^{*}$).
\par
Since an ODP is a hypersurface singularity, 
the dualizing sheaf of a Calabi-Yau $n$-fold 
with at most one ODP is trivial by (i).

\subsubsection
{Deformations of Calabi-Yau varieties with 
at most one ODP} 
\par
Let $X$ be a Calabi-Yau $n$-fold with 
at most one ODP. 

\begin{defn}
Let $(S,0)$ be a complex space with marked point
and let $\cal X$ be a complex space.
A proper, surjective, flat holomorphic map
$\pi\colon{\cal X}\to S$ is called 
a {\it deformation} of $X$
if $\pi^{-1}(0)\cong X$. 
If $\cal X$ and $S$ are smooth and if a general fiber
of $\pi\colon{\cal X}\to S$ is smooth,
the deformation $\pi\colon({\cal X},X)\to(S,0)$ is 
called a {\it smoothing} of $X$.
If there exists a smoothing of $X$, $X$ is said to be
{\it smoothable}. 
\end{defn}

We refer to \cite[Example\,5.8]{Namikawa94}
for an example of a non-smoothable Calabi-Yau 
threefold with a unique ODP as its singular set.
\par
Since $H^{0}(X,\Theta_{X})=0$ 
(cf. \cite[pp.432, l.23]{Namikawa94}),
there exists a deformation germ
${\frak p}\colon({\frak X},X)\to(\rm{Def}(X),[X])$ 
of $X$ with the universal property:
Every deformation germ $\pi\colon({\cal X},X)\to(S,0)$ 
is induced from ${\frak p}\colon{\frak X}\to\rm{Def}(X)$ 
by a unique holomorphic map 
$f\colon(S,0)\to(\rm{Def}(X),[X])$.
This local universal deformation of $X$ is
called the {\it Kuranishi family} of $X$.
The Kuranishi family is unique up to an isomorphism.
The base space $(\rm{Def}(X),[X])$
is called the {\it Kuranishi space} of $X$.
By \cite{Kawamata92}, \cite{Ran92}, \cite{Tian87}, 
\cite{Tian92}, \cite{Todorov89}, 
${\rm Def}(X)$ is smooth. 
We denote by $T_{{\rm Def}(X),[X]}$ the tangent
space of ${\rm Def}(X)$ at $[X]$. See 
\cite{Douady74}, \cite{Grauert74}, \cite{Kuranishi62} 
for more details about the Kuranishi family.
\par
For a deformation $\pi\colon({\cal X},X)\to(S,0)$, 
the fiber $X_{s}$ $(s\in S)$ is a Calabi-Yau $n$-fold 
with at most one ODP if $s\in S$ is sufficiently 
close to $0$ (cf. \cite[Prop.\,6.1]{Namikawa94},
\cite[Prop.\,4.2]{Tian92}).
\par
In the rest of Subsection 2.1,
we assume that $X$ is a smoothable Calabi-Yau $n$-fold
with at most one ODP.
Let $\pi\colon({\cal X},X)\to(S,0)$ be a smoothing. 
The critical locus of $\pi$ is defined by
$$
\Sigma_{\pi}:=\{x\in{\cal X};\,
d\pi_{x}=0\}.
$$
The discriminant locus of $\pi\colon{\cal X}\to S$
is the subvariety of $S$ defined by 
$$
{\cal D}
:=
\pi(\Sigma_{\pi})
=
\{s\in S;\,{\rm Sing}(X_{s})\not=\emptyset\}.
$$

\begin{lem}
Let $N+1=\dim S$.
For $p\in{\rm Sing}(X)$, there exists a neighborhood 
$V_{p}\cong\varDelta^{n+1}\times\varDelta^{N}$ of $p$ 
in $\cal X$ such that 
$$
\pi|_{V_{p}}(z,w)
=(z_{0}^{2}+\cdots+z_{n}^{2},w_{1},\ldots,w_{N}),
\qquad
z=(z_{0},\ldots,z_{n}),
\quad
w=(w_{1},\ldots,w_{N}).
$$
In particular, if ${\rm Sing}(X)\not=\emptyset$,
$\cal D$ is a divisor of $S$ smooth at $0$.
\end{lem}

\begin{pf}
Let $p\in{\rm Sing}(X)$. Let $s=(s_{0},\ldots,s_{N})$
be a system of coordinates near $0\in S$.
By e.g. \cite[pp.103, (6.7)]{Looijenga84}, 
there exists $f_{p}\in{\cal O}_{S,0}$ such that 
$$
{\cal O}_{{\cal X},p}\cong
{\cal O}_{{\Bbb C}^{n+1}\times S,(0,0)}/
(z_{0}^{2}+\cdots+z_{n}^{2}+f_{p}(s)),
\qquad
\pi(z,s)=s.
$$
Since $\cal X$ is smooth, we get $df_{p}(0)\not=0$. 
Hence we can assume that $f_{p}(s)=s_{0}$ after
a suitable change of the coordinates of $S$.
\end{pf}

\subsubsection
{The Kodaira-Spencer map}
\par
For a smoothing $\pi\colon({\cal X},X)\to(S,0)$,
the short exact sequence of sheaves on $X$
$$
0\longrightarrow
\pi^{*}\Omega^{1}_{S}|_{X}
\longrightarrow
\Omega^{1}_{\cal X}|_{X}
\longrightarrow
\Omega^{1}_{X}
\longrightarrow0
$$
induces the long exact sequence:
$$
\cdots\longrightarrow
\rm{Hom}_{{\cal O}_{X}}(\pi^{*}\Omega^{1}_{S}|_{X},
{\cal O}_{X})
\longrightarrow
{\rm Ext}^{1}_{{\cal O}_{X}}(\Omega^{1}_{X},{\cal O}_{X})
\longrightarrow
{\rm Ext}^{1}_{{\cal O}_{X}}(\Omega^{1}_{\cal X}|_{X},
{\cal O}_{X})
\longrightarrow\cdots
$$

\begin{defn}
The {\it Kodaira-Spencer map} of 
$\pi\colon({\cal X},X)\to(S,0)$ is the coboundary map
$$
\rho_{0}\colon T_{0}S=
\rm{Hom}_{{\cal O}_{X}}(\pi^{*}\Omega^{1}_{S}|_{X},
{\cal O}_{X})
\to
{\rm Ext}^{1}_{{\cal O}_{X}}(\Omega^{1}_{X},{\cal O}_{X}).
$$ 
\end{defn}

\begin{prop} 
The Kodaira-Spencer map 
$\rho_{[X]}\colon T_{{\rm Def}(X),[X]}\to
{\rm Ext}^{1}_{{\cal O}_{X}}(\Omega^{1}_{X},{\cal O}_{X})$
for the Kuranishi family of $X$
is an isomorphism.
\end{prop}

\begin{pf}
See \cite{Kawamata92}, \cite{Ran92}, \cite{Tian87}, 
\cite{Tian92}, \cite{Todorov89}. 
\end{pf}

Let 
$$
r\colon
{\rm Ext}^{1}_{{\cal O}_{X}}(\Omega^{1}_{X},{\cal O}_{X})
\ni\alpha\to\alpha|_{X_{\rm reg}}\in 
{\rm Ext}^{1}_{{\cal O}_{X}}(\Omega^{1}_{X_{\rm reg}},
{\cal O}_{X_{\rm reg}})=H^{1}(X_{\rm reg},\Theta_{X})
$$ 
be the restriction map. 
Since $n\geq3$,
$r\colon
{\rm Ext}^{1}_{{\cal O}_{X}}(\Omega^{1}_{X},{\cal O}_{X})
\to H^{1}(X_{\rm reg},\Theta_{X})$ is an isomorphism
by \cite[Th.\,2]{Schlessinger71} and 
\cite[Prop.\,1.1]{Tian92}.

\begin{lem}
Under the natural identification
$H^{0}(X_{\rm reg},\pi^{*}\Theta_{S}|_{X_{\rm reg}})
\cong T_{0}S$ via $\pi$,
the composition 
$r\circ\rho_{0}\colon T_{0}S\to 
H^{1}(X_{\rm reg},\Theta_{X})$ 
is the coboundary map of the long exact sequence 
of cohomologies associated with the short exact sequence 
of sheaves
\begin{equation}
0\longrightarrow
\Theta_{X_{\rm reg}}
\longrightarrow
\Theta_{\cal X}|_{X_{\rm reg}}
\longrightarrow
\pi^{*}\Theta_{S}|_{X_{\rm reg}}
\longrightarrow0.
\end{equation}
\end{lem}

\begin{pf}
The commutative diagram of the short exact sequences
of sheaves
\begin{equation*}
\begin{CD}
0
@>>> \pi^{*}\Omega^{1}_{S}|_{X} 
@>>> \Omega^{1}_{\cal X}|_{X}
@>>> \Omega^{1}_{X} 
@>>>0
\\
@. @V r VV @V r VV @V r VV @.
\\
0
@>>> \pi^{*}\Omega^{1}_{S}|_{X_{\rm reg}}
@>>>
\Omega^{1}_{\cal X}|_{X_{\rm reg}}
@>>> 
\Omega^{1}_{X}|_{X_{\rm reg}} 
@>>>0
\end{CD}
\end{equation*}
induces the commutative diagram of exact sequences
\begin{equation*}
\begin{CD}
\,
@>>>
{\rm Hom}_{{\cal O}_{X}}(\pi^{*}\Omega^{1}_{S}|_{X},
{\cal O}_{X})
@>>> 
{\rm Ext}^{1}_{{\cal O}_{X}}(\Omega^{1}_{X},{\cal O}_{X})
@>>>\,
\\
@. @V r VV @V r VV @.
\\
\,@>>>
{\rm Hom}_{{\cal O}_{X}}(\pi^{*}\Omega^{1}_{S}
|_{X_{\rm reg}},{\cal O}_{X_{\rm reg}})
@>>> 
{\rm Ext}^{1}_{{\cal O}_{X}}(\Omega^{1}_{X_{\rm reg}},
{\cal O}_{X_{\rm reg}})
@>>>\,
\end{CD}
\end{equation*}
where the first (resp. second) vertical arrow is 
isomorphic by the normality of ${\rm Sing}(X)$ 
(resp. \cite[Prop.\,1.2]{Tian92}).
Since $\pi^{*}\Omega^{1}_{S}|_{X_{\rm reg}}$,
$\Omega^{1}_{\cal X}|_{X_{\rm reg}}$,
$\Omega^{1}_{X_{\rm reg}}$ are locally free,
the second line is the long exact sequence of 
cohomologies associated with (2.1).
\end{pf}

\begin{lem}
Suppose $X$ is smoothable.
Then the Kuranishi family of $X$ is 
a smoothing of $X$.
\end{lem}

\begin{pf} 
Since the assertion is obvious when $X$ is smooth, 
we assume that $X$ has a unique ODP $p$.
Since $X$ is smoothable, a general fiber of 
the Kuranishi family of $X$ is smooth. 
We must prove the smoothness of the total space 
$\frak X$ of the Kuranishi family of $X$. 
Since ${\rm Sing}(X)=\{p\}$, it suffices to prove 
the smoothness of $\frak X$ at $p$.
\par
Let 
${\rm Def}(X,p)\cong({\Bbb C},0)$
be the Kuranishi space of the ODP $(X,p)$
(cf. \cite[Chap.\,6 C]{Looijenga84}).
The universal deformation of $X$ induces 
a holomorphic map of germs
$f\colon{\rm  Def}(X)\to{\rm Def}(X,p)$.
The existence of a smoothing 
of $X$ implies the surjectivity of the differential 
of $f$ at $[X]$. Hence $f$ may be regarded as
a part of a system of coordinates of ${\rm Def}(X)$ 
at $[X]$. Since
\begin{equation}
{\cal O}_{{\frak X},p}\cong
{\cal O}_{{\Bbb C}^{n+1}\times{\rm Def}(X),(0,[X])}/
(z_{0}^{2}+\cdots+z_{n}^{2}+f)
\end{equation}
by e.g. \cite[pp.103, (6.7)]{Looijenga84},
this implies the smoothness of $\frak X$ at $p$. 
\end{pf}

Let ${\frak p}\colon({\frak X},X)\to({\rm Def}(X),[X])$
be the Kuranishi family of $X$.

\begin{prop}
There exist a pointed projective variety $(B,0)$, 
a projective variety $\frak Z$, 
and a surjective flat holomorphic map
$f\colon{\frak Z}\to B$ such that the deformation germ
$f\colon({\frak Z},f^{-1}(0))\to(B,0)$ is isomorphic to
${\frak p}\colon({\frak X},X)\to({\rm Def}(X),[X])$.
In particular, the map 
${\frak p}\colon{\frak X}\to{\rm Def}(X)$ 
is projective.
\end{prop}

\begin{pf}
See \cite[pp.\,441, l.7-l.12]{Namikawa94}.
\end{pf}

\subsubsection
{The Serre duality for Calabi-Yau varieties 
with at most one ODP}
\par
Let 
$$
\langle\cdot,\cdot\rangle\colon
H^{n-1}(X,\Omega^{1}_{X}\otimes K_{X})
\times{\rm Ext}^{1}_{{\cal O}_{X}}
(\Omega^{1}_{X}\otimes K_{X},K_{X})
\to
H^{n}(X,K_{X})\cong\Bbb C
$$
be the Yoneda product. Since $X$ is compact, 
the Yoneda product is a perfect pairing by 
\cite[Th.\,4.1 and Th.\,4.2]{Banica76}. 
Hence we get by Proposition 2.5
$$
H^{n-1}(X,\Omega^{1}_{X}\otimes K_{X})
=
{\rm Ext}^{1}_{{\cal O}_{X}}
(\Omega^{1}_{X}\otimes K_{X},K_{X})^{\lor}
=(T_{{\rm Def}(X),[X]})^{\lor}
=\Omega^{1}_{{\rm Def}(X),[X]}.
$$
If $X$ is smooth, then 
${\rm Ext}^{1}_{{\cal O}_{X}}
(\Omega^{1}_{X}\otimes K_{X},K_{X})=
H^{1}(X,\Theta_{X})$ and the Yoneda product
is given by the ordinary Serre duality pairing 
\cite[Th.\,4.2]{Banica76}. 
\par
Let
$H^{n-1}_{c}(X_{\rm reg},
\Omega^{1}_{X}\otimes K_{X})$ be
the cohomology with compact support.

\begin{lem}
The natural map
$H^{n-1}_{c}(X_{\rm reg},\Omega^{1}_{X}\otimes K_{X})
\to H^{n-1}(X_{\rm reg},
\Omega^{1}_{X}\otimes K_{X})$
is an isomorphism. Under this isomorphism,
the Yoneda product $\langle\cdot,\cdot\rangle$
coincides with the Serre duality pairing 
on the regular part of $X$:
$$
H^{n-1}_{c}(X_{\rm reg},
\Omega^{1}_{X}\otimes K_{X})
\times H^{1}(X_{\rm reg},\Theta_{X})
\to
H^{n}_{c}(X_{\rm reg},\,K_{X})\cong\Bbb C.
$$
\end{lem}

\begin{pf}
Since
${\rm Ext}^{1}_{{\cal O}_{X}}
(\Omega^{1}_{X},{\cal O}_{X})=
{\rm Ext}^{1}_{{\cal O}_{X_{\rm reg}}}
(\Omega^{1}_{X_{\rm reg}},{\cal O}_{X_{\rm reg}})$
by \cite[Prop.\,1.1]{Tian92}, 
the Serre duality for open manifolds
\cite[Th.\,4.1 and Th.\,4.2]{Banica76} yields that
$$
H^{n-1}(X,\Omega^{1}_{X}\otimes K_{X})
=
{\rm Ext}^{1}_{{\cal O}_{X}}
(\Omega^{1}_{X},{\cal O}_{X})^{\lor}
=
{\rm Ext}^{1}_{{\cal O}_{X_{\rm reg}}}
(\Omega^{1}_{X_{\rm reg}},{\cal O}_{X_{\rm reg}})
=
H^{n-1}_{c}(X_{\rm reg},
\Omega^{1}_{X}\otimes K_{X})
$$
and that the Yoneda product pairing
$$
H^{n-1}_{c}(X_{\rm reg},
\Omega^{1}_{X}\otimes K_{X})
\times{\rm Ext}^{1}_{{\cal O}_{X_{\rm reg}}}
(\Omega^{1}_{X_{\rm reg}}\otimes K_{X},K_{X})
\to
H^{n}_{c}(X_{\rm reg},\,K_{X})
$$
is perfect. Since $X_{\rm reg}$ is smooth, 
${\rm Ext}^{1}_{{\cal O}_{X_{\rm reg}}}
(\Omega^{1}_{X_{\rm reg}}\otimes
K_{X},K_{X})=H^{1}(X_{\rm reg},\Theta_{X})$ 
and the Yoneda product pairing 
$\langle\cdot,\cdot\rangle$ coincides with
the Serre duality pairing.
\end{pf}

\subsection
{}
{\bf The locally-freeness of 
the direct image sheaves: the case $n=3$}
\par
Let $n\geq3$.
Let $X$ be a smoothable Calabi-Yau $n$-fold with 
at most one ODP.
Let $\pi\colon({\cal X},X)\to(S,0)$ be a smoothing 
of $X$. Set
$\Omega^{1}_{{\cal X}/S}:=
\Omega^{1}_{\cal X}/\pi^{*}\Omega^{1}_{S}$. 

\begin{lem}
The sheaf $\Omega^{1}_{{\cal X}/S}$ is 
a flat ${\cal O}_{S}$-module.
\end{lem}

\begin{pf}
Since 
$\Omega^{1}_{{\cal X}/S,\,x}\cong
{\cal O}_{{\cal X},\,x}^{\oplus n}$ for 
$x\in\cal X\setminus\Sigma_{\pi}$, it suffices to 
prove the assertion for $x\in\Sigma_{\pi}$.
Let 
$({\frak Y},o)\to({\rm Def}(A_{1}),0)\cong({\Bbb C},0)$ 
be the Kuranishi family of an ODP $o$. 
There exists a map
$f\colon(S,\pi(x))\to({\rm Def}(A_{1}),0)$ such that
$({\cal X},x)\to(S,\pi(x))$ is induced from
$({\frak Y},o)\to({\rm Def}(A_{1}),0)$ by $f$. 
Let 
$p\colon{\cal X}={\frak Y}\times_{\rm Def(A_{1})}S
\to{\frak Y}$ be the projection. Since
$\Omega^{1}_{{\cal X}/S,\,x}=
p^{*}\Omega^{1}_{{\frak Y}/{\rm Def}(A_{1})}$,
the assertion follows from the fact that
$\Omega^{1}_{{\frak Y}/{\rm Def}(A_{1}),o}$ is a flat
${\cal O}_{{\rm Def}(A_{1}),0}$-module (cf.
\cite[p.\,13, l.\,28--p.\,14, l.\,1]{Namikawa02}).
\end{pf}

Let us consider the case $S={\rm Def}(X)$. Let
${\frak p}\colon({\frak X},X)\to({\rm Def}(X),[X])$ 
be the Kuranishi family of $X$.

\begin{thm}
If $n=3$,
the function ${\rm Def}(X)\ni s\to
h^{q}(X_{s},\Omega^{1}_{X_{s}})\in\Bbb Z$ 
is constant for all $q\geq0$. 
In particular,  
$R^{q}{\frak p}_{*}\Omega^{1}_{{\frak X}/{\rm Def}(X)}$ 
is a locally free
${\cal O}_{{\rm Def}(X)}$-module on ${\rm Def}(X)$
for all $q\geq0$.
\end{thm}

The proof of this theorem is divided 
into the four lemmas below.

\begin{lem}
If $n\geq3$, the function ${\rm Def}(X)\ni s\to
h^{n-1}(X_{s},\Omega^{1}_{X_{s}})\in\Bbb Z$ 
is constant. In particular,  
$R^{n-1}{\frak p}_{*}\Omega^{1}_{{\frak X}/{\rm Def}(X)}$ 
is a locally free
${\cal O}_{{\rm Def}(X)}$-module on ${\rm Def}(X)$.
\end{lem}

\begin{pf}
Since $K_{X}\cong{\cal O}_{X}$, we have
$$
T_{{\rm Def}(X),[X]}
\cong
{\rm Ext}^{1}_{{\cal O}_{X}}
(\Omega^{1}_{X},{\cal O}_{X})
=
{\rm Ext}^{1}_{{\cal O}_{X}}
(K_{X}\otimes\Omega^{1}_{X},K_{X})
=
H^{n-1}(X,K_{X}\otimes\Omega^{1}_{X})^{\lor},
$$
where the first isomorphism follows from 
Proposition 2.5, the second equality follows from 
the triviality of $K_{X}$, 
and the third equality follows from 
the Serre duality 
\cite[Chap.\,III Th.\,7.6 (b) (iii)]{Hartshorne77}. 
Notice that we can apply the Serre duality to $X$, 
because $X$ has at most one ODP
and hence $X$ is Cohen-Macaulay 
\cite[Chap.\,II Th.\,8.21, Prop.\,8.23]{Hartshorne77}.
Since $K_{X}\cong{\cal O}_{X}$, we get
$h^{n-1}(X,\Omega^{1}_{X})=
\dim T_{{\rm Def}(X),[X]}$.
The smoothness of ${\rm Def}(X)$ at $[X]$ 
implies that the function on ${\rm Def}(X)$
$$
{\rm Def}(X)\ni s\to\dim T_{{\rm Def}(X),s}=
\dim T_{{\rm Def}(X_{s}),[X_{s}]}=
h^{n-1}(X_{s},\Omega^{1}_{X_{s}})\in\Bbb Z
$$ 
is constant, for the Zariski tangent space coincides 
with the usual tangent space for smooth varieties.
Notice that the first equality
$\dim T_{{\rm Def}(X),s}=
\dim T_{{\rm Def}(X_{s}),[X_{s}]}$
follows from \cite[Sect.\,8.2]{Douady74}.
Since $\Omega^{1}_{{\frak X}/{\rm Def}(X)}$ 
is a flat ${\cal O}_{{\rm Def}(X)}$-module by
Lemmas 2.7 and 2.10,
$R^{n-1}{\frak p}_{*}\Omega^{1}_{{\frak X}/{\rm Def}(X)}$ 
is locally free by 
\cite[Chap.\,3, Th.\,4.12 (ii)]{Banica76}.
\end{pf}

\begin{lem}
If $n=3$, then $h^{3}(X_{s},\Omega^{1}_{X_{s}})=0$
for all $s\in{\rm Def}(X)$. In particular, 
$R^{3}\pi_{*}\Omega^{1}_{{\frak X}/{\rm Def}(X)}=0$.
\end{lem}

\begin{pf}
See \cite[p.\,432, l.23]{Namikawa94}.
\end{pf}

\begin{lem}
If $n=3$, the function 
${\rm Def}(X)\ni s\to 
h^{1}(X_{s},\Omega^{1}_{X_{s}})\in\Bbb Z$
is constant. In particular,
$R^{1}{\frak p}_{*}\Omega^{1}_{{\frak X}/{\rm Def}(X)}$ 
is a locally free
${\cal O}_{{\rm Def}(X)}$-module. 
\end{lem}

\begin{pf}
Since $\Omega^{1}_{{\frak X}/{\rm Def}(X)}$ is
a flat ${\cal O}_{{\rm Def}(X)}$-module, the function
${\rm Def}(X)\ni s\to\chi(X_{s},\Omega^{1}_{X_{s}})
\in\Bbb Z$ is constant, where 
$\chi(X_{s},\Omega^{1}_{X_{s}})$ denotes
the Euler characteristic of $\Omega^{1}_{X_{s}}$. 
Since $h^{q}(X_{s},\Omega^{1}_{X_{s}})$ is 
independent of $s\in{\rm Def}(X)$ for all $q\not=1$ 
by Lemmas 2.12 and 2.13, this implies that
$h^{1}(X_{s},\Omega^{1}_{X_{s}})$ is independent of
$s\in{\rm Def}(X)$.
\end{pf}

\begin{lem}
If $n=3$, then $R^{1}{\frak p}_{*}\Omega^{1}_{\frak X}$
is locally free. 
Moreover, the restriction map 
$R^{1}{\frak p}_{*}\Omega^{1}_{\frak X}\to
R^{1}{\frak p}_{*}\Omega^{1}_{{\frak X}/{\rm Def}(X)}$ 
is an isomorphism of ${\cal O}_{{\rm Def}(X)}$-modules.
\end{lem}

\begin{pf}
Set $N:=\dim{\rm Def}(X)$.
The short exact sequence of sheaves on $\frak X$
$$
0\to
{\cal O}_{\frak X}^{\oplus N}\cong
{\frak p}^{*}\Omega^{1}_{{\rm Def}(X)}
\to
\Omega^{1}_{\frak X}\to
\Omega^{1}_{{\frak X}/{\rm Def}(X)}
\to0
$$
induces the long exact sequence of direct images
$$
\cdots\longrightarrow
R^{1}{\frak p}_{*}{\frak p}^{*}\Omega^{1}_{{\rm Def}(X)}
\longrightarrow
R^{1}{\frak p}_{*}\Omega^{1}_{\frak X}
\longrightarrow
R^{1}{\frak p}_{*}\Omega^{1}_{{\frak X}/{\rm Def}(X)}
\longrightarrow
R^{2}{\frak p}_{*}{\frak p}^{*}\Omega^{1}_{{\rm Def}(X)}
\longrightarrow\cdots
$$
Since
$R^{1}{\frak p}_{*}{\frak p}^{*}\Omega^{1}_{{\rm Def}(X)}
=(R^{1}{\frak p}_{*}{\cal O}_{\frak X})^{\oplus N}=0$ 
and
$R^{2}{\frak p}_{*}{\frak p}^{*}\Omega^{1}_{{\rm Def}(X)}
=(R^{2}{\frak p}_{*}{\cal O}_{\frak X})^{\oplus N}=0$
by Definition 2.1 (ii), the second assertion follows
from the above exact sequence. 
\par
By the same argument as above, we see that 
the restriction map 
$H^{1}(X_{s},\Omega^{1}_{\frak X}|_{X_{s}})\to
H^{1}(X_{s},\Omega^{1}_{X_{s}})$ is an isomorphism
for all $s\in{\rm Def}(X)$. Hence 
$h^{1}(X_{s},\Omega^{1}_{\frak X}|_{X_{s}})$
is independent of $s\in{\rm Def}(X)$ by Lemma 2.14.
This, together with 
\cite[Chap.\,3, Th.\,4.12 (ii)]{Banica76}
proves the first assertion.
\end{pf}

Theorem 2.11 follows from Lemmas 2.12, 2.13, 2.14,
and 2.15. $\Box$
\newline
\par
Let $H^{2}(X,{\Bbb Z})_{{\rm Def}(X)}$ be 
the constant sheaf on ${\rm Def}(X)$ with stalk 
$H^{2}(X,\Bbb Z)$. By \cite[Prop.\,6.1]{Namikawa94},
$R^{2}{\frak p}_{*}{\Bbb Z}$ is isomorphic to
the constant sheaf $H^{2}(X,{\Bbb Z})_{{\rm Def}(X)}$.
\par
Since $R^{1}{\frak p}_{*}{\cal O}_{\frak X}=
R^{2}{\frak p}_{*}{\cal O}_{\frak X}=0$ by
Definition 2.1 (ii), 
the exponential sequence on $\frak X$ induces 
the exact sequence of direct images
\begin{equation}
\begin{CD}
0=R^{1}{\frak p}_{*}{\cal O}_{\frak X}
\longrightarrow
R^{1}{\frak p}_{*}{\cal O}_{\frak X}^{*}
@>\cong>>
R^{2}{\frak p}_{*}{\Bbb Z}
\longrightarrow
R^{2}{\frak p}_{*}{\cal O}_{\frak X}=0.
\end{CD}
\end{equation}
\par
For a holomorphic line bundle
${\cal L}\in H^{1}({\frak X},{\cal O}_{\frak X}^{*})$,
the Dolbeault cohomology class of the Chern form
$c_{1}({\cal L},h)\in 
H^{1}({\frak X},\Omega^{1}_{\frak X})$ is independent
of the choice of a Hermitian metric $h$ on $\cal L$,
which we will denote by ${\frak C}_{1}({\cal L})$.
Since every element of $H^{2}({\frak X},{\Bbb Z})$
is represented uniquely as the Chern class of
an element of $H^{1}({\frak X},{\cal O}_{\frak X}^{*})$
by the isomorphism (2.3), we define the map
$j\colon H^{2}(X,{\Bbb Z})\to
H^{1}({\frak X},\Omega^{1}_{\frak X})$ by
$$
j\left(c_{1}({\cal L})|_{X}\right)
:={\frak C}_{1}({\cal L}),
\qquad
{\cal L}\in H^{1}({\frak X},{\cal O}_{\frak X}^{*}).
$$
We regard ${\frak C}_{1}({\cal L})$ as an element of 
$H^{0}({\rm Def}(X),
R^{1}{\frak p}_{*}\Omega^{1}_{{\frak X}/{\rm Def}(X)})$ 
after Lemma 2.15.
Since $H^{2}(X,{\Bbb Z})$ is finitely generated,
the map $j$ extends to a homomorphism of
${\cal O}_{{\rm Def}(X)}$-modules
$$
j\colon 
H^{2}(X,{\Bbb Z})_{{\rm Def}(X)}\otimes_{\Bbb Z}
{\cal O}_{{\rm Def}(X)}
\to
R^{1}{\frak p}_{*}\Omega^{1}_{{\frak X}/{\rm Def}(X)}.
$$

\begin{lem}
The homomorphism $j$ is an isomorphism of
${\cal O}_{{\rm Def}(X)}$-modules.
\end{lem}

\begin{pf}
Since
$H^{2}(X,{\Bbb Z})_{{\rm Def}(X)}\otimes_{\Bbb Z}
{\cal O}_{{\rm Def}(X)}$ and 
$R^{1}{\frak p}_{*}\Omega^{1}_{{\frak X}/{\rm Def}(X)}$ 
are locally free by Lemma 2.15, it suffices to prove 
that $j|_{X}\colon H^{2}(X,{\Bbb C})\to
H^{1}(X,\Omega^{1}_{X})$ is an isomorphism.
Since 
$h^{2}(X,{\Bbb C})=h^{2}(X_{s},{\Bbb C})$
by \cite[Prop.\,6.1]{Namikawa94} and since
$h^{1}(X_{s},\Omega^{1}_{X_{s}})=
h^{1}(X,\Omega^{1}_{X})$ by Lemma 2.14, we get
$h^{2}(X,{\Bbb C})=h^{1}(X,\Omega^{1}_{X})$.
Since
$j|_{X}$ is surjective by \cite[Lemma 2.2]{Namikawa94},  
it is an isomorphism.
\end{pf}

%%%%%%%%%%%%%%%%%%%%%%%%%%%%%%%%%%%%%%%%%%%%%%%%%%%%%%%%%%%
%
%                 Section 3
%
%%%%%%%%%%%%%%%%%%%%%%%%%%%%%%%%%%%%%%%%%%%%%%%%%%%%%%%%%%%

\section
{\bf Quillen metrics}
\par
Throughout Section 3, we fix the following notation:
Let $X$ be a complex manifold.
Let $(F,h_{F})$ be a holomorphic Hermitian vector bundle 
on $X$, which we also write $\overline{F}=(F,h_{F})$
for simplicity.

\subsection
{}
{\bf Analytic torsion and BCOV torsion}
\par
In Subsection\,3.1, assume that $X$ is 
a compact K\"ahler manifold with 
K\"ahler metric $g_{X}$ and 
with K\"ahler form $\gamma_{X}$.
We set $\overline{X}=(X,g_{X})$. 
Define $\overline{\Omega}_{X}^{p}$ to be
the holomorphic vector bundle $\Omega_{X}^{p}$ 
equipped with the Hermitian metric induced from $g_{X}$.
\par
Let $A^{p,q}_{X}(F)$ be the vector space of $F$-valued 
smooth $(p,q)$-forms on $X$. Set
$S_{F}=\bigoplus_{q\geq0}A^{0,q}_{X}(F)$.
Let $\langle\cdot,\cdot\rangle$ be the Hermitian metric 
on $(\bigwedge T^{*(0,1)}X)\otimes F$ induced from 
$g_{X}$ and $h_{F}$. 
The volume form of $\overline{X}$ is defined by 
$dv_{X}=\gamma_{X}^{\dim X}/(\dim X)!$. 
The {\it $L^{2}$-metric} is the
Hermitian metric on $S_{F}$ defined by 
$$
(s,s')_{L^{2}}:=
\frac{1}{(2\pi)^{\dim X}}
\int_{X}\langle s(x),s'(x)\rangle_{x}\,dv_{X}(x),
\qquad s,s'\in S_{F}.
$$
\par
Let $\bar{\partial}_{F}$ be the Dolbeault operator 
acting on $S_{F}$ and let $\bar{\partial}^{*}_{F}$ be
the formal adjoint of $\bar{\partial}_{F}$ with 
respect to $(\cdot,\cdot)_{L^{2}}$. Then
$\square_{F}=
(\bar{\partial}_{F}+\bar{\partial}^{*}_{F})^{2}$
is the corresponding $\bar{\partial}$-Laplacian. 
Let $\sigma(\square_{F})$ be the spectrum 
of $\square_{F}$ and let $E_{F}(\lambda)$ be 
the eigenspace of $\square_{F}$ 
with respect to the eigenvalue $\lambda$. 
\par
Let $N$ and $\epsilon$ be the operators on $S_{F}$
defined by $N=q$ and $\epsilon=(-1)^{q}$
on $A^{0,q}_{X}(F)$. Then $N$ and $\epsilon$ 
preserve $E_{F}(\lambda)$. 
\par
The zeta function
$$
\zeta_{\overline{F}}(s):=
\sum_{\lambda\in\sigma(\square_{F})\setminus\{0\}}
\lambda^{-s}\,{\rm Tr}\,[\epsilon N|_{E_{F}(\lambda)}].
$$
converges absolutely for $s\in{\Bbb C}$ with
${\rm Re}\,s\gg1$. 
By \cite[II, Th.\,2.16, (2.98)]{BGS88}, 
$\zeta_{\overline{F}}(s)$ has a meromorphic
continuation to the complex plane, which is
holomorphic at $s=0$. 

\begin{defn}
(i) 
The {\it analytic torsion} of 
$(\overline{X},\overline{F})$ is defined by
$$
\tau(\overline{X},\overline{F}):=
\exp(-\zeta_{\overline{F}}'(0)).
$$
\newline
{(ii) }
The {\it BCOV torsion} of $\overline{X}$ is defined by
$$
{\cal T}_{\rm BCOV}(\overline{X}):=
\prod_{p\geq0}
\tau(\overline{X},\overline{\Omega}^{p}_{X})^{(-1)^{p}p}
=\exp[
-\sum_{p\geq0}(-1)^{p}p\,
\zeta'_{\overline{\Omega}_{X}^{p}}(0)].
$$
\end{defn}

We refer the reader to
\cite{BGS88}, \cite{RaySinger73}
for more details about analytic torsion.

\subsection
{}
{\bf Quillen metrics}
\par

\begin{defn}
(i) 
The {\it determinant of the cohomologies of $F$}
is the complex line defined by
$$
\lambda(F):=
\bigotimes_{q\geq0}(\det H^{q}(X,F))^{(-1)^{q}}.
$$
\newline{(ii) }
The {\it BCOV line }\rm is the complex line 
$\lambda(\Omega^{\bullet}_{X})$ defined by
$$
\lambda(\Omega^{\bullet}_{X}):=
\bigotimes_{p\geq0}\lambda(\Omega^{p}_{X})^{(-1)^{p}p}
=\bigotimes_{p,q\geq0}
(\det H^{q}(X,\Omega^{p}_{X}))^{(-1)^{p+q}p}.
$$
\end{defn}

Set 
$K^{q}(\overline{X},\overline{F})=
\ker\square_{F}\cap A^{0,q}_{X}(F)$. Then
$K^{q}(\overline{X},\overline{F})$ inherits 
a metric from $(\cdot,\cdot)_{L^{2}}$. 
By Hodge theory, we have an isomorphism 
$H^{q}(X,F)\cong K^{q}(\overline{X},\overline{F})$. 
We define $h_{H^{q}(X,F)}$ to be
the metric on $H^{q}(X,F)$ induced from 
the $L^{2}$-metric
on $K^{q}(\overline{X},\overline{F})$ by 
this isomorphism. 
\par
Let $\|\cdot\|_{L^{2},\lambda(F)}$
be the Hermitian metric on $\lambda(F)$ 
induced from $\{h_{H^{q}(X,F)}\}_{q\geq0}$.

\begin{defn}
(i) 
The {\it Quillen metric } on $\lambda(F)$ is 
defined by
$$
\|\alpha\|^{2}_{Q,\lambda(F)}
:=
\tau(\overline{X},\overline{F})\cdot
\|\alpha\|_{L^{2},\lambda(F)}^{2},
\qquad
\alpha\in\lambda(F).
$$
\newline{(ii) }
The {\it Quillen metric } on 
$\lambda(\Omega_{X}^{\bullet})$ 
is defined by
$$
\|\cdot\|_{Q,\lambda(\Omega_{X}^{\bullet})}^{2}
:=
\bigotimes_{p\geq0}
\|\cdot\|_{Q,\lambda(\Omega_{X}^{p})}^{(-1)^{p}2p}
=
{\cal T}_{\rm BCOV}(\overline{X})
\cdot
\bigotimes_{p\geq0}
\|\cdot\|_{L^{2},\lambda(\Omega_{X}^{p})}^{(-1)^{p}2p}.
$$
\end{defn}

Let $(F_{1},h_{F_{1}}),\cdots,(F_{l},h_{F_{l}})$
be holomorphic Hermitian vector bundles on $X$,
and let $\|\cdot\|^{2}_{Q,\lambda(F_{k})}$ be
the Quillen metric on $\lambda(F_{k})$.
For $\otimes_{k=1}^{l}\alpha_{k}\in
\bigotimes_{k=1}^{l}\lambda(F_{k})$, we set
$\|\otimes_{k=1}^{l}\alpha_{k}
\|_{Q,\otimes_{k}\lambda(F_{k})}^{2}:=
\prod_{k=1}^{l}
\|\alpha_{k}\|_{Q,\lambda(F_{k})}^{2}$.
When the line $\lambda(F)$ is clear from
the context, we write $\|\cdot\|_{Q}$ for
$\|\cdot\|_{Q,\lambda(F)}$.
We refer the reader to
\cite{BGS88}, \cite{BismutLebeau91}, \cite{Quillen85},
\cite{Soule92}
for more details about Quillen metrics.

\subsection
{}{\bf The Serre duality}
\par
Let $n:=\dim X$.
By the Serre duality, the following pairing 
on the Dolbeault cohomology groups is perfect:
$$
H^{q}(X,\Omega_{X}^{p})\times H^{n-q}(X,\Omega_{X}^{n-p})
\ni(\alpha,\beta)\to
\left(\frac{\sqrt{-1}}{2\pi}\right)^{n}
\int_{X}\alpha\wedge\beta\in{\Bbb C}.
$$
\par
Let $\{\psi_{i}\}$ be an arbitrary basis of 
$H^{q}(X,\Omega_{X}^{p})$, 
and let $\{\psi_{i}^{\lor}\}$ the dual basis of 
$H^{n-q}(X,\Omega_{X}^{n-p})$ with respect to 
the Serre duality pairing. 
Then the element of
$\det H^{p}(X,\Omega_{X}^{q})
\otimes\det H^{n-p}(X,\Omega_{X}^{n-q})$ defined by
\begin{equation}
{\bf 1}_{(p,q),(n-p,n-q)}:=
\bigwedge_{i}\psi_{i}\otimes
\bigwedge_{i}\psi_{i}^{\lor}
\end{equation}
is independent of the choice of a basis $\{\psi_{i}\}$
and is called the {\it canonical element}.
Similarly, the following element of
$\lambda(\Omega_{X}^{p})\otimes
\lambda(\Omega_{X}^{n-p})^{(-1)^{n}}$ is also called
the canonical element:
$$
{\bf 1}_{p,n-p}
={\bf 1}_{p,n-p}(X)
:=
\bigotimes_{q=0}^{n}{\bf 1}_{(p,q),(n-p,n-q)}
\in
\lambda(\Omega_{X}^{p})\otimes
\lambda(\Omega_{X}^{n-p})^{(-1)^{n}}.
$$
Then ${\bf 1}_{(p,q),(n-p,n-q)}
={\bf 1}^{-1}_{(p,q),(n-p,n-q)}$ by (3.1).

\par
Let $1_{\Bbb C}$ be the trivial Hermitian structure
on $\Bbb C$, i.e., 
$1_{\Bbb C}(a)=|a|^{2}$ for $a\in{\Bbb C}$.

\begin{prop}
The following identity holds:
\begin{equation}
\|{\bf 1}_{p,n-p}\|_{L^{2}}=
\|{\bf 1}_{p,n-p}\|_{Q}=1.
\end{equation}
In particular, 
the canonical element ${\bf 1}_{p,n-p}$
induces the following canonical
isometries of the Hermitian lines:
\begin{align}
\left(\lambda(\Omega_{X}^{p})\otimes
\lambda(\Omega_{X}^{n-p})^{(-1)^{n}},
\|\cdot\|_{L^{2},\lambda(\Omega_{X}^{p})\otimes
\lambda(\Omega_{X}^{n-p})^{(-1)^{n}}}\right)
&\cong
({\Bbb C},1_{\Bbb C}),
\\
\left(\lambda(\Omega_{X}^{p})\otimes
\lambda(\Omega_{X}^{n-p})^{(-1)^{n}},
\|\cdot\|_{Q,\lambda(\Omega_{X}^{p})\otimes
\lambda(\Omega_{X}^{n-p})^{(-1)^{n}}}\right)
&\cong
({\Bbb C},1_{\Bbb C}).
\end{align}
\end{prop}

\begin{pf}
Let 
$\{\phi_{i}\}$ be a unitary basis of 
$H^{q}(X,\Omega_{X}^{p})$ with respect to 
the $L^{2}$-metric. The dual basis of
$\{\phi_{i}\}$ with respect to the Serre duality 
pairing is given by $\{\bar{*}\phi_{i}\}$, 
where $*\colon A^{p,q}_{X}\to A^{n-q,n-p}_{X}$
is the Hodge $*$-operator with respect to 
the metric $g_{X}$.
By setting $\psi_{i}=\phi_{i}$ in (3.1), 
we get the first equality
\begin{equation}
\|{\bf 1}_{(p,q),(n-p,n-q)}\|_{L^{2}}=1,
\end{equation}
which yields the isometry (3.3).
\par
Let $\zeta_{p,q}(s)$ be the spectral zeta function
of the $\bar{\partial}$-Laplacian acting 
on $A^{p,q}_{X}$. Since 
$\bar{*}^{-1}\square_{p,q}\bar{*}=\square_{n-p,n-q}$,
we have $\zeta_{p,q}(s)=\zeta_{n-p,n-q}(s)$, 
which yields that
\begin{equation}
\tau(\overline{X},\overline{\Omega}_{X}^{p})
=
\tau(\overline{X},
\overline{\Omega}_{X}^{n-p})^{(-1)^{n+1}}.
\end{equation}
The second isometry (3.4) follows from 
(3.3) and (3.6).
\end{pf}

For more details about the Serre duality for 
Quillen metrics, we refer to \cite[(9)]{GilletSoule91}.

\subsection
{}{\bf Characteristic classes}
\par
In Subsections 3.4 and 3.5, 
we do {\it not } assume that
$X$ is compact K\"ahler.

\subsubsection
{Chern forms}
\par
For a square matrix $A$, set 
${\rm Td}(A):=\det\left(\frac{A}{I-\exp(-A)}\right)$
and ${\rm ch}(A):={\rm Tr}[e^{A}]$.
Let $R(\overline{F})$ be the curvature of 
$\overline{F}=(F,h_{F})$ with respect to 
the holomorphic Hermitian connection. 
The real closed forms on $X$ defined by
$$
{\rm Td}(F,h_{F}):=
{\rm Td}\left(-\frac{1}{2\pi\sqrt{-1}}
R(\overline{F})\right),
\qquad\qquad
{\rm ch}(F,h_{F}):=
{\rm ch}\left(-\frac{1}{2\pi\sqrt{-1}}
R(\overline{F})\right)
$$
are called the {\it Todd form} and the 
{\it Chern character form} of $\overline{F}$, 
respectively. 
\par
Let $c_{i}(F,h_{F})$ be
the $i$-th Chern form of $(F,h_{F})$.

\subsubsection
{Bott-Chern classes}
\par
Let 
${\cal E}:0\to E_{0}\to E_{1}\to\cdots\to E_{m}\to0$ 
be an exact sequence of holomorphic vector bundles 
on $X$, equipped with Hermitian metrics 
$h_{i}$, $i=0,\ldots,m$. We set 
$\overline{\cal E}:=({\cal E},\{h_{i}\}_{i=0}^{m})$.
By \cite[I, Th.\,1.29 and Eqs. (0.5), (1.124)]{BGS88}, 
one has the Bott-Chern secondary class
$\widetilde{\rm ch}(\overline{\cal E})\in
\bigoplus_{p\geq0}A^{p,p}(X)/
{\rm Im}\,\partial+{\rm Im}\,\bar{\partial}$
associated to 
the Chern character and $\overline{\cal E}$
such that
$$
dd^{c}\,\widetilde{\rm ch}(\overline{\cal E})
=
\sum_{i=0}^{m}(-1)^{i+1}{\rm ch}(E_{i},h_{i}).
$$
\par
Consider the case where $m=1$ and $E_{0}=E_{1}=E$.
Let $h'$ and $h$ be Hermitian metrics of $E_{0}$
and $E_{1}$, respectively.
By \cite[I, Th.\,1.27]{BGS88} or
\cite[Sect.\,1.2.4]{GilletSoule90},
one has the Bott-Chern secondary class
$\widetilde{\rm ch}(E;\,h,h')\in
\bigoplus_{p\geq0}A^{p,p}(X)/
{\rm Im}\,\partial+{\rm Im}\,\bar{\partial}$
such that 
$$
dd^{c}\widetilde{\rm ch}(E;\,h,h')=
{\rm ch}(E,h)-{\rm ch}(E,h').
$$
When ${\rm rk}(E)=1$, we have the following explicit
formula by 
\cite[I, (1.2.5.1), (1.3.1.2)]{GilletSoule90}:
\begin{equation}
\widetilde{\rm ch}(E;\,h,h')=
\sum_{m=1}^{\infty}\frac{1}{m!}\sum_{a+b=m-1}
c_{1}(E,h)^{a}c_{1}(E,h')^{b}\,
\log\left(\frac{h'}{h}\right).
\end{equation}
Similarly, let
$\widetilde{\rm Td}(E;h,h')\in
\bigoplus_{p\geq0}A^{p,p}(X)/
{\rm Im}\,\partial+{\rm Im}\,\bar{\partial}$
denote the Bott-Chern secondary class associated 
to the Todd form such that
$$
dd^{c}\,\widetilde{\rm Td}(E;h,h')
=
{\rm Td}(E,h)-{\rm Td}(E,h').
$$
For more details about Bott-Chern classes, we refer to
\cite{BGS88}, \cite{GilletSoule90}, \cite{Soule92}.

\subsection
{}{\bf The curvature formulas}
\par
Let $S$ be a complex manifold and
let $\pi\colon X\to S$ be a proper surjective 
holomorphic submersion. Then every fiber of $\pi$ is 
a compact complex manifold.
The map $\pi\colon X\to S$ is said to be 
{\it locally K\"ahler} if for every $s\in S$ 
there is an open subset $U\ni s$ such that 
$\pi^{-1}(U)$ possesses a K\"ahler metric.
We set $X_{s}=\pi^{-1}(s)$ for $s\in S$.
\par
Let $TX/S:=\ker\pi_{*}\subset TX$ be the relative 
holomorphic tangent bundle of the family 
$\pi\colon X\to S$. 
Set $\Omega^{p}_{X/S}:=\bigwedge^{p}(TX/S)^{\lor}$
and $K_{X/S}:=K_{X}\otimes(\pi^{*}K_{S})^{-1}
=\Omega^{\dim X-\dim S}_{X/S}$.
\par
A $C^{\infty}$ Hermitian metric on $TX/S$ is 
said to be {\it fiberwise K\"ahler}
if the induced metric on $X_{s}$ is K\"ahler 
for all $s\in S$. 
By Kodaira-Spencer, 
there exists an fiberwise K\"ahler metric 
on $TX/S$ if and only if 
every $X_{s}$ possesses a K\"ahler metric.
\par
Assume that every fiber $X_{s}$ possesses 
a K\"ahler metric.
Let $g_{X/S}$ be a fiberwise K\"ahler metric on $TX/S$. 
Set $g_{s}=g_{X/S}|_{X_{s}}$ and
$\overline{X}_{s}=(X_{s},g_{s})$ for $s\in S$.
We define $\overline{\Omega}^{p}_{X_{s}}$ to be
the holomorphic vector bundle $\Omega_{X_{s}}^{p}$ 
equipped with the Hermitian metric induced from $g_{s}$. 
When $p=0$, $\overline{\Omega}^{0}_{X_{s}}$ is defined
as the trivial line bundle ${\cal O}_{X_{s}}$ 
equipped with the trivial Hermitian metric.
\par
Since $\dim H^{q}(X_{s},\Omega_{X_{s}}^{p})$ 
is locally constant, the direct image sheaf 
$R^{q}\pi_{*}\Omega^{p}_{X/S}$ is locally free 
for all $p,q\geq0$ and is identified with the 
corresponding holomorphic vector bundle over $S$. 
Set
$$
\lambda(\Omega_{X/S}^{\bullet}):=
\bigotimes_{p,q\geq0}
(\det R^{q}\pi_{*}\Omega^{p}_{X/S})^{(-1)^{p+q}p}.
$$
Via the natural fiberwise identification
$\lambda(\Omega_{X/S}^{\bullet})|_{s}=
\lambda(\Omega_{X_{s}}^{\bullet})$ for all $s\in S$,
$\lambda(\Omega_{X/S}^{\bullet})$ is equipped with 
the Hermitian metric 
$\|\cdot\|_{\lambda(\Omega_{X/S}^{\bullet}),Q}$ 
defined by
$$
\|\cdot\|_{Q,\lambda(\Omega_{X/S}^{\bullet})}(s):=
\|\cdot\|_{Q,\lambda(\Omega_{X_{s}}^{\bullet})},
\qquad
s\in S,
$$
which is smooth by \cite[III, Cor.\,3.9]{BGS88}.
We set
$\lambda(\Omega_{X/S}^{\bullet})_{Q}:=
(\lambda(\Omega_{X/S}^{\bullet}),
\|\cdot\|_{Q,\lambda(\Omega_{X/S}^{\bullet})})$.
\par
Since 
$\dim K^{q}(\overline{X}_{s},
\overline{\Omega}^{p}_{X_{s}})$ is locally constant, 
there exists a $C^{\infty}$ vector bundle 
${\cal K}^{p,q}(X/S)$ over $S$ such that
${\cal K}^{p,q}(X/S)_{s}=K^{q}(\overline{X}_{s},
\overline{\Omega}^{p}_{X_{s}})$ for all $s\in S$. 
Then the fiberwise isomorphism
$H^{q}(X_{s},\Omega^{p}_{X_{s}})\cong 
K^{q}(\overline{X}_{s},\overline{\Omega}^{p}_{X_{s}})$ 
via Hodge theory induces an isomorphism of 
$C^{\infty}$ vector bundles
$R^{q}\pi_{*}\Omega^{p}_{X/S}\cong{\cal K}^{p,q}(X/S)$.
Let $h_{R^{q}\pi_{*}\Omega^{p}_{X/S}}$ be 
the $C^{\infty}$ Hermitian metric on 
$R^{q}\pi_{*}\Omega^{p}_{X/S}$ induced from 
the $L^{2}$-metric on ${\cal K}^{p,q}(X/S)$ 
by this isomorphism.
We define 
$\overline{R^{q}\pi_{*}\Omega^{p}_{X/S}}:=
(R^{q}\pi_{*}\Omega^{p}_{X/S},
h_{R^{q}\pi_{*}\Omega^{p}_{X/S}})$.
\par
Let ${\cal T}_{\rm BCOV}(X/S)$ be the function on $S$ 
defined by
$$
{\cal T}_{\rm BCOV}(X/S)(s):=
{\cal T}_{\rm BCOV}(\overline{X}_{s})=
\prod_{p\geq0}\tau(\overline{X}_{s},
\overline{\Omega}^{p}_{X_{s}})^{(-1)^{p}p},
\qquad
s\in S.
$$
\par 
For a differential form $\varphi$, $[\varphi]^{(p,q)}$
denotes the component of bidegree $(p,q)$ of $\varphi$.

\begin{thm}
Assume that the map $\pi\colon X\to S$ is 
locally K\"ahler and set $n=\dim X-\dim S$. 
Then ${\cal T}_{\rm BCOV}(X/S)$ 
lies in $C^{\infty}(S)$, and the following equation 
of $C^{\infty}$ $(1,1)$-forms on $S$ holds:
$$
\begin{aligned}
c_{1}(\lambda(\Omega_{X/S}^{\bullet})_{Q})
&=
-dd^{c}\log{\cal T}_{\rm BCOV}(X/S)
+\sum_{q\geq0}(-1)^{p+q}p\,
c_{1}(\overline{R^{q}\pi_{*}\Omega^{p}_{X/S}})\\
&=
-\frac{1}{12}\pi_{*}\left[c_{1}(TX/S,g_{X/S})\,
c_{n}(TX/S,g_{X/S})\right]^{(1,1)}.
\end{aligned}
$$
\end{thm}

\begin{pf}
See \cite[pp.\,374]{BCOV94} and \cite[Th.\,0.1]{BGS88}.
\end{pf}

%%%%%%%%%%%%%%%%%%%%%%%%%%%%%%%%%%%%%%%%%%%%%%%%%%%%%%%%%%%
%
%                 Section 4
%
%%%%%%%%%%%%%%%%%%%%%%%%%%%%%%%%%%%%%%%%%%%%%%%%%%%%%%%%%%%

\section
{\bf The BCOV invariant of Calabi-Yau manifolds}
\par
Throughout Section 4, we fix the following notation:
Let $X$ be a {\it smooth} Calabi-Yau $n$-fold. Let 
${\frak p}\colon({\frak X},X)\to({\rm Def}(X),[X])$
be the Kuranishi family of $X$.
\par
Let $g$ be a K\"ahler metric on $X$ with 
K\"ahler form $\gamma$. We define 
${\rm Vol}(X,\gamma):=
(2\pi)^{-n}\int_{X}\gamma^{n}/n!=\|1\|_{L^{2}}^{2}$.
Notice that our definition of ${\rm Vol}(X,\gamma)$
is different from the ordinary one because of
the factor $(2\pi)^{-n}$.
We set $c_{i}(X,\gamma):=c_{i}(TX,g)$ and
$\chi(X):=\int_{X}c_{n}(X,\gamma)$.
Let $\eta\in H^{0}(X,\Omega^{n}_{X})\setminus\{0\}$.

\subsection
{}{\bf The BCOV Hermitian line}
\par
Recall that the $L^{2}$-norm on 
$H^{0}(X,\Omega^{n}_{X})$ is independent of 
the choice of a K\"ahler metric $g$ because
$$
\|\eta\|_{L^{2}}^{2}=
(2\pi)^{-n}(\sqrt{-1})^{n^{2}}
\int_{X}\eta\wedge\bar{\eta}.
$$
After \cite[Sect.\,5.1]{Yoshikawa04}, 
we make the following:

\begin{defn}
(i) For $\overline{X}=(X,\gamma)$, define
${\cal A}(\overline{X})={\cal A}(X,\gamma)\in{\Bbb R}$
by
$$
{\cal A}(\overline{X}):=
{\rm Vol}(X,\gamma)^{\frac{\chi(X)}{12}}\,
\exp\left[-\frac{1}{12}\int_{X}
\log\left(
\frac{(\sqrt{-1})^{n^{2}}\eta\wedge\bar{\eta}}
{\gamma^{n}/n!}
\cdot
\frac{{\rm Vol}(X,\gamma)}{\|\eta\|_{L^{2}}^{2}}
\right)\,c_{n}(X,\gamma)\right].
$$
(ii)
The {\it BCOV metric} is the Hermitian structure
$\|\cdot\|_{\lambda(\Omega^{\bullet}_{X})}$ on
$\lambda(\Omega^{\bullet}_{X})$ defined by
$$
\|\cdot\|_{\lambda(\Omega^{\bullet}_{X})}^{2}:=
{\cal A}(\overline{X})\cdot
\|\cdot\|_{Q,\lambda(\Omega^{\bullet}_{X})}^{2}.
$$
\newline{(iii) }
The {\it BCOV Hermitian line} is defined by
$$
\overline{\lambda(\Omega^{\bullet}_{X})}:=
(\lambda(\Omega^{\bullet}_{X}),
\|\cdot\|_{\lambda(\Omega^{\bullet}_{X})}).
$$
\end{defn}

\begin{rem}
By Yau \cite{Yau78}, 
every K\"ahler class on $X$ contains 
a unique Ricci-flat  K\"ahler form. 
If $\kappa$ is a Ricci-flat K\"ahler form on $X$, 
then
$$
\frac{\kappa^{n}/n!}
{(\sqrt{-1})^{n^{2}}\eta\wedge\bar{\eta}}
=
\frac{{\rm Vol}(X,\kappa)}{\|\eta\|_{L^{2}}^{2}},
$$
and hence
$\log{\cal A}(X,\kappa)=
\frac{\chi(X)}{12}\,\log{\rm Vol}(X,\kappa)$
in this case.
\end{rem}

\subsection
{}{\bf The Weil-Petersson metric and the Hodge metric}
\par
To compute the curvature of the BCOV Hermitian line
bundles, let us recall the definitions of 
the Weil-Petersson metric \cite{Tian87}
and the Hodge metric \cite{Lu99}, \cite{Lu01}.
\par
By Proposition 2.5, the homomorphism of
${\cal O}_{{\rm Def}(X)}$-modules on 
${\rm Def}(X)$ induced by the Kodaira-Spencer map
$$
\rho_{{\rm Def}(X)}\colon
\Theta_{{\rm Def}(X)}\to
R^{1}{\frak p}_{*}\Theta_{{\frak X}/{\rm Def}(X)}
$$
is an isomorphism, which is called 
the {\it Kodaira-Spencer isomorphism } in this paper.
\par
Since 
$H^{n-1}(X_{s},\Omega^{1}_{X_{s}})\subset 
H^{n}(X_{s},\Bbb C)$ 
consists of primitive cohomology classes for 
all $s\in{\rm Def}(X)$, 
the $L^{2}$-metric on 
$R^{1}{\frak p}_{*}\Omega^{n-1}_{{\frak X}/{\rm Def}(X)}$ 
is independent
of the choice of a fiberwise-K\"ahler metric on 
$T{\frak X}/{\rm Def}(X)$ by e.g.
\cite[Th.\,6.32]{Voisin02}. 
We will often denote the $L^{2}$-metric 
$h_{R^{1}{\frak p}_{*}
\Omega^{n-1}_{{\frak X}/{\rm Def}(X)}}$ on
$R^{1}{\frak p}_{*}\Omega^{n-1}_{{\frak X}/{\rm Def}(X)}$
by $(\cdot,\cdot)_{L^{2}}$. 
Then
$$
(\xi,\zeta)_{L^{2}}
=
-(2\pi)^{-n}(\sqrt{-1})^{n^{2}}
\int_{X}\xi\wedge\overline{\zeta},
\qquad
\xi,\zeta\in H^{1}(X,\Omega^{n-1}_{X}).
$$
\par
For $s\in{\rm Def}(X)$,
let 
$\rho_{s}\colon T_{{\rm Def}(X),s}\to 
H^{1}(X_{s},\Theta_{X_{s}})$
be the Kodaira-Spencer map, and
let $\eta_{s}\in H^{0}(X_{s},\Omega^{n}_{X_{s}})
\setminus\{0\}$.
Let $\iota(\cdot)$ be the interior product.

\begin{defn}
The {\it Weil-Petersson metric } $g_{\rm WP}$ 
on ${\rm Def}(X)$
is defined by
$$
g_{\rm WP}(u,v):=
-\frac{\int_{X_{s}}
\iota(\rho_{s}(u))\eta_{s}\wedge
\overline{\iota(\rho_{s}(v))\eta_{s}}}
{\int_{X_{s}}\eta_{s}\wedge\overline{\eta}_{s}}
=
\frac{(\iota(\rho_{s}(u))\eta_{s},
\iota(\rho_{s}(v))\eta_{s})_{L^{2}}}
{\|\eta_{s}\|_{L^{2}}^{2}}
$$
for $u,v\in T_{{\rm Def}(X),s}$.
Let $\omega_{\rm WP}$ be the K\"ahler form 
of $g_{\rm WP}$. 
\end{defn}

Let $\eta_{{\frak X}/{\rm Def}(X)}$ be a local
basis of ${\frak p}_{*}K_{{\frak X}/{\rm Def}(X)}$.
By e.g. \cite[Th.\,2]{Tian87}, we have
\begin{equation}
\omega_{\rm WP}
=
-dd^{c}\log\|\eta_{{\frak X}/{\rm Def}(X)}\|_{L^{2}}^{2}
=
c_{1}({\frak p}_{*}K_{{\frak X}/{\rm Def}(X)},
\|\cdot\|_{L^{2}}).
\end{equation}

\begin{prop}
The Kodaira-Spencer map
$\rho_{{\rm Def}(X)}$ induces an isometry of 
the following holomorphic Hermitian vector bundles 
on ${\rm Def}(X)$:
$$
(\Theta_{{\rm Def}(X)},g_{\rm WP})\otimes
({\frak p}_{*}K_{{\frak X}/{\rm Def}(X)},
\|\cdot\|_{L^{2}})
\cong
(R^{1}{\frak p}_{*}\Omega^{n-1}_{{\frak X}/{\rm Def}(X)},
h_{R^{1}{\frak p}_{*}
\Omega^{n-1}_{{\frak X}/{\rm Def}(X)}}).
$$
In particular, $\rho_{{\rm Def}(X)}$ induces an isometry 
of the following holomorphic Hermitian line bundles 
on ${\rm Def}(X)$:
$$
\begin{aligned}
\,&
(\det R^{1}{\frak p}_{*}
\Omega^{n-1}_{{\frak X}/{\rm Def}(X)},
\det 
h_{R^{1}{\frak p}_{*}
\Omega^{n-1}_{{\frak X}/{\rm Def}(X)}})
\\
&
\cong
(\det\Theta_{{\rm Def}(X)},\det g_{\rm WP})\otimes
({\frak p}_{*}K_{{\frak X}/{\rm Def}(X)},
\|\cdot\|_{L^{2}})^{\otimes h^{1,n-1}(X)}.
\end{aligned}
$$
\end{prop}

\begin{pf}
The Kodaira-Spencer isomorphism is given by
$$
\Theta_{{\rm Def}(X)}\otimes
{\frak p}_{*}K_{{\frak X}/{\rm Def}(X)}
\ni u\otimes\eta
\to
\iota(\rho_{{\rm Def}(X)}(u))\eta\in
R^{1}{\frak p}_{*}\Omega^{n-1}_{{\frak X}/{\rm Def}(X)}.
$$
Hence
$(\iota(\rho_{{\rm Def}(X)}(u))\,\eta,
\iota(\rho_{{\rm Def}(X)}(v))\,\eta)_{L^{2}}=
g_{\rm WP}(u,v)\cdot\|\eta\|_{L^{2}}^{2}$
by Definition 4.3.
\end{pf}

\begin{defn}
The Ricci form of the Weil-Petersson metric is
the Chern form of the Hermitian line bundle
$(\det\Theta_{{\rm Def}(X)},\det g_{\rm WP})$:
$$
{\rm Ric}\,\omega_{\rm WP}:=
c_{1}(\det\Theta_{{\rm Def}(X)},\det g_{\rm WP}).
$$
\end{defn}

\begin{prop}
The following identities hold:
$$
c_{1}(\det R^{n-p}\pi_{*}
\Omega^{p}_{{\frak X}/{\rm Def}(X)},\|\cdot\|_{L^{2}})
=
\left\{
\begin{array}{ll}
-\omega_{\rm WP}& (p=0)
\\
-{\rm Ric}\,\omega_{\rm WP}-
h^{1,n-1}(X)\,\omega_{\rm WP}
&(p=1)
\\
{\rm Ric}\,\omega_{\rm WP}+
h^{1,n-1}(X)\,\omega_{\rm WP}
&(p=n-1)
\\
\omega_{\rm WP}
&(p=n).
\end{array}
\right.
$$
\end{prop}

\begin{pf}
The assertion for $p=0,n$ follows from (4.1).
The assertion for $p=1,n-1$ follows from 
Proposition 4.4 and the Serre duality.
\end{pf}

See \cite[Sect.\,2]{FangLu03} for a generalization
of Proposition 4.6.
In the case $n=3$, the following positivity result 
for ${\rm Ric}\,\omega_{\rm WP}+
(h^{1,2}(X)+3)\,\omega_{\rm WP}$ 
shall be crucial in Sect.\,7.

\begin{prop}
When $n=3$, the $(1,1)$-form 
${\rm Ric}\,\omega_{\rm WP}+
(h^{1,2}(X)+3)\,\omega_{\rm WP}$
is a K\"ahler form on ${\rm Def}(X)$.
\end{prop}

\begin{pf}
See \cite[Th.\,1.1]{Lu01}.
\end{pf}

\begin{defn}
When $n=3$, the {\it Hodge form} on ${\rm Def}(X)$ 
is the positive $(1,1)$-form on ${\rm Def}(X)$ 
defined as
$$
\omega_{\rm H}
:=
{\rm Ric}\,\omega_{\rm WP}+
(h^{1,2}(X)+3)\,\omega_{\rm WP}.
$$
The corresponding K\"ahler metric on the Kuranishi 
space ${\rm Def}(X)$ is called the {\it Hodge metric} 
on ${\rm Def}(X)$.
\end{defn}

The Hodge metric is related to the invariant
Hermitian metric on the period domain for
Calabi-Yau threefolds as follows.
Let $X$ be a polarized smooth Calabi-Yau threefold.
Let $D$ be the classifying space for
the polarized Hodge structures of weight $3$
on $H^{3}(X,{\Bbb Z})/{\rm Torsion}$
defined by Griffiths e.g. 
\cite[Sect.\,2]{Griffiths84}. 
Let $F^{i}$ $(i=1,2,3)$ be 
the Hodge bundles on $D$. Let $\omega_{D}$ be 
the invariant Hermitian metric of $D$. Let 
$f\colon{\rm Def}(X)\to D$
be the period map. Then we have
\newline{(a)}
$\omega_{\rm WP}=f^{*}(c_{1}(F^{3},\|\cdot\|_{L^{2}}))$
\cite{Voisin02};
\newline{(b)}
Up to a constant, $\omega_{\rm H}=f^{*}(\omega_{D})$
\cite{Lu99}. In particular, $\omega_{\rm H}$ is 
always K\"ahlerian.
\par
We refer to e.g. \cite{Griffiths84} for more 
details about the classifying space $D$.

\subsection
{}{\bf The curvature formula for 
the BCOV Hermitian line bundles}
\par
Let $\pi\colon({\cal X},X)\to(S,0)$ be a flat deformation 
of $X$. Set $X_{s}=\pi^{-1}(s)$ for $s\in S$.
Let $g_{{\cal X}/S}$ be a fiberwise-K\"ahler metric 
on $T{\cal X}/S$. Then the line bundle
$\lambda(\Omega^{\bullet}_{{\cal X}/S})$ on $S$
is equipped with the BCOV metric
$\|\cdot\|_{\lambda(\Omega^{\bullet}_{{\cal X}/S})}$
with respect to $g_{{\cal X}/S}$.
\par
Let $\mu\colon(S,0)\to({\rm Def}(X),[X])$ be the 
holomorphic map such that the family
$\pi\colon({\cal X},X)\to(S,0)$ is induced from 
the Kuranishi family by $\mu$. Then we have
$$
c_{1}(\pi_{*}\omega_{{\cal X}/S},\|\cdot\|_{L^{2}})
=
\mu^{*}\omega_{\rm WP}
$$
near $s=0$. 
Let $\eta_{{\cal X}/S}$ be a local basis of 
$\pi_{*}\omega_{{\cal X}/S}$ 
and set
$$
\omega_{{\rm WP},\,{\cal X}/S}
:=
\mu^{*}\omega_{\rm WP}
=
-dd^{c}\log\|\eta_{{\cal X}/S}\|_{L^{2}}^{2}
=
c_{1}(\pi_{*}\omega_{{\cal X}/S},\|\cdot\|_{L^{2}}).
$$

\begin{thm}
The following identity of $(1,1)$-forms on $(S,0)$ holds:
$$
c_{1}(\overline{\lambda(\Omega^{\bullet}_{{\cal X}/S})})
=
\frac{\chi(X)}{12}\,\omega_{{\rm WP},\,{\cal X}/S}.
$$
\end{thm}

\begin{pf}
We follow \cite[Sect.\,5.2]{Yoshikawa04}.
Since the assertion is of local nature,
it suffices to prove it when
$S\cong\varDelta^{\dim S}$. Then
$\pi_{*}K_{{\cal X}/S}\cong{\cal O}_{S}$.
Let
$\eta_{{\cal X}/S}\in H^{0}(S,
\pi_{*}K_{{\cal X}/S})$
be a nowhere vanishing holomorphic section.
For $s\in S$,
set $\eta_{s}=\eta_{{\cal X}/S}|_{X_{s}}$.
Then $\eta_{s}\in 
H^{0}(X_{s},K_{X_{s}})\setminus\{0\}$
and
$\eta_{{\cal X}/S}$ are identified with 
the family of holomorphic $n$-forms 
$\{\eta_{s}\}_{s\in S}$ 
varying holomorphically in $s\in S$.
Define $\|\eta_{{\cal X}/S}\|_{L^{2}}^{2}\in
C^{\infty}(S)$ by
$$
\|\eta_{{\cal X}/S}\|_{L^{2}}^{2}(s)=
\|\eta_{s}\|_{L^{2}}^{2},
\qquad
s\in S.
$$
\par
Set $g_{s}=g_{{\cal X}/S}|_{X_{s}}$. 
Then $g_{{\cal X}/S}$ is identified with 
the family of K\"ahler metrics $\{g_{s}\}_{s\in S}$. 
Let $\gamma_{s}$ be the K\"ahler form of $h_{s}$.
Let $\gamma_{{\cal X}/S}=\{\gamma_{s}\}_{s\in S}$ 
be the family of K\"ahler forms associated to 
$g_{{\cal X}/S}$.
\par
Define the $C^{\infty}$ functions
${\rm Vol}({\cal X}/S)$ and ${\cal A}({\cal X}/S)$
on $S$ by
$$
{\rm Vol}({\cal X}/S)(s)={\rm Vol}(X_{s},\gamma_{s}),
\qquad
A({\cal X}/S)(s)=A(X_{s},\gamma_{s}),
\qquad
s\in S.
$$
\par
Let $c_{i}({\cal X}/S)$ be the $i$-th Chern form 
of the holomorphic Hermitian vector bundle 
$(T{\cal X}/S,g_{{\cal X}/S})$.
Since
$$
c_{1}({\cal X}/S)=
-c_{1}(K_{{\cal X}/S},\det g_{{\cal X}/S}^{-1})=
dd^{c}\log\left(
\frac{(\sqrt{-1})^{n^{2}}\eta_{{\cal X}/S}\wedge
\overline{\eta}_{{\cal X}/S}}
{\gamma_{{\cal X}/S}^{n}/n!}\right),
$$
the following identity of $(1,1)$-forms 
on $\cal X$ holds:
\begin{equation}
\begin{aligned}
c_{1}({\cal X}/S)
&=
-\pi^{*}\left\{\omega_{{\rm WP},{\cal X}/S}+
dd^{c}\log{\rm Vol}({\cal X}/S)\right\}
\\
&\quad
+dd^{c}\log\left\{
\frac{(\sqrt{-1})^{n^{2}}\eta_{{\cal X}/S}\wedge
\bar{\eta}_{{\cal X}/S}}
{\gamma_{{\cal X}/S}^{n}/n!}\cdot
\pi^{*}\left(\frac{{\rm Vol}({\cal X}/S)}
{\|\eta_{{\cal X}/S}\|_{L^{2}}^{2}}
\right)\right\}.
\end{aligned}
\end{equation}
Then we get 
\begin{equation}
\begin{aligned}
\,&
-\frac{1}{12}\pi_{*}
\left[c_{1}({\cal X}/S)\,c_{n}({\cal X}/S)\right]
\\
&=
-\frac{1}{12}\pi_{*}\left[
-\pi^{*}\left\{\omega_{{\rm WP},{\cal X}/S}+
dd^{c}\log{\rm Vol}({\cal X}/S)\right\}\,
c_{n}({\cal X}/S)\right]
\\
&\quad+
\pi_{*}\left[-\frac{1}{12}
dd^{c}\log\left\{
\frac{(\sqrt{-1})^{n^{2}}\eta_{{\cal X}/S}\wedge
\bar{\eta}_{{\cal X}/S}}
{\gamma_{{\cal X}/S}^{n}/n!}\cdot
\pi^{*}\left(\frac{{\rm Vol}({\cal X}/S)}
{\|\eta_{{\cal X}/S}\|_{L^{2}}^{2}}
\right)\right\}\,
c_{n}({\cal X}/S)\right]
\\
&=
\frac{\chi(X)}{12}\,\omega_{{\rm WP},\,{\cal X}/S}
+dd^{c}\log A({\cal X}/S),
\end{aligned}
\end{equation}
where the first equality follows from (4.2), 
and the second one follows from the projection formula 
and the commutativity of $dd^{c}$ and $\pi_{*}$.
\par
Since the map $\pi\colon{\cal X}\to S$ is locally 
projective by Proposition 2.8, we may apply 
Theorem 3.5 to the family $\pi\colon{\cal X}\to S$. 
Then we deduce from (4.3) that
$$
\begin{aligned}
c_{1}(\overline{\lambda(\Omega_{X/S}^{\bullet})})
&=
c_{1}(\lambda(\Omega_{X/S}^{\bullet})_{Q})
-dd^{c}\log A({\cal X}/S)
\\
&=
-\frac{1}{12}\pi_{*}[c_{1}({\cal X}/S)\,
c_{n}({\cal X}/S)]-dd^{c}\log A({\cal X}/S)
\\
&
=
\frac{\chi(X)}{12}\,\omega_{{\rm WP},\,{\cal X}/S}.
\end{aligned}
$$
This completes the proof of Theorem 4.9.
\end{pf}

\begin{thm}
Let $X$ be a smooth Calabi-Yau $n$-fold. 
The Hermitian metric
$\|\cdot\|_{\lambda(\Omega^{\bullet}_{X})}$ on
$\lambda(\Omega^{\bullet}_{X})$ is independent of 
the choice of a K\"ahler metric on $X$. 
In particular, the BCOV Hermitian line
$\overline{\lambda(\Omega^{\bullet}_{X})}$
is an invariant of $X$.
\end{thm}

\begin{pf}
Let $\sigma\in
\lambda(\Omega_{X}^{\bullet})\setminus\{0\}$.
Let ${\cal X}=X\times{\Bbb P}^{1}\to{\Bbb P}^{1}$ 
be the trivial family over ${\Bbb P}^{1}$.
Let $\gamma_{0}$, $\gamma_{\infty}$ be arbitrary
K\"ahler forms on $X$. Let 
$\gamma_{{\cal X}/{\Bbb P}^{1}}=
\{\gamma_{t}\}_{t\in{\Bbb P}^{1}}$
be a $C^{\infty}$-family of K\"ahler forms on $X$ 
connecting $\gamma_{0}$ and $\gamma_{\infty}$. 
Since $\omega_{{\rm WP},{\cal X}/{\Bbb P}^{1}}=0$,
$\log\|\sigma\|_{\lambda(
\Omega_{{\cal X}/{\Bbb P}^{1}}^{\bullet})}^{2}$ 
is a harmonic function on ${\Bbb P}^{1}$ by 
Theorem 4.9. Hence 
$\|\sigma\|_{\lambda(
\Omega_{{\cal X}/{\Bbb P}^{1}}^{\bullet})}$ 
is a constant function on ${\Bbb P}^{1}$. 
This proves Theorem 4.10.
\end{pf}

\subsection
{}{\bf The BCOV invariant of Calabi-Yau threefolds}
\par
In Subsection 4.4, we fix $n=3$.
Hence $X$ is a smooth Calabi-Yau threefold. 
Set $b_{2}(X):=\dim H^{2}(X,{\Bbb R})$.
Let $c_{X}(\cdot,\cdot,\cdot)$ be
the cubic form on $H^{2}(X,{\Bbb R})$ induced from 
the cup-product:
$$
c_{X}(\alpha,\beta,\gamma):=
\frac{1}{(2\pi)^{3}}
\int_{X}\alpha\wedge\beta\wedge\gamma,
\qquad
\alpha,\beta,\gamma\in H^{2}(X,{\Bbb R}).
$$

\subsubsection
{The covolume of the cohomology lattice}
\par
Let $\kappa$ be a K\"ahler class on $X$.  Let
$\langle\cdot,\cdot\rangle_{L^{2},\kappa}$ be
the $L^{2}$-inner product on $H^{2}(X,{\Bbb R})$ 
with respect to $\kappa$, and let
$\langle\cdot,\cdot\rangle_{L^{2},\det\kappa}$ be
the induced $L^{2}$-inner product 
on $\det H^{2}(X,{\Bbb R})$.
Set $H^{2}(X,{\Bbb Z})_{\rm fr}:=H^{2}(X,{\Bbb Z})/
{\rm Torsion}$.

\begin{defn}
For a basis
$\{{\bf e}_{1},\ldots,{\bf e}_{b_{2}(X)}\}$ of
$H^{2}(X,{\Bbb Z})_{\rm fr}$ over $\Bbb Z$,
set
$$
{\rm Vol}_{L^{2}}(H^{2}(X,{\Bbb Z}),\kappa):
=
\det\left(
\langle{\bf e}_{i},{\bf e}_{j}\rangle_{L^{2},\kappa}
\right)=
\langle
{\bf e}_{1}\wedge\cdots\wedge{\bf e}_{b_{2}(X)},
{\bf e}_{1}\wedge\cdots\wedge{\bf e}_{b_{2}(X)}
\rangle_{L^{2},\det\kappa}.
$$
\end{defn}

Obviously, ${\rm Vol}_{L^{2}}(H^{2}(X,{\Bbb Z}),\kappa)$ 
is independent of the choice of a ${\Bbb Z}$-basis of
$H^{2}(X,{\Bbb Z})_{\rm fr}$; 
it is the volume of the real torus
$H^{2}(X,{\Bbb R})/H^{2}(X,{\Bbb Z})_{\rm fr}$ 
with respect to
$\langle\cdot,\cdot\rangle_{L^{2},\kappa}$.
We can write
${\rm Vol}_{L^{2}}(H^{2}(X,{\Bbb Z}),\kappa)$
in terms of the cubic form $c_{X}$ as follows:
\par
Let $L$ be the operator on $H^{\bullet}(X,{\Bbb R})$ 
defined by
$L(\varphi)=\kappa\wedge\varphi$ for 
$\varphi\in H^{\bullet}(X,{\Bbb R})$.

\begin{lem}
The following identity holds
$$
\langle\alpha,\beta\rangle_{L^{2},\kappa}=
\frac{3}{2}\frac{c_{X}(\alpha,\kappa,\kappa)\,
c_{X}(\beta,\kappa,\kappa)}
{c_{X}(\kappa,\kappa,\kappa)}-
c_{X}(\alpha,\beta,\kappa),
\qquad
\alpha,\beta\in H^{2}(X,{\Bbb R}).
$$
In particular,
${\rm Vol}_{L^{2}}(H^{2}(X,{\Bbb Z}),\kappa)\in{\Bbb Q}$
if $\kappa\in H^{2}(X,{\Bbb Q})$.
\end{lem}

\begin{pf}
Let $\varphi\in H^{1,1}(X,{\Bbb R})=H^{2}(X,{\Bbb R})$. 
By \cite[Lemma\,6.31]{Voisin02}, one has the orthogonal
decomposition
$H^{1,1}(X,{\Bbb R})=\ker(L^{2})\oplus{\Bbb R}\kappa$
with respect to
$\langle\cdot,\cdot\rangle_{L^{2},\kappa}$.
Since
\begin{equation}
\langle\varphi,\varphi\rangle_{L^{2},\kappa}
=
\left\{
\begin{array}{ll}
-c_{X}(\varphi,\varphi,\kappa)
&(\varphi\in\ker(L^{2}))
\\
\frac{1}{2}c_{X}(\varphi,\varphi,\kappa)
&(\varphi\in{\Bbb R}\kappa)
\end{array}
\right.
\end{equation}
by \cite[Th.\,6.32]{Voisin02}, we get the decomposition
\begin{equation}
\varphi
=
\left(
\varphi-
\frac{c_{X}(\varphi,\kappa,\kappa)}
{c_{X}(\kappa,\kappa,\kappa)}\kappa
\right)+
\frac{c_{X}(\varphi,\kappa,\kappa)}
{c_{X}(\kappa,\kappa,\kappa)}\kappa
\in
\ker(L^{2})\oplus{\Bbb R}\kappa.
\end{equation}
By (4.4), (4.5), we get
$$
\begin{aligned}
\langle\alpha,\beta\rangle_{L^{2},\kappa}
&=
-c_{X}\left(
\alpha-
\frac{c_{X}(\alpha,\kappa,\kappa)}
{c_{X}(\kappa,\kappa,\kappa)}\kappa\,,\,
\beta-
\frac{c_{X}(\beta,\kappa,\kappa)}
{c_{X}(\kappa,\kappa,\kappa)}\kappa\,,\,\kappa
\right)
\\
&\quad+
\frac{1}{2}c_{X}\left(
\frac{c_{X}(\alpha,\kappa,\kappa)}
{c_{X}(\kappa,\kappa,\kappa)}\kappa\,,\,
\frac{c_{X}(\beta,\kappa,\kappa)}
{c_{X}(\kappa,\kappa,\kappa)}\kappa\,,\,\kappa
\right)
\\
&=
\frac{3}{2}\frac{c_{X}(\alpha,\kappa,\kappa)\,
c_{X}(\beta,\kappa,\kappa)}
{c_{X}(\kappa,\kappa,\kappa)}
-c_{X}(\alpha,\beta,\kappa).
\end{aligned}
$$
This proves the lemma.
\end{pf}

\subsubsection
{The BCOV invariant}
\par
Let us introduce the main object of this paper.

\begin{defn}
For a K\"ahler form $\gamma$ on $X$, 
the {\it BCOV invariant} of $(X,\gamma)$ is the real
number defined by
$$
\begin{aligned}
\tau_{\rm BCOV}(X,\gamma):
&=
{\rm Vol}(X,\gamma)^{-3}\,
{\rm Vol}_{L^{2}}(H^{2}(X,{\Bbb Z}),[\gamma])^{-1}\,
{\cal A}(X,\gamma)\,{\cal T}_{\rm BCOV}(X,\gamma)
\\
&=
{\rm Vol}(X,\gamma)^{\frac{\chi(X)}{12}-3}\,
{\rm Vol}_{L^{2}}(H^{2}(X,{\Bbb Z}),[\gamma])^{-1}
\\
&\quad\times
\exp\left[-\frac{1}{12}\int_{X}
\log\left(\frac{\sqrt{-1}\,\eta\wedge\bar{\eta}}
{\gamma^{3}/3!}
\cdot
\frac{{\rm Vol}(X,\gamma)}{\|\eta\|_{L^{2}}^{2}}
\right)\,c_{3}(X,\gamma)\right]\,
{\cal T}_{\rm BCOV}(X,\gamma).
\end{aligned}
$$
\end{defn}

In the rest of Section 4, we derive
a variational formula for the BCOV invariant.

\subsubsection
{The curvature formula for the BCOV invariant}
\par
Let $\pi\colon({\cal X},X)\to(S,0)$ be a flat deformation 
of $X$ which is induced from the Kuranishi family
by a holomorphic map $\mu\colon(S,0)\to({\rm Def}(X),[X])$.
Let $\omega_{{\rm H},{\cal X}/S}$ be the $(1,1)$-form
on $S$ induced from the Hodge form on ${\rm Def}(X)$
via $\mu$:
$$
\omega_{{\rm H},{\cal X}/S}:=\mu^{*}\omega_{\rm H}.
$$
\par 
Let $g_{{\cal X}/S}$ be a fiberwise-K\"ahler metric 
on $T{\cal X}/S$. Let
$\gamma_{s}$ be the K\"ahler form of 
$g_{{\cal X}/S}|_{X_{s}}$.
Let $\tau_{\rm BCOV}({\cal X}/S)$ 
be the function on $S$ defined by
$$
\tau_{\rm BCOV}({\cal X}/S)(s):=
\tau_{\rm BCOV}(X_{s},\gamma_{s}),
\qquad s\in S.
$$

\begin{thm}
The following identity of $(1,1)$-forms on $(S,0)$ holds
$$
\begin{aligned}
dd^{c}\log\tau_{\rm BCOV}({\cal X}/S)
&
=
-\frac{\chi(X)}{12}\,\omega_{{\rm WP},{\cal X}/S}
-\omega_{{\rm H},{\cal X}/S}
\\
&
=
-\left(h^{1,2}(X)+\frac{\chi(X)}{12}+3\right)
\,\mu^{*}\omega_{\rm WP}
-\mu^{*}{\rm Ric}\,\omega_{\rm WP}.
\end{aligned}
$$
\end{thm}

\begin{pf}
We follow \cite[Th.\,5.6]{Yoshikawa04}.
Let ${\cal A}({\cal X}/S)$ and 
${\cal T}_{\rm BCOV}({\cal X}/S)$ be 
the $C^{\infty}$ functions on $S$ defined by
$$
{\cal A}({\cal X}/S)(s):=
{\cal A}(X_{s},\gamma_{s}),
\qquad
{\cal T}_{\rm BCOV}({\cal X}/S)(s):=
{\cal T}_{\rm BCOV}(X_{s},\gamma_{s})
$$
for $s\in S$.
By Theorems 3.5 and 4.9, we get
$$
\begin{aligned}
\,&
-dd^{c}\log[{\cal A}({\cal X}/S)\,
{\cal T}_{\rm BCOV}({\cal X}/S)]+
\sum_{p,q\geq0}(-1)^{p+q}p\,
c_{1}(\det R^{q}\pi_{*}\Omega^{p}_{{\cal X}/S},
\|\cdot\|_{L^{2},g_{{\cal X}/S}})
\\
&=
\frac{\chi(X)}{12}\mu^{*}\omega_{\rm WP}.
\end{aligned}
$$
Since $R^{q}\pi_{*}\Omega^{p}_{{\cal X}/S}\not=0$
if and only if $p+q=3$ or $p=q$, we deduce from
Proposition 4.6 that
\begin{equation}
\begin{aligned}
\,&
-dd^{c}\log[{\cal A}({\cal X}/S)\,
{\cal T}_{\rm BCOV}({\cal X}/S)]
+\sum_{p>0}p\,
c_{1}(\det R^{p}\pi_{*}\Omega^{p}_{{\cal X}/S},
\|\cdot\|_{L^{2},g_{{\cal X}/S}})
\\
&-\left(\mu^{*}{\rm Ric}\,\omega_{\rm WP}+
h^{1,2}(X)\,\mu^{*}\omega_{\rm WP}\right)
-3\mu^{*}\omega_{\rm WP}
\\
&=
\frac{\chi(X)}{12}\mu^{*}\omega_{\rm WP}.
\end{aligned}
\end{equation}
\par
Define a function
${\rm Vol}_{L^{2}}(H^{2}({\cal X}/S,{\Bbb Z}))$
on $S$ by
$$
{\rm Vol}_{L^{2}}(H^{2}({\cal X}/S,{\Bbb Z}))(s)
:=
{\rm Vol}_{L^{2}}
(H^{2}(X_{s},{\Bbb Z}),[\gamma_{s}]),
\qquad
s\in S.
$$
\par
Since $\pi\colon{\cal X}\to S$ is induced from the Kuranishi
family, there exist holomorphic line bundles 
${\cal L}_{1},\ldots,{\cal L}_{b_{2}(X)}$ on $\cal X$ 
by Lemma 2.16
such that $c_{1}({\cal L}_{i})|_{X}={\bf e}_{i}$
for $1\leq i\leq b_{2}(X)$, and such that 
${\frak C}_{1}({\cal L}_{1})\wedge\cdots\wedge 
{\frak C}_{1}({\cal L}_{b_{2}(X)})$ is a nowhere 
vanishing holomorphic section of 
$R^{1}\pi_{*}\Omega^{1}_{{\cal X}/S}$.
Then
\begin{equation}
\|{\frak C}_{1}({\cal L}_{1})\wedge\cdots\wedge 
{\frak C}_{1}({\cal L}_{b_{2}(X)})
\|_{L^{2},g_{{\cal X}/S}}^{2}
=
{\rm Vol}_{L^{2}}(H^{2}({\cal X}/S,{\Bbb Z})).
\end{equation}
\par
By the Serre duality and (3.5),
${\bf 1}_{(1,1),(2,2)}\otimes
({\frak C}_{1}({\cal L}_{1})\wedge\cdots\wedge 
{\frak C}_{1}({\cal L}_{b_{2}(X)}))^{-1}$ 
is a nowhere vanishing holomorphic section of
$R^{2}\pi_{*}\Omega^{2}_{{\cal X}/S}$ such that
\begin{equation}
\|{\bf 1}_{(1,1),(2,2)}\otimes
({\frak C}_{1}({\cal L}_{1})\wedge\cdots\wedge 
{\frak C}_{1}({\cal L}_{b_{2}(X)}))^{-1}
\|_{L^{2},g_{{\cal X}/S}}^{2}
=
{\rm Vol}_{L^{2}}
(H^{2}({\cal X}/S,{\Bbb Z}))^{-1}.
\end{equation}
\par
Let ${\rm Vol}({\cal X}/S,\gamma_{{\cal X}/S})$ be
the function on $S$ defined by
$$
{\rm Vol}({\cal X}/S)(s)
:=
{\rm Vol}(X_{s},\gamma_{s}).
$$
Then 
$\frac{\gamma_{{\cal X}/S}^{3}}{3!{\rm Vol}({\cal X}/S)}$
is a nowhere vanishing holomorphic section of
$R^{3}\pi_{*}\Omega^{3}_{{\cal X}/S}$ such that
\begin{equation}
\left\|
\frac{\gamma_{{\cal X}/S}^{3}}{3!{\rm Vol}({\cal X}/S)}
\right\|_{L^{2},g_{{\cal X}/S}}^{2}
=
{\rm Vol}({\cal X}/S)^{-1}.
\end{equation}
Substituting (4.7), (4.8), (4.9) into (4.6), 
we get the equation:
\begin{equation}
\begin{aligned}
\,&
-dd^{c}\log[{\cal A}({\cal X}/S)\,
{\cal T}_{\rm BCOV}({\cal X}/S)]
+dd^{c}\log\,{\rm Vol}_{L^{2}}(H^{2}({\cal X}/S,{\Bbb Z}))
+3dd^{c}\log\,{\rm Vol}({\cal X}/S)
\\
&=
\left(h^{1,2}(X)+\frac{\chi(X)}{12}+3\right)
\,\mu^{*}\omega_{\rm WP}
+\mu^{*}{\rm Ric}\,\omega_{\rm WP}.
\end{aligned}
\end{equation}
The theorem follows from the definition of 
the BCOV invariant and (4.10).
\end{pf}

\begin{rem}
If we follow the mirror symmetry and if
$X^{\lor}$ is the mirror Calabi-Yau threefold of $X$,
the coefficient of $\mu^{*}\omega_{\rm WP}$ in (4.10) 
is compatible with that of \cite[Eq.\,(14)]{BCOV93}
since $h^{1,1}(X^{\lor})=h^{1,2}(X)$ and
$\chi(X^{\lor})=-\chi(X)$.
\end{rem}

For a higher dimensional analogue of Theorem 4.14,
we refer to \cite{FangLu03}.

\begin{thm}
The BCOV invariant $\tau_{\rm BCOV}(X,\gamma)$ is 
independent of the choice of a K\"ahler metric on $X$. 
In particular, $\tau_{\rm BCOV}(X,\gamma)$
is an invariant of $X$.
\end{thm}

\begin{pf}
Let ${\cal X}=X\times{\Bbb P}^{1}\to{\Bbb P}^{1}$ 
be the trivial family over ${\Bbb P}^{1}$.
Let $\gamma_{0}$, $\gamma_{\infty}$ be arbitrary
K\"ahler forms on $X$. Let 
$\gamma_{{\cal X}/{\Bbb P}^{1}}=
\{\gamma_{t}\}_{t\in{\Bbb P}^{1}}$
be a $C^{\infty}$-family of K\"ahler forms on $X$ 
connecting $\gamma_{0}$ and $\gamma_{\infty}$. 
Since $\mu^{*}\omega_{\rm WP}$ and
$\mu^{*}{\rm Ric}(\omega_{\rm WP})$ are independent
of $t$,
$\log\tau_{\rm BCOV}({\cal X}/{\Bbb P}^{1})$ 
is a harmonic function on ${\Bbb P}^{1}$ by 
Theorem 4.14. Hence 
$\tau_{\rm BCOV}({\cal X}/{\Bbb P}^{1})$ 
is a constant function on ${\Bbb P}^{1}$. 
\end{pf}

After Theorem 4.16,
we shall write $\tau_{\rm BCOV}(X)$ for
$\tau_{\rm BCOV}(X,\gamma)$ in the rest of this paper.

%%%%%%%%%%%%%%%%%%%%%%%%%%%%%%%%%%%%%%%%%%%%%%%%%%%%%%%%%%%
%
%                 Section 5
%
%%%%%%%%%%%%%%%%%%%%%%%%%%%%%%%%%%%%%%%%%%%%%%%%%%%%%%%%%%%

\section
{\bf The singularity of the Quillen metric 
on the BCOV bundle}
\par
In Section 5, we fix the following notation:
Let $\cal X$ be a compact K\"ahler manifold of 
dimension $n+1$ and let $S$ be a compact Riemann surface.
Let $\pi\colon{\cal X}\to S$ be a surjective holomorphic 
map, and we do {\it not} assume that a general fiber 
of $\pi$ is Calabi-Yau. 
\par
Let $\Sigma_{\pi}$ be the critical locus of $\pi$,
and set 
$$
{\cal D}:=\pi(\Sigma_{\pi}),
\qquad
S^{o}:=S\setminus{\cal D},
\qquad
{\cal X}^{o}:=\pi^{-1}(S^{o}),
\qquad
\pi^{o}:=\pi|_{{\cal X}^{o}}.
$$
Then $\pi^{o}\colon{\cal X}^{o}\to S^{o}$ is a holomorphic 
family of compact complex manifolds, and
$\Omega^{1}_{{\cal X}^{o}/S^{o}}$ is a holomorphic
vector bundle of rank $n$ over ${\cal X}^{o}$. 
\par
As in Sections 3 and 4,
we have the holomorphic line bundles on $S^{o}$:
$$
\lambda(\Omega^{p}_{{\cal X}^{o}/S^{o}})
=
\otimes_{q=0}^{n}
(\det R^{q}\pi_{*}\Omega^{p}_{{\cal X}^{o}/S^{o}})^{(-1)^{q}},
\qquad
\lambda(\Omega_{{\cal X}^{o}/S^{o}}^{\bullet})
=
\otimes_{p=0}^{n}
\lambda(\Omega^{p}_{{\cal X}^{o}/S^{o}})^{(-1)^{p}p}.
$$
In this section, we construct holomorphic extensions of 
$\lambda(\Omega^{p}_{{\cal X}^{o}/S^{o}})$ and
$\lambda(\Omega^{\bullet}_{{\cal X}^{o}/S^{o}})$ from
$S^{o}$ to $S$, and we study the singularity of 
the corresponding Quillen metrics.

\subsection
{}
{\bf The K\"ahler extension of 
the determinant line bundles}
\par
Since $\Omega_{{\cal X}/S}^{1}=
\Omega^{1}_{\cal X}/\pi^{*}\Omega^{1}_{S}$,
we have the following complex of coherent sheaves 
on $\cal X$, which is acyclic on ${\cal X}$
(cf. \cite[p.94 l.12-l.16]{Looijenga84}):
$$
0\longrightarrow
\pi^{*}\Omega^{1}_{S}
\longrightarrow
\Omega_{\cal X}^{1}
\longrightarrow
\Omega_{{\cal X}/S}^{1}
\longrightarrow0.
$$

\begin{defn}
(i) 
For $p>0$, let ${\cal E}^{p}_{{\cal X}/S}$ be 
the complex of holomorphic vector bundles on $\cal X$ 
defined by
$$
{\cal E}^{p}_{{\cal X}/S}\colon
(\pi^{*}\Omega^{1}_{S})^{\otimes p}\longrightarrow
\Omega_{\cal X}^{1}\otimes
(\pi^{*}\Omega^{1}_{S})^{\otimes (p-1)}
\longrightarrow\cdots\longrightarrow
\Omega_{\cal X}^{p-1}\otimes\pi^{*}\Omega_{S}^{1}
\longrightarrow
\Omega_{\cal X}^{p},
$$
where the maps
$\Omega^{i}_{\cal X}\otimes
(\pi^{*}\Omega^{1}_{S})^{\otimes(p-i)}\to
\Omega^{i+1}_{\cal X}\otimes
(\pi^{*}\Omega^{1}_{S})^{\otimes(p-i-1)}$ are given by
$$
\omega\otimes(\pi^{*}\xi)^{\otimes(p-i)}
\mapsto(\omega\wedge\pi^{*}\xi)\otimes
(\pi^{*}\xi)^{\otimes(p-i-1)},
\qquad
\omega\in\Omega^{i}_{\cal X},
\quad
\xi\in\Omega^{1}_{S}.
$$
For $p=0$, set
${\cal E}^{0}_{{\cal X}/S}\colon 
0\rightarrow{\cal O}_{{\cal X}}\rightarrow 0$.
\newline
{(ii) }
For $p\geq 0$, let ${\cal F}^{p}_{{\cal X}/S}$ be 
the complex of coherent sheaves on $\cal X$ 
defined by
$$
\CD
{\cal F}^{p}_{{\cal X}/S}\colon
0\; @>>>\;{\cal E}^{p}_{{\cal X}/S}\;
@>r>>\;\Omega_{{\cal X}/S}^{p}\;@>>>\;0,
\endCD
$$
where 
$r\colon\Omega^{p}_{\cal X}\to\Omega^{p}_{{\cal X}/S}$
is the quotient map for $p>0$ and the identity map
for $p=0$. 
\end{defn}

Since ${\rm rk}(\pi^{*}\Omega^{1}_{S})=1$,
${\cal F}^{p}_{{\cal X}/S}$ is acyclic on 
${\cal X}\setminus\Sigma_{\pi}$ for $p>1$ and
on ${\cal X}$ for $p=0,1$.

\begin{defn}
(i) Let $\lambda({\cal E}^{p}_{{\cal X}/S})$ be
the holomorphic line bundle on $S$ defined by
$$
\lambda({\cal E}^{p}_{{\cal X}/S})
:=
\bigotimes_{i=0}^{p}
\lambda(\Omega_{\cal X}^{p-i}\otimes
(\pi^{*}\Omega^{1}_{S})^{\otimes i})^{(-1)^{i}}.
$$
(ii) Let $\lambda(\Omega_{{\cal X}/S}^{\bullet})$
be the holomorphic line bundle on $S$ defined by
$$
\lambda(\Omega_{{\cal X}/S}^{\bullet}):=
\bigotimes_{p\geq0}
\lambda({\cal E}^{p}_{{\cal X}/S})^{(-1)^{p}p}.
$$
We call $\lambda({\cal E}^{p}_{{\cal X}/S})$ and
$\lambda(\Omega_{{\cal X}/S}^{\bullet})$
the {\it K\"ahler extensions} of
$\lambda(\Omega_{{\cal X}^{o}/S^{o}}^{p})$ and
$\lambda(\Omega_{{\cal X}^{o}/S^{o}}^{\bullet})$ 
from $S^{o}$ to $S$, respectively.
\end{defn}

Since ${\cal F}^{p}_{{\cal X}/S}$ is acyclic on 
${\cal X}\setminus\Sigma_{\pi}$, we have the canonical
isomorphisms of holomorphic line bundles on $S^{o}$:
$$
\lambda(\Omega^{p}_{{\cal X}^{o}/S^{o}})
\cong
\lambda({\cal E}^{p}_{{\cal X}/S})|_{S^{o}},
\qquad
\lambda(\Omega_{{\cal X}^{o}/S^{o}}^{\bullet})
\cong
\lambda(\Omega_{{\cal X}/S}^{\bullet})|_{S^{o}}.
$$
\par
Let $g_{\cal X}$ be a K\"ahler metric on $\cal X$.
Let $g_{{\cal X}/S}:=g_{\cal X}|_{T{\cal X}/S}$ be 
the Hermitian metric on 
$T{\cal X}/S|_{{\cal X}\setminus\Sigma_{\pi}}$ 
induced from $g_{\cal X}$. 
Then $g_{{\cal X}/S}$ (resp. $g_{\cal X}$) induces 
the Hermitian metric $g_{\Omega^{p}_{{\cal X}/S}}$ 
(resp. $g_{\Omega^{p}_{\cal X}}$) on
$\Omega^{p}_{{\cal X}/S}|_{{\cal X}\setminus\Sigma_{\pi}}$ 
(resp. $\Omega^{p}_{\cal X}$) for all $p\geq0$. 
\par
Following Bismut \cite{Bismut97} and 
Yoshikawa \cite{Yoshikawa05}, we determine 
the singularity of the Quillen metric on
$\lambda(\Omega^{p}_{{\cal X}^{o}/S^{o}})$ 
near ${\cal D}$ with respect to the K\"ahler extension
and with respect to the metrics $g_{{\cal X}/S}$, 
$g_{\Omega^{p}_{{\cal X}/S}}$.

\subsection
{}
{\bf Three Quillen metrics on the extended BCOV bundles}
\par
Let $0\in{\cal D}$.
Let $({\cal U},t)$ be a coordinate neighborhood 
of $0$ in $S$ centered at $0$ such that
${\cal U}\cong\varDelta$ and
${\cal U}\cap{\cal D}=\{0\}$.
We set ${\cal U}^{o}:={\cal U}\setminus{\cal D}
={\cal U}\setminus\{0\}$.
\par
Let $k_{S}$ be a Hermitian metric on $\Omega^{1}_{S}$ 
such that $k_{S}(dt,dt)=1$ on $\cal U$.
Then $\pi^{*}k_{S}$ is a Hermitian metric 
on $\pi^{*}\Omega^{1}_{S}$. 
Let $g_{\pi^{*}\Omega^{1}_{S}}$ be 
the Hermitian metric on 
$\pi^{*}\Omega^{1}_{S}|_{{\cal X}\setminus\Sigma_{\pi}}$ 
induced from $g_{\Omega^{1}_{\cal X}}$ by the inclusion
$\pi^{*}\Omega^{1}_{S}\subset\Omega^{1}_{\cal X}$.
Since
$$
\pi^{*}k_{S}(d\pi,d\pi)=\pi^{*}\{k_{S}(dt,dt)\}=1,
\qquad
g_{\pi^{*}\Omega^{1}_{S}}(d\pi,d\pi)=
g_{\Omega^{1}_{\cal X}}(d\pi,d\pi)=\|d\pi\|^{2}
$$
on $\pi^{-1}(\cal U)$, 
the following identity holds on $\pi^{-1}({\cal U})$:
$$
g_{\pi^{*}\Omega^{1}_{S}}=\|d\pi\|^{2}\,\pi^{*}k_{S}.
$$
\par
We define three Quillen metrics on
the K\"ahler extension
$\lambda({\cal E}^{p}_{{\cal X}/S})|_{{\cal U}^{o}}$
as follows.

\begin{defn}
(i) 
Let 
$\|\cdot\|^{2}_{\lambda(\Omega^{p}_{{\cal X}^{o}/S^{o}}),
Q,g_{{\cal X}/S}}$ be the Quillen metric on 
$\lambda(\Omega^{p}_{{\cal X}^{o}/S^{o}})|_{{\cal U}^{o}}$ 
with respect to
$g_{{\cal X}/S}$ and $g_{\Omega^{p}_{{\cal X}/S}}$.
Let 
$\|\cdot\|^{2}_{\lambda({\cal E}^{p}_{{\cal X}/S}),Q,
g_{{\cal X}/S}}$
be the Quillen metric on 
$\lambda({\cal E}^{p}_{{\cal X}/S})|_{{\cal U}^{o}}$ 
induced from
$\|\cdot\|^{2}_{\lambda(\Omega^{p}_{{\cal X}^{o}/S^{o}}),
Q,g_{{\cal X}/S}}$ by the canonical isomorphism 
$\lambda(\Omega_{{\cal X}^{o}/S^{o}}^{p})|_{{\cal U}^{o}}
\cong
\lambda({\cal E}^{p}_{{\cal X}/S})|_{{\cal U}^{o}}$:
$$
\|\cdot\|^{2}_{\lambda({\cal E}^{p}_{{\cal X}/S}),Q,
g_{{\cal X}/S}}
:=
\|\cdot\|^{2}_{\lambda(\Omega^{p}_{{\cal X}^{o}/S^{o}}),
Q,g_{{\cal X}/S}}.
$$
\newline{(ii) }
Let 
$\|\cdot\|^{2}_{\lambda(\Omega^{p-i}_{\cal X}\otimes
(\pi^{*}\Omega^{1}_{S})^{\otimes i}),Q,\pi^{*}k_{S}}$ 
be the Quillen metric on 
$\lambda(\Omega^{p-i}_{\cal X}\otimes
(\pi^{*}\Omega^{1}_{S})^{\otimes i})|_{{\cal U}^{o}}$ 
with respect to $g_{{\cal X}/S}$ and 
$g_{\Omega^{p-i}_{\cal X}}\otimes\pi^{*}k_{S}$. 
Set
$$
\|\cdot\|^{2}_{\lambda({\cal E}^{p}_{{\cal X}/S}),
Q,\pi^{*}k_{S}}
:=
\bigotimes_{i=0}^{p}
\|\cdot\|^{2(-1)^{i}}_{\lambda(\Omega^{p-i}_{\cal X}\otimes
(\pi^{*}\Omega^{1}_{S})^{\otimes i}),Q,\pi^{*}k_{S}}.
$$
\newline{(iii) }
Let 
$\|\cdot\|^{2}_{\lambda(\Omega^{p-i}_{\cal X}\otimes
(\pi^{*}\Omega^{1}_{S})^{\otimes i}),Q,
h_{\pi^{*}\Omega^{1}_{S}}}$ be the Quillen metric on 
$\lambda(\Omega^{p-i}_{\cal X}\otimes
(\pi^{*}\Omega^{1}_{S})^{\otimes i})|_{{\cal U}^{o}}$ 
with respect to $g_{{\cal X}/S}$ and
$g_{\Omega^{p-i}_{\cal X}}\otimes
g_{\pi^{*}\Omega^{1}_{S}}$. 
Set
$$
\|\cdot\|^{2}_{\lambda({\cal E}^{p}_{{\cal X}/S}),
Q,g_{\pi^{*}\Omega^{1}_{S}}}
:=
\bigotimes_{i=0}^{p}
\|\cdot\|^{2(-1)^{i}}_{\lambda(\Omega^{p-i}_{\cal X}
\otimes(\pi^{*}\Omega^{1}_{S})^{\otimes i}),Q,
g_{\pi^{*}\Omega^{1}_{S}}}.
$$
\end{defn}

When $p=0$, we have the following relations
$$
\|\cdot\|^{2}_{\lambda({\cal E}^{0}_{{\cal X}/S}),Q,
g_{{\cal X}/S}}=
\|\cdot\|^{2}_{\lambda({\cal E}^{0}_{{\cal X}/S}),
Q,\pi^{*}k_{S}}=
\|\cdot\|^{2}_{\lambda({\cal E}^{0}_{{\cal X}/S}),
Q,g_{\pi^{*}\Omega^{1}_{S}}}=
\|\cdot\|^{2}_{\lambda({\cal O}_{\cal X}),
Q,g_{{\cal X}/S}}.
$$
We shall prove that 
$\log\|\cdot\|^{2}_{\lambda({\cal E}^{p}_{{\cal X}/S}),
Q,g_{{\cal X}/S}}$ has logarithmic singularities 
at $0\in{\cal D}$, whose coefficients are determined 
by the resolution data of the Gauss map.

\subsection
{}{\bf The Gauss maps and their resolutions}
\par
Let 
$\varPi\colon{\Bbb P}(\Omega^{1}_{\cal X})\to{\cal X}$
be the projective bundle associated with the holomorphic
cotangent bundle $\Omega^{1}_{\cal X}$.
Let $\varPi^{\lor}\colon{\Bbb P}(T{\cal X})\to{\cal X}$
be the projective bundle associated with the holomorphic
tangent bundle $T{\cal X}$. Then the fiber 
${\Bbb P}(T_{x}{\cal X})^{\lor}$ is the set of all 
hyperplanes of $T_{x}{\cal X}$ containing 
$0_{x}\in T_{x}{\cal X}$. We have 
${\Bbb P}(\Omega^{1}_{\cal X})\cong
{\Bbb P}(T{\cal X})^{\lor}$.
\par
We define the Gauss maps
$\nu\colon{\cal X}\setminus\Sigma_{\pi}\to
{\Bbb P}(\Omega^{1}_{\cal X})$ and
$\mu\colon{\cal X}\setminus\Sigma_{\pi}\to
{\Bbb P}(T{\cal X})^{\lor}$
by
$$
\nu(x):=[d\pi_{x}]=
\left[
\sum_{i=0}^{n}\frac{\partial\pi}{\partial z_{i}}(x)\,
dz_{i}\right],
\qquad
\mu(x):=[T_{x}X_{\pi(x)}].
$$
Then $\nu=\mu$ under the canonical isomorphism
${\Bbb P}(\Omega^{1}_{\cal X})\cong
{\Bbb P}(T{\cal X})^{\lor}$.
\par
Let 
$L:={\cal O}_{{\Bbb P}(\Omega^{1}_{\cal X})}(-1)\subset
\varPi^{*}\Omega^{1}_{\cal X}$
be the tautological line bundle over 
${\Bbb P}(\Omega^{1}_{\cal X})$, and set
$Q:=\varPi^{*}\Omega^{1}_{\cal X}/L$. Then we have 
the following exact sequences $\cal S$
of holomorphic vector bundles 
on ${\Bbb P}(\Omega^{1}_{\cal X})$:
$$
{\cal S}\colon
0\longrightarrow
L\longrightarrow 
\varPi^{*}\Omega^{1}_{\cal X}\longrightarrow
Q\longrightarrow0.
$$
Let $p\leq n$.
Since ${\rm rk}(L)=1$, this induces
the following exact sequence of holomorphic vector 
bundles on ${\Bbb P}(\Omega^{1}_{\cal X})$: 
$$
{\cal K}^{p}\colon
0
\longrightarrow
L^{p}
\longrightarrow
\varPi^{*}\Omega^{1}_{\cal X}\otimes L^{p-1}
\longrightarrow\cdots\longrightarrow
\varPi^{*}\Omega^{p-1}_{\cal X}\otimes L
\longrightarrow 
\varPi^{*}\Omega^{p}_{\cal X}
\longrightarrow
\bigwedge^{p}Q
\longrightarrow0,
$$
where $\varPi^{*}\Omega^{p}_{\cal X}\to
\bigwedge^{p}Q$ is the quotient map and
$\varPi^{*}\Omega^{i}_{\cal X}\otimes L^{p-i}
\to\varPi^{*}\Omega^{i+1}_{\cal X}\otimes L^{p-i-1}$ 
is given by
$\omega\otimes\sigma^{\otimes(p-i)}\mapsto
(\omega\wedge\sigma)\otimes\sigma^{\otimes(p-i-1)}$ 
for $\omega\in\varPi^{*}\Omega^{1}_{\cal X}$ 
and $\sigma\in L$. Then
$$
{\cal F}^{p}_{{\cal X}/S}=\nu^{*}{\cal K}^{p}.
$$
\par
Similarly,
let $H:={\cal O}_{{\Bbb P}(\Omega^{1}_{\cal X})}(1)$,
and let $U$ be the universal hyperplane bundle of
$(\varPi^{\lor})^{*}T{\cal X}$. Then the dual of 
${\cal S}$ is given by
$$
{\cal S}^{\lor}\colon
0\longrightarrow U\longrightarrow
(\varPi^{\lor})^{*}T{\cal X}\longrightarrow
H\longrightarrow 0.
$$
Since $T_{x}{\cal X}/S=\{v\in T_{x}{\cal X};\,
d\pi_{x}(v)=0\}$, we have
$$
T{\cal X}/S=\mu^{*}U.
$$
\par
Let $g_{U}$ be the Hermitian metric on $U$ induced from 
$(\varPi^{\lor})^{*}g_{\cal X}$, and let $g_{H}$ be
the Hermitian metric on $H$ induced from
$(\varPi^{\lor})^{*}g_{\cal X}$ by 
the $C^{\infty}$-isomorphism $H\cong U^{\perp}$.
\par
Let $g_{L}$ be the Hermitian metric on $L$ induced from
$\varPi^{*}g_{\Omega^{1}_{\cal X}}$ by the inclusion 
$L\subset\varPi^{*}\Omega^{1}_{\cal X}$. 
Let $g_{Q}$ be the Hermitian metric on $Q$ induced from
$\varPi^{*}g_{\Omega^{1}_{\cal X}}$ 
by the $C^{\infty}$-isomorphism $Q\cong L^{\perp}$.
We consider the Hermitian metric 
$g_{\varPi^{*}\Omega^{i}_{\cal X}\otimes L^{p-i}}$
on $\varPi^{*}\Omega^{i}_{\cal X}\otimes L^{p-i}$ 
induced from $\varPi^{*}g_{\Omega^{1}_{\cal X}}$, 
$g_{L}$, and we consider
the Hermitian metric $g_{\wedge^{p}Q}$
on $\bigwedge^{p}Q$ induced from $g_{Q}$.
We define $\overline{\cal K}^{p}$ to be
the exact sequence ${\cal K}^{p}$ equipped with 
the Hermitian metrics
$\{g_{\varPi^{*}\Omega^{i}_{\cal X}\otimes L^{p-i}}\}$
and $g_{\wedge^{p}Q}$.
Then we have the following isomorphisms of
Hermitian vector bundles over
${\cal X}\setminus\Sigma_{\pi}$:
\begin{equation}
\overline{\cal F}^{p}_{{\cal X}/S}=
\nu^{*}\overline{\cal K}^{p},
\qquad
(T{\cal X}/S,g_{{\cal X}/S})=
\mu^{*}(U,g_{U}).
\end{equation}
Since $d\pi$ is a nowhere vanishing holomorphic section 
of $\nu^{*}L|_{{\cal X}\setminus\Sigma_{\pi}}$, 
we get the following equation 
on ${\cal X}\setminus\Sigma_{\pi}$
$$
-dd^{c}\log\|d\pi\|^{2}=\nu^{*}c_{1}(L,g_{L}).
$$
\par
Since $\Sigma_{\pi}$ is a proper analytic subset of
$\cal X$,
the Gauss maps 
$\nu\colon{\cal X}\setminus\Sigma_{\pi}\to
{\Bbb P}(\Omega^{1}_{\cal X})$ 
and
$\mu\colon{\cal X}\setminus\Sigma_{\pi}\to
{\Bbb P}(T{\cal X})^{\lor}$ extend to 
meromorphic maps
$\nu\colon{\cal X}\dashrightarrow
{\Bbb P}(\Omega^{1}_{\cal X})$
and
$\mu\colon
{\cal X}\dashrightarrow{\Bbb P}(T{\cal X})^{\lor}$
by e.g. \cite[Th.\,4.5.3]{NoguchiOchiai90}.
By Hironaka, there exist a projective algebraic manifold 
$\widetilde{\cal X}$, a divisor of normal crossing 
$E\subset{\cal X}$, a birational holomorphic map
$q\colon\widetilde{\cal X}\to{\cal X}$, 
and holomorphic maps
$\widetilde{\nu}\colon\widetilde{\cal X}\to
{\Bbb P}(\Omega^{1}_{\cal X})$ and
$\widetilde{\mu}\colon\widetilde{\cal X}\to
{\Bbb P}(T{\cal X})^{\lor}$ satisfying
the following conditions:
\newline{(i) }
$q|_{\widetilde{\cal X}\setminus q^{-1}(\Sigma_{\pi})}
\colon
\widetilde{\cal X}\setminus q^{-1}(\Sigma_{\pi})
\to{\cal X}\setminus\Sigma_{\pi}$ is an isomorphism;
\newline{(ii) }
$q^{-1}(\Sigma_{\pi})=E$;
\newline{(iii) }
$\widetilde{\nu}=\nu\circ q$ and
$\widetilde{\mu}=\mu\circ q$ on
$\widetilde{\cal X}\setminus E$.
\par
By (iii), we have $\widetilde{\nu}=\widetilde{\mu}$
under the canonical isomorphism
${\Bbb P}^{1}(\Omega^{1}_{\cal X})=
{\Bbb P}(T{\cal X})^{\lor}$.
\par
We set $\widetilde{\pi}:=\pi\circ q$ and
$\widetilde{X}_{t}:=\widetilde{\pi}^{-1}(t)$
for $t\in S$. Similarly, we set
$E_{b}:=E\cap\widetilde{X}_{b}$ for $b\in{\cal D}$.
Then $E=\amalg_{b\in{\cal D}}E_{b}$, because
$E=q^{-1}(\Sigma_{\pi})\subset
\widetilde{\pi}^{-1}({\cal D})$.

\subsection
{}
{\bf The singularity of Quillen metrics}
\par
After Barlet \cite{Barlet82}, we define
a subspace of $C^{0}({\cal U})$ by
$$
{\cal B}({\cal U})
:=
C^{\infty}({\cal U})
\oplus
\bigoplus_{r\in{\Bbb Q}\cap(0,1]}
\bigoplus_{k=0}^{n}
|t|^{2r}(\log|t|)^{k}\cdot C^{\infty}({\cal U}).
$$
A function
$\varphi(t)\in{\cal B}({\cal U})$ has 
an asymptotic expansion at $0\in{\cal D}$,
i.e., there exist
$r_{1},\ldots,r_{m}\in{\Bbb Q}\cap(0,1]$ and
$f_{0},f_{l,k}\in C^{\infty}({\cal U})$,
$l=1,\ldots,m$, $k=0,\ldots,n$, such that
$\varphi(t)=f_{0}(t)+
\sum_{l=1}^{m}\sum_{k=0}^{n}
|t|^{2r_{l}}(\log|t|)^{k}\,f_{l,k}(t)$
as $t\to0$.
In what follows,
if $f(t),g(t)\in 
C^{\infty}({\cal U}^{o})$ satisfies
$f(t)-g(t)\in{\cal B}({\cal U})$, 
we write 
$$
f\equiv_{\cal B}g.
$$ 
\par
The purpose of Section 5 is to prove the following:

\begin{thm}
Let $\sigma_{p}$ be a nowhere vanishing $C^{\infty}$
section of the K\"ahler extension
$\lambda({\cal E}^{p}_{{\cal X}/S})|_{\cal U}$.
Then 
$$
\begin{aligned}
\,&
\log
\|\sigma_{p}\|^{2}_{\lambda({\cal E}^{p}_{{\cal X}/S}),
Q,g_{{\cal X}/S}}
\equiv_{\cal B}
\\
&
\left(
\int_{E_{0}}\sum_{j=0}^{p}(-1)^{p-j}\,
\widetilde{\mu}^{*}
\left\{
{\rm Td}(U)\,
\frac{{\rm Td}(c_{1}(H))-e^{-(p-j)c_{1}(H)}}
{c_{1}(H)}\right\}\,
q^{*}{\rm ch}(\Omega^{j}_{\cal X})
\right)\,\log|t|^{2}.
\end{aligned}
$$
\end{thm}

\par
The proof of Theorem 5.4 is divided into the following
three intermediary results, whose proofs shall be 
given in the subsections below:

\begin{prop}
The following identity of functions on $\cal U$
holds
$$
\log(
\|\cdot\|^{2}_{\lambda({\cal E}^{p}_{{\cal X}/S}),
Q,g_{{\cal X}/S}}/
\|\cdot\|^{2}_{\lambda({\cal E}^{p}_{{\cal X}/S}),Q,
g_{\pi^{*}\Omega^{1}_{S}}})\equiv_{\cal B}0.
$$
\end{prop}

\begin{prop}
The following identity of functions on $\cal U$
holds
$$
\begin{aligned}
\,
&
\log\left(
\frac{\|\cdot\|^{2}_{\lambda({\cal E}^{p}_{{\cal X}/S}),
Q,g_{\pi^{*}\Omega^{1}_{S}}}}
{\|\cdot\|^{2}_{\lambda({\cal E}^{p}_{{\cal X}/S}),
Q,\pi^{*}k_{S}}}\right)
\equiv_{\cal B}
\\
&
\left(\int_{E_{0}}\sum_{j=0}^{p}(-1)^{p-j}\,
\widetilde{\mu}^{*}\left\{
{\rm Td}(U)\,
\frac{1-e^{-(p-j)c_{1}(H)}}{c_{1}(H)}\right\}\,
q^{*}{\rm ch}(\Omega^{j}_{\cal X})\right)\,
\log|t|^{2}.
\end{aligned}
$$
\end{prop}

\begin{prop}
The following identity of functions on $\cal U$
holds
$$
\begin{aligned}
\,&
\log
\|\sigma_{p}(t)
\|^{2}_{\lambda({\cal E}^{p}_{{\cal X}/S}),
Q,\pi^{*}k_{S}}
\equiv_{\cal B}
\\
&
\left(
\int_{E_{0}}\sum_{j=0}^{p}(-1)^{p-j}
\widetilde{\mu}^{*}\left\{{\rm Td}(U)\,
\frac{{\rm Td}(c_{1}(H))-1}{c_{1}(H)}\right\}\,
q^{*}{\rm ch}(\Omega^{j}_{\cal X})
\right)\,
\log|t|^{2}.
\end{aligned}
$$
\end{prop}

{\bf Proof of Theorem 5.4.} 
By Propositions 5.5, 5.6, and 5.7, we get
$$
\begin{aligned}
\,&
\log
\|\sigma_{p}
\|^{2}_{\lambda({\cal E}^{p}_{{\cal X}/S}),Q,
g_{{\cal X}/S}}
=
\\
&
\log\left(
\frac{\|\cdot\|^{2}_{\lambda({\cal E}^{p}_{{\cal X}/S}),
Q,g_{{\cal X}/S}}}
{\|\cdot\|^{2}_{\lambda({\cal E}^{p}_{{\cal X}/S}),Q,
g_{\pi^{*}\Omega^{1}_{S}}}}\right)
+
\log\left(
\frac{\|\cdot\|^{2}_{\lambda({\cal E}^{p}_{{\cal X}/S}),
Q,g_{\pi^{*}\Omega^{1}_{S}}}}
{\|\cdot\|^{2}_{\lambda({\cal E}^{p}_{{\cal X}/S}),
Q,\pi^{*}k_{S}}}\right)
+
\log
\|\sigma_{p}
\|^{2}_{\lambda({\cal E}^{p}_{{\cal X}/S}),Q,\pi^{*}k_{S}}
\\
&
\equiv_{\cal B}
\left(\int_{E_{0}}\sum_{j=0}^{p}(-1)^{p-j}\,
\widetilde{\mu}^{*}\left\{{\rm Td}(U)\,
\frac{1-e^{-(p-j)c_{1}(H)}}{c_{1}(H)}\right\}\,
q^{*}{\rm ch}(\Omega^{j}_{\cal X})\right)\,
\log|t|^{2}
\\
&
\quad+
\left(
\int_{E_{0}}\sum_{j=0}^{p}(-1)^{p-j}
\widetilde{\mu}^{*}\left\{{\rm Td}(U)\,
\frac{{\rm Td}(c_{1}(H))-1}{c_{1}(H)}\right\}\,
q^{*}{\rm ch}(\Omega^{j}_{\cal X})
\right)\,
\log|t|^{2}
\\
&
\equiv_{\cal B}
\left(
\int_{E_{0}}\sum_{j=0}^{p}(-1)^{p-j}\,
\widetilde{\mu}^{*}\left\{
{\rm Td}(U)\,
\frac{{\rm Td}(c_{1}(H))-e^{-(p-j)c_{1}(H)}}
{c_{1}(H)}\right\}\,
q^{*}{\rm ch}(\Omega^{j}_{\cal X})
\right)\,
\log|t|^{2}.
\end{aligned}
$$
This proves the theorem.
$\Box$

\subsection
{}{\bf Proof of Proposition 5.5}
\par
Let
$g_{\Omega_{\cal X}^{i}\otimes
(\pi^{*}\Omega^{1}_{S})^{\otimes(p-i)}}$ 
be the Hermitian metric on 
$\Omega_{\cal X}^{i}\otimes
(\pi^{*}\Omega^{1}_{S})^{\otimes(p-i)}$ 
induced from $g_{\cal X}$, $g_{\pi^{*}\Omega^{1}_{S}}$.
We define
$\overline{\cal F}^{p}_{{\cal X}/S}$ to be
the complex of holomorphic vector bundles 
${\cal F}^{p}_{{\cal X}/S}$
equipped with the Hermitian metrics
$g_{\Omega_{\cal X}^{i}\otimes
(\pi^{*}\Omega^{1}_{S})^{\otimes(p-i)}}$ 
on $\Omega_{\cal X}^{i}\otimes
(\pi^{*}\Omega^{1}_{S})^{\otimes(p-i)}$
and $g_{\Omega^{p}_{{\cal X}/S}}$ 
on $\Omega^{p}_{{\cal X}/S}$.
\par
Let $\pi_{*}$ (resp. $\widetilde{\pi}_{*}$) be
the integration along the fibers of $\pi$
(resp. $\widetilde{\pi}$). 
For a $C^{\infty}$ differential form $\psi$ on 
$\widetilde{\cal X}$, one has
$\widetilde{\pi}_{*}(\psi)^{(0,0)}\in
{\cal B}({\cal U})$ by 
\cite[Th.\,4bis]{Barlet82}.
\par
Since ${\cal F}^{p}_{{\cal X}/S}$ is acyclic
on ${\cal X}^{o}$, 
the following identity of $C^{\infty}$ functions
on $S^{o}$ holds by the anomaly formula 
\cite[Th.\,0.3]{BGS88}:
\begin{equation}
\log\left(
\frac{\|\cdot\|^{2}_{\lambda(
{\cal E}^{p}_{{\cal X}/S}),Q,
g_{{\cal X}/S}}}
{\|\cdot\|^{2}_{\lambda({\cal E}^{p}_{{\cal X}/S}),Q,
g_{\pi^{*}\Omega^{1}_{S}}}}\right)
=
\pi_{*}\left({\rm Td}(T{\cal X}/S,g_{{\cal X}/S})\,
\widetilde{\rm ch}
(\overline{\cal F}^{p}_{{\cal X}/S})\right)^{(0,0)}.
\end{equation}
By (5.1), the following identity of $C^{\infty}$ 
differential forms on ${\cal X}\setminus\Sigma_{\pi}$ 
holds:
$$
{\rm Td}(T{\cal X}/S,g_{{\cal X}/S})\,
\widetilde{\rm ch}(\overline{\cal F}^{p}_{{\cal X}/S})
|_{{\cal X}\setminus\Sigma_{\pi}}=
\mu^{*}{\rm Td}(U,g_{U})\,
\nu^{*}\widetilde{\rm ch}(\overline{\cal K}^{p}).
$$
Since $q_{*}=(q^{-1})^{*}$ on
$\widetilde{\cal X}\setminus q^{-1}(\Sigma_{\pi})$, 
this yields the following identity 
on ${\cal X}\setminus\Sigma_{\pi}$:
$$
{\rm Td}(T{\cal X}/S,g_{{\cal X}/S})\,
\widetilde{\rm ch}(\overline{\cal F}^{p}_{{\cal X}/S})
|_{{\cal X}\setminus\Sigma_{\pi}}
=
(q)_{*}\left\{
\widetilde{\mu}^{*}{\rm Td}(U,g_{U})\,
\widetilde{\nu}^{*}\widetilde{\rm ch}(\overline{\cal K}^{p})
\right\}.
$$
Hence we get the following equation of 
$C^{\infty}$ functions on $S^{o}$:
\begin{equation}
\pi_{*}\left({\rm Td}(T{\cal X}/S,g_{{\cal X}/S})\,
\widetilde{\rm ch}(\overline{\cal F}^{p}_{\cal X})
|_{{\cal X}\setminus\Sigma_{\pi}}\right)^{(0,0)}
=
\left[\widetilde{\pi}_{*}\left\{
\widetilde{\mu}^{*}{\rm Td}(U,g_{U})\,
\widetilde{\nu}^{*}\widetilde{\rm ch}
(\overline{\cal K}^{p})\right\}\right]^{(0,0)}.
\end{equation}
Since 
$\{\widetilde{\mu}^{*}{\rm Td}(U,g_{U})\,
\widetilde{\nu}^{*}\widetilde{\rm ch}
(\overline{\cal K}^{p})\}^{(n,n)}$ 
is a $C^{\infty}$ $(n,n)$-form on $\widetilde{\cal X}$ 
and since the projection
$\widetilde{\pi}\colon\widetilde{\cal X}\to S$ 
is proper and holomorphic, 
the right hand side of (5.3) lies in
${\cal B}({\cal U})$ by \cite[Th.\,4bis]{Barlet82},
which, together with (5.2), (5.3), yields the result. 
$\Box$

\subsection
{}{\bf Proof of Proposition 5.6}
\par
For $0\leq i\leq p$, 
we deduce from the anomaly formula 
\cite[Th.\,0.3]{BGS88} that
\begin{equation}
\begin{aligned}
\,
&
\log\left(
\frac{\|\cdot\|^{2}_{\lambda(\Omega^{i}_{\cal X}
\otimes(\pi^{*}\Omega^{1}_{S})^{\otimes(p-i)}),Q,
g_{\pi^{*}\Omega^{1}_{S}}}}
{\|\cdot\|^{2}_{\lambda(\Omega^{i}_{\cal X}\otimes
(\pi^{*}\Omega^{1}_{S})^{\otimes(p-i)}),Q,\pi^{*}k_{S}}}
\right)
\\
&=
\pi_{*}\left({\rm Td}(T{\cal X}/S,g_{{\cal X}/S})\,
{\rm ch}(\Omega^{i}_{\cal X},g_{\cal X})\,
\widetilde{\rm ch}
((\pi^{*}\Omega^{1}_{S})^{\otimes(p-i)};
\,\pi^{*}k_{S},\,g_{\pi^{*}\Omega^{1}_{S}})
\right)^{(0,0)}
\\
&=
\pi_{*}\left({\rm Td}(T{\cal X}/S,g_{{\cal X}/S})\,
{\rm ch}(\Omega^{i}_{\cal X},g_{\cal X})\,
\widetilde{\rm ch}
((\pi^{*}\Omega^{1}_{S})^{\otimes(p-i)};
\,\pi^{*}k_{S},\,\|d\pi\|^{2}\pi^{*}k_{S})
\right)^{(0,0)}.
\end{aligned}
\end{equation}
Since 
$\nu^{*}c_{1}(L,g_{L})|_{{\cal X}\setminus\Sigma_{\pi}}
=-dd^{c}\log\|d\pi\|^{2}$ 
and $c_{1}(\Omega^{1}_{S},k_{S})=0$ on $\cal U$,
we deduce from (3.7) that
\begin{equation}
\begin{aligned}
\,&
\left.
\widetilde{\rm ch}
((\pi^{*}\Omega^{1}_{S})^{\otimes l};\,
\pi^{*}k_{S}^{\otimes l},\,
\|d\pi\|^{2l}\pi^{*}k_{S}^{\otimes l})
\right|_{\pi^{-1}({\cal U})\setminus\Sigma_{\pi}}
\\
&=
\sum_{m=1}^{\infty}\frac{1}{m!}\sum_{a+b=m-1}
c_{1}\left((\pi^{*}\Omega^{1}_{S})^{\otimes l},
\pi^{*}k_{S}^{\otimes l}\right)^{a}\,
c_{1}\left((\pi^{*}\Omega^{1}_{S})^{\otimes l},
\|d\pi\|^{2l}\pi^{*}k_{S}^{\otimes l}\right)^{b}\,
\log\|d\pi\|^{2l}
\\
&=
\sum_{m=1}^{\infty}\frac{1}{m!}
(-dd^{c}\log\|d\pi\|^{2l})^{m-1}\log\|d\pi\|^{2l}
=
\frac{e^{l\nu^{*}c_{1}(L,g_{L})}-1}
{\nu^{*}c_{1}(L,g_{L})}\,\log\|d\pi\|^{2}.
\end{aligned}
\end{equation}
By substituting (5.5) and 
${\rm Td}(T{\cal X}/S,g_{{\cal X}/S})=
\mu^{*}{\rm Td}(U,g_{U})$
into (5.4), we get
\begin{equation}
\begin{aligned}
\,
&
\left.\log\left(
\frac{\|\cdot\|^{2}_{\lambda(\Omega^{i}_{\cal X}\otimes
(\pi^{*}\Omega^{1}_{S})^{\otimes(p-i)}),Q,
g_{\pi^{*}\Omega^{1}_{S}}}}
{\|\cdot\|^{2}_{\lambda(\Omega^{i}_{\cal X}\otimes
(\pi^{*}\Omega^{1}_{S})^{\otimes(p-i)}),Q,\pi^{*}k_{S}}}
\right)\right|_{{\cal U}^{o}}
\\
&=
\pi_{*}\left\{
\mu^{*}{\rm Td}(U,g_{U})\,
{\rm ch}(\Omega^{i}_{\cal X},g_{\cal X})\,
\frac{e^{(p-i)\nu^{*}c_{1}(L,g_{L})}-1}
{\nu^{*}c_{1}(L,g_{L})}\,\log\|d\pi\|^{2}\right\}^{(0,0)}
\\
&=
\widetilde{\pi}_{*}\left\{
\widetilde{\mu}^{*}{\rm Td}(U,g_{U})\,
q^{*}{\rm ch}(\Omega^{i}_{\cal X},g_{\cal X})\,
\frac{e^{(p-i)\widetilde{\nu}^{*}c_{1}(L,g_{L})}-1}
{\widetilde{\nu}^{*}c_{1}(L,g_{L})}\,
q^{*}(\log\|d\pi\|^{2})\right\}^{(0,0)},
\end{aligned}
\end{equation}
which yields that
\begin{equation}
\begin{aligned}
\,
&
\left.\log\left(
\frac{\|\cdot\|^{2}_{\lambda({\cal E}^{p}_{{\cal X}/S}),
Q,g_{\pi^{*}\Omega^{1}_{S}}}}
{\|\cdot\|^{2}_{\lambda({\cal E}^{p}_{{\cal X}/S}),
Q,\pi^{*}k_{S}}}\right)\right|_{{\cal U}^{o}}
=
\\
&
\sum_{j=0}^{p}(-1)^{p-j}
\log\left(
\frac{\|\cdot\|^{2}_{\lambda(\Omega^{j}_{\cal X}
\otimes
(\pi^{*}\Omega^{1}_{S})^{\otimes(p-j)}),Q,
g_{\pi^{*}\Omega^{1}_{S}}}}
{\|\cdot\|^{2}_{\lambda(\Omega^{j}_{\cal X}\otimes
(\pi^{*}\Omega^{1}_{S})^{\otimes(p-j)}),Q,\pi^{*}k_{S}}}
\right)=
\\
&
\widetilde{\pi}_{*}\left[
q^{*}(\log\|d\pi\|^{2})\,\sum_{j=0}^{p}(-1)^{p-j}
\widetilde{\mu}^{*}{\rm Td}(U,g_{U})\,
\widetilde{\nu}^{*}\left\{
\frac{e^{(p-j)c_{1}(L,g_{L})}-1}{c_{1}(L,g_{L})}
\right\}
q^{*}{\rm ch}(\Omega^{j}_{\cal X},g_{\cal X})\,
\right]^{(0,0)}.
\end{aligned}
\end{equation}

\begin{lem}
Let $\varphi$ be a $\partial$ and 
$\bar{\partial}$-closed $C^{\infty}$ 
differential form on $\widetilde{\cal X}$. 
Let $(F,\|\cdot\|)$ be a holomorphic Hermitian 
line bundle on $\widetilde{\cal X}$. 
Let $s$ be a holomorphic section of $F$ with
${\rm div}(s)\subset\bigcup_{b\in{\cal D}}
\widetilde{X}_{b}$. Then the following identity
of functions on ${\cal U}$ holds 
$$
\widetilde{\pi}_{*}
((\log\|s\|^{2})\,\varphi)^{(0,0)}|_{\cal U}
\equiv_{\cal B}
\left(\int_{{\rm div}(s)\cap\widetilde{X}_{0}}
\varphi\right)\,\log|t|^{2}.
$$
In particular,
$$
\widetilde{\pi}_{*}
(q^{*}(\log\|d\pi\|^{2})\,\varphi)^{(0,0)}
|_{\cal U}
\equiv_{\cal B}
\left(\int_{E_{0}}\varphi\right)\,\log|t|^{2}.
$$
\end{lem}

\begin{pf}
See \cite[Lemma\,4.4 and Cor.\,4.6]{Yoshikawa05}
\end{pf}

Since
$\sum_{j=0}^{p}(-1)^{p-j}
\widetilde{\mu}^{*}{\rm Td}(U,g_{U})\,
\widetilde{\nu}^{*}\{
\frac{e^{(p-j)c_{1}(L,g_{L})}-1}{c_{1}(L,g_{L})}\}
q^{*}{\rm ch}(\Omega^{j}_{\cal X},g_{\cal X})$
is a $C^{\infty}$ differential form on $\widetilde{X}$
and since 
$\widetilde{\nu}^{*}c_{1}(L)=
-\widetilde{\mu}^{*}c_{1}(H)$ 
in $H^{2}(\widetilde{\pi}^{-1}({\cal U}),{\Bbb Z})$,
Proposition 5.6 follows from (5.7) and Lemma 5.8.
\qed

\subsection
{}{\bf Proof of Proposition 5.7}
\par
We need the following result:

\begin{thm}
Let $\xi\to{\cal X}$ be a holomorphic vector bundle 
on $\cal X$ equipped with a Hermitian metric $h_{\xi}$.
Let $\lambda(\xi)=\det R\pi_{*}\xi$ be the determinant
of the cohomologies of $\xi$ equipped with the Quillen
metric $\|\cdot\|^{2}_{\lambda(\xi),Q}$
with respect to $g_{{\cal X}/S}$ and $\xi$. 
Let $s$ be a nowhere vanishing holomorphic
section of $\lambda(\xi)|_{\cal U}$.
Then
$$
\log\|s\|^{2}_{Q,\lambda(\xi)}
\equiv_{\cal B}
\left(\int_{E_{0}}\widetilde{\mu}^{*}\left\{
{\rm Td}(U)\,
\frac{{\rm Td}(c_{1}(H))-1}{c_{1}(H)}\right\}\,
q^{*}{\rm ch}(\xi)\right)\log|t|^{2}.
$$
\end{thm}

\begin{pf}
See \cite[Th.\,1.1]{Yoshikawa05}.
\end{pf}

Let $\sigma_{(p,j)}$ be a nowhere vanishing 
$C^{\infty}$ section of 
$\lambda(\Omega^{j}_{\cal X}\otimes
(\pi^{*}\Omega^{1}_{S})^{\otimes(p-j)})|_{\cal U}$.
Then
$$
\sigma_{p}:=
\otimes_{j=0}^{p}\sigma_{(p,j)}^{(-1)^{p-j}}
$$
is a nowhere vanishing $C^{\infty}$ section of
$\lambda({\cal E}^{p}_{{\cal X}/S})|_{\cal U}$.
Since $\pi^{*}\Omega^{1}_{S}$ is trivial 
near $E_{0}$ and since
$$
\log
\|\cdot\|^{2}_{\lambda({\cal E}^{p}_{{\cal X}/S}),
Q,\pi^{*}k_{S}}=
\sum_{j=0}^{p}(-1)^{p-j}
\log
\|\cdot\|^{2}_{\lambda(\Omega^{j}_{\cal X}\otimes
(\pi^{*}\Omega^{1}_{S})^{\otimes(p-j)}),Q,\pi^{*}k_{S}},
$$
we deduce from Theorem 5.9 that
$$
\begin{aligned}
\,&
\log
\|\sigma_{p}\|^{2}_{
\lambda({\cal E}^{p}_{{\cal X}/S}),Q,\pi^{*}k_{S}}
|_{\cal U}
\\
&
=\sum_{j=0}^{p}(-1)^{p-j}
\log
\|\sigma_{(p,j)}\|^{2}_{
\lambda(\Omega^{j}_{\cal X}\otimes
(\pi^{*}\Omega^{1}_{S})^{\otimes(p-j)}),Q,\pi^{*}k_{S}}
|_{\cal U}
\\
&
\equiv_{\cal B}
\sum_{j=0}^{p}(-1)^{p-j}
\left(
\int_{E_{0}}\widetilde{\mu}^{*}\left\{
{\rm Td}(U)\,
\frac{{\rm Td}(c_{1}(H))-1}{c_{1}(H)}\right\}\,
q^{*}{\rm ch}(\Omega^{j}_{\cal X}\otimes
(\pi^{*}\Omega^{1}_{S})^{\otimes(p-j)})\right)\,
\log|t|^{2}
\\
&
\equiv_{\cal B}
\left(
\int_{E_{0}}\sum_{j=0}^{p}(-1)^{p-j}
\widetilde{\mu}^{*}\left\{{\rm Td}(U)\,
\frac{{\rm Td}(c_{1}(H))-1}{c_{1}(H)}\right\}\,
q^{*}{\rm ch}
(\Omega^{j}_{\cal X})\right)\,
\log|t|^{2}.
\end{aligned}
$$
This completes the proof of Proposition 5.7.
\qed

\subsection
{}
{\bf An extension of Theorem 5.4}
\par
Let $h_{\pi^{-1}({\cal U})}$ be a K\"ahler metric
on $\pi^{-1}({\cal U})$, and let $h_{{\cal X}/S}$
be the Hermitian metric on $T{\cal X}/S$ induced
from $h_{\pi^{-1}({\cal U})}$. We do not assume
that $h_{\pi^{-1}({\cal U})}$ extends to a K\"ahler
metric on ${\cal X}$.

\begin{thm}
Let $\sigma_{p}$ be a nowhere vanishing $C^{\infty}$
section of the K\"ahler extension
$\lambda({\cal E}^{p}_{{\cal X}/S})|_{\cal U}$.
Then 
$$
\begin{aligned}
\,&
\log
\|\sigma_{p}\|^{2}_{\lambda({\cal E}^{p}_{{\cal X}/S}),
Q,h_{{\cal X}/S}}|_{\cal U}
\equiv_{\cal B}
\\
&
\left(
\int_{E_{0}}\sum_{j=0}^{p}(-1)^{p-j}\,
\widetilde{\mu}^{*}
\left\{
{\rm Td}(U)\,
\frac{{\rm Td}(c_{1}(H))-e^{-(p-j)c_{1}(H)}}
{c_{1}(H)}\right\}\,
q^{*}{\rm ch}(\Omega^{j}_{\cal X})
\right)\,\log|t|^{2}.
\end{aligned}
$$
\end{thm}

\begin{pf}
By the anomaly formula \cite[Ths.\,0.2 and 0.3]{BGS88}, 
we have on ${\cal U}^{o}$
\begin{equation}
\begin{aligned}
\,
&
\log
\left(
\|\cdot\|^{2}_{
\lambda({\cal E}^{p}_{{\cal X}/S}),Q,h_{{\cal X}/S}}/
\|\cdot\|^{2}_{
\lambda({\cal E}^{p}_{{\cal X}/S}),Q,g_{{\cal X}/S}}
\right)
\\
&
=
\sum_{q}(-1)^{q}q\,
\pi_{*}\left(
\widetilde{\rm Td}(T{\cal X}/S;\,
g_{{\cal X}/S},h_{{\cal X}/S})\,
{\rm ch}
(\Omega^{q}_{{\cal X}/S},
h_{\Omega^{q}_{{\cal X}/S}})\right)^{(0,0)}
\\
&
\quad+
\sum_{q}(-1)^{q}q\,
\pi_{*}\left(
{\rm Td}(T{\cal X}/S,
g_{{\cal X}/S})\,
\widetilde{\rm ch}
(\Omega^{q}_{{\cal X}/S};\,
g_{\Omega^{q}_{{\cal X}/S}},
h_{\Omega^{q}_{{\cal X}/S}})
\right)^{(0,0)}.
\end{aligned}
\end{equation}
\par
Let $h_{U}$ be the Hermitian metric on $U$ induced 
from $(\varPi^{\lor})^{*}h_{\pi^{-1}(\cal U)}$.
Let $h_{\Omega^{1}_{\cal X}}$ be the Hermitian metric
on $\Omega^{1}_{\cal X}|_{\pi^{-1}({\cal U})}$
induced from $h_{\pi^{-1}(\cal U)}$.
Let $h_{\Omega^{q}_{{\cal X}/S}}$ be the Hermitian
metric on $\Omega^{q}_{{\cal X}/S}$ induced from
$h_{\Omega^{1}_{\cal X}}$.
Let $h_{\wedge^{q}Q}$ be the Hermitian metric
on $\wedge^{q}Q$ induced from 
$\varPi^{*}h_{\Omega^{1}_{\cal X}}$.
Then we have the following isomorphisms of
holomorphic Hermitian vector bundles over
${\cal X}\setminus\Sigma_{\pi}$:
\begin{equation}
(T{\cal X}/S,h_{{\cal X}/S})=\mu^{*}(U,h_{U}),
\qquad
(\Omega^{q}_{{\cal X}/S},h_{\Omega^{q}_{{\cal X}/S}})
=\nu^{*}(\wedge^{q}Q,h_{\wedge^{q}Q}).
\end{equation}
By (5.1), (5.8), (5.9),  we get
\begin{equation}
\begin{aligned}
\,
&
\log
\left(
\|\cdot\|^{2}_{
\lambda({\cal E}^{p}_{{\cal X}/S}),Q,h_{{\cal X}/S}}/
\|\cdot\|^{2}_{
\lambda({\cal E}^{p}_{{\cal X}/S}),Q,g_{{\cal X}/S}}
\right)
\\
&
=
\sum_{q}(-1)^{q}q\,
\widetilde{\pi}_{*}\left(
\widetilde{\mu}^{*}
\widetilde{\rm Td}(U;\,g_{U},h_{U})\,
\widetilde{\nu}^{*}{\rm ch}
(\wedge^{q}Q,h_{\wedge^{q}Q})\right)^{(0,0)}
\\
&
\quad+
\sum_{q}(-1)^{q}q\,
\widetilde{\pi}_{*}\left(
\widetilde{\mu}^{*}{\rm Td}(U,g_{U})\,
\widetilde{\nu}^{*}\widetilde{\rm ch}
(\wedge^{q}Q;\,g_{\wedge^{q}Q},h_{\wedge^{q}Q})
\right)^{(0,0)}
\equiv_{\cal B}0.
\end{aligned}
\end{equation}
Here the right hand side of (5.10) lies in
${\cal B}({\cal U})$ by \cite[Th.\,4bis]{Barlet82},
because
$$
\widetilde{\mu}^{*}
\widetilde{\rm Td}(U;h_{U},g_{U})\,
\widetilde{\nu}^{*}{\rm ch}
(\wedge^{q}Q,h_{\wedge^{q}Q}),
\qquad
\widetilde{\mu}^{*}{\rm Td}(U,g_{U})\,
\widetilde{\nu}^{*}\widetilde{\rm ch}
(\wedge^{q}Q;h_{\wedge^{q}Q},g_{\wedge^{q}Q})
$$
are $C^{\infty}$ differential forms on
$\widetilde{\pi}^{-1}({\cal U})$.
The result follows from Th.\,5.4 and (5.10).
\end{pf}

\subsection
{}
{\bf The case of ODP}
\par
In Subsection 5.9, 
we assume that $\Sigma_{\pi}\cap X_{0}$ consists of 
non-degenerate critical points. 
Hence ${\rm Sing}(X_{0})$ consists of ODP's. 
For $y\in{\cal X}$, let ${\frak m}_{y}$ be
the maximal ideal of the local ring 
${\cal O}_{{\cal X},y}$. Then there exists 
a neighborhood of $X_{0}$ in $\cal X$ on which
${\cal I}_{\Sigma_{\pi}}=\oplus_{y\in{\rm Sing}(X_{0})}
{\frak m}_{y}$. 
Let $q\colon\widetilde{\cal X}\to{\cal X}$
be the blowing-up of the discrete set
$\Sigma_{\pi}\cap X_{0}$, and
set $E_{y}:=q^{-1}(y)$ for $y\in{\rm Sing}(X_{0})$.
Then $E_{0}=\amalg_{y\in{\rm Sing}(X_{0})}E_{y}$
and $E_{y}\cong{\Bbb P}^{n}$.
\par
Since $\Sigma_{\pi}$ is discrete, we may identify
${\Bbb P}(\Omega^{1}_{\cal X})$ and
${\Bbb P}(T{\cal X})$ with the trivial projective
bundle on a neighborhood of $\Sigma_{\pi}\cap X_{0}$
by fixing a system of coordinates near 
$\Sigma_{\pi}\cap X_{0}$. Under this trivialization,
we consider the Gauss maps $\nu$ and $\mu$ 
only on a small neighborhood of $\Sigma_{\pi}\cap X_{0}$. 
Then we have the following on a neighborhood of
each $y\in\Sigma_{\pi}\cap X_{0}$:
$$
\mu(z)=\nu(z)=
\left(\frac{\partial\pi}{\partial z_{0}}(z):\cdots:
\frac{\partial\pi}{\partial z_{n}}(z)\right).
$$
Since $\pi$ is non-degenerate at every 
$y\in\Sigma_{\pi}\cap X_{0}$, 
we may assume by Morse's lemma that
$\pi(z)=z_{0}^{2}+\cdots+z_{n}^{2}$ 
near $\Sigma_{\pi}\cap X_{0}$. Hence the composition
$\nu\circ q\colon\widetilde{\cal X}\setminus E_{0}\to
{\Bbb P}^{n}$ extends to a holomorphic map
$\widetilde{\nu}:=\nu\circ q\colon\widetilde{\cal X}\to
{\Bbb P}^{n}$ such that
$$
\widetilde{\nu}|_{E}=\widetilde{\mu}|_{E}=
{\rm id}_{E}.
$$ 
\par
For $n\in{\Bbb N}$ and $0\leq p\leq n$, 
set
$$
\delta(n,p):=
\sum_{j=0}^{p}(-1)^{j}
\binom{n+1}{j}\frac{(p-j+1)^{n+2}-(p-j)^{n+2}}{(n+2)!}.
$$
\par
For a formal power series $f(x)\in{\Bbb C}[[x]]$,
we define $f(x)|_{x^{m}}$ to be the coefficient of 
$x^{m}$ of $f(x)$. 
Recall that the metric 
$h_{{\cal X}/S}$ is defined only on $T{\cal X}/S
|_{{\pi}^{-1}({\cal U})\setminus\Sigma_{\pi}}$.

\begin{thm}
Let $\sigma_{p}$ be a nowhere vanishing $C^{\infty}$
section of the K\"ahler extension
$\lambda({\cal E}^{p}_{{\cal X}/S})|_{\cal U}$.
Then the following identity of functions on
$\cal U$ holds
$$
(-1)^{p}\log
\|\sigma_{p}(t)\|^{2}_{
\lambda({\cal E}^{p}_{{\cal X}/S}),Q,h_{{\cal X}/S}}
\equiv_{\cal B}
(-1)^{n}\delta(n,p)\,\#{\rm Sing}(X_{0})\,
\log|t|^{2}.
$$
\end{thm}

\begin{pf}
In Theorem 5.10, we can identify $U$ (resp. $L$)
with the universal hyperplane bundle 
(resp. tautological line bundle) on ${\Bbb P}^{n}$.
Then $H=L^{-1}$. Set $x:=c_{1}(H)$. Hence
$\int_{{\Bbb P}^{n}}x^{n}=1$.
From the exact sequence 
$0\to U\to{\Bbb C}^{n+1}\to H\to 0$, 
we get
${\rm Td}(U)={\rm Td}^{-1}(x)=
(1-e^{-x})/x$.
Since $q(E_{0})$ consists of a point,
we get
$q^{*}\Omega^{j}_{\cal X}|_{E_{0}}=
{\Bbb C}^{\oplus{\binom{n+1}{p}}}$.
By substituting this and the equation
$q^{*}{\rm ch}(\Omega^{j}_{\cal X})|_{E_{0}}
=\binom{n+1}{p}$
into the formula in Theorem 5.10, we get
\begin{equation}
\begin{aligned}
\,&
\int_{E}\sum_{j=0}^{p}(-1)^{p-j}\,
\widetilde{\mu}^{*}\left\{
{\rm Td}(U)\,
\frac{{\rm Td}(c_{1}(H))-e^{-(p-j)c_{1}(H)}}
{c_{1}(H)}\right\}\,
q^{*}{\rm ch}(\Omega^{j}_{\cal X})
\\
&=
\#{\rm Sing}(X_{0})\,
\left.\sum_{j=0}^{p}(-1)^{p-j}\,
\frac{1}{{\rm Td}(x)}\cdot
\frac{{\rm Td}(x)-e^{-(p-j)x}}{x}
\cdot\binom{n+1}{j}\right|_{x^{n}}
\\
&=
\#{\rm Sing}(X_{0})\,
\left.\sum_{j=0}^{p}(-1)^{p-j}\,\binom{n+1}{j}
\left\{\frac{(e^{-x}-1)e^{-(p-j)x}}{x^{2}}+
\frac{1}{x}\right\}\right|_{x^{n}}
\\
&=
\#{\rm Sing}(X_{0})\,
\sum_{j=0}^{p}(-1)^{p-j}\,\binom{n+1}{j}
\{e^{-(p-j+1)x}-e^{-(p-j)x}\}|_{x^{n+2}}
\\
&=
(-1)^{n-p}\,\delta(n,p)\,\#{\rm Sing}(X_{0}).
\end{aligned}
\end{equation}
The result follows from Theorem 5.4 and (5.11).
\end{pf}

\begin{lem}
The following identities hold:
$$
\delta(3,p)+\delta(3,3-p)=1
\quad(0\leq p\leq3),
\qquad
\sum_{p=0}^{3}p\,\delta(3,p)=
\frac{19}{4}.
$$
\end{lem}

\begin{pf}
By the definition of $\delta(n,p)$, we get
$$
\delta(3,0)=\frac{1}{120},
\quad
\delta(3,1)=\frac{27}{120},
\quad
\delta(3,2)=\frac{93}{120},
\quad
\delta(3,3)=\frac{119}{120},
$$
which yields the result.
\end{pf}

\par
Set
$$
\sigma:=\otimes_{p=0}^{n}\sigma_{p}^{(-1)^{p}p}.
$$
Then $\sigma$ is a nowhere vanishing $C^{\infty}$
section of $\lambda(\Omega_{{\cal X}/S}^{\bullet})$ 
near ${\cal D}$.

\begin{thm}
When $n=3$,
$$
\log
\|\sigma(t)\|^{2}_{
\lambda(\Omega^{\bullet}_{{\cal X}/S}),
Q,h_{{\cal X}/S}}
\equiv_{\cal B}
-\frac{19}{4}\,\#{\rm Sing}(X_{0})\,\log|t|^{2}.
$$
\end{thm}

\begin{pf}
By Theorem 5.11, we get
$$
\begin{aligned}
\log
\|\sigma\|^{2}_{\lambda(\Omega^{\bullet}_{{\cal X}/S}),
Q,g_{{\cal X}/S}}|_{\cal U}
&
=
\sum_{p=0}^{3}(-1)^{p}p\,\log
\|\sigma_{p}\|^{2}_{\lambda({\cal E}^{p}_{{\cal X}/S}),
Q,g_{{\cal X}/S}}|_{\cal U}
\\
&
\equiv_{\cal B}
(-1)^{3}\sum_{p=0}^{3}p\,\delta(3,p)\,
\#{\rm Sing}(X_{0})\,\log|t|^{2}.
\end{aligned}
$$
This, together with the second identity of Lemma 5.12, 
yields the result.
\end{pf}

\begin{rem}
In our subsequent paper \cite{FLY05},
we shall determine the behavior of
$\log
\|\sigma(t)\|^{2}_{\lambda(\Omega^{\bullet}_{{\cal X}/S}),
Q,h_{{\cal X}/S}}$ as $t\to0$ for arbitrary relative
dimension $n$.
\end{rem}

%%%%%%%%%%%%%%%%%%%%%%%%%%%%%%%%%%%%%%%%%%%%%%%%%%%%%%%%%%%
%
%                 Section 6
%
%%%%%%%%%%%%%%%%%%%%%%%%%%%%%%%%%%%%%%%%%%%%%%%%%%%%%%%%%%%

\section
{\bf The cotangent sheaf of the Kuranishi space}
\par
Let $X$ be a smoothable Calabi-Yau $n$-fold 
with only one ODP as its singular set. Let
${\frak p}\colon({\frak X},X)\to({\rm Def}(X),[X])$
be the Kuranishi family of $X$ with discriminant
locus $\frak D$. Then $\frak X$, ${\rm Def}(X)$, 
and $\frak D$ are smooth by Lemmas 2.3 and 2.7.

\begin{lem}
The dualizing sheaf $K_{\frak X}$ of $\frak X$ 
is trivial.
In particular, the relative dualizing sheaf
$K_{{\frak X}/{\rm Def}(X)}=
K_{\frak X}\otimes
(\frak p^{*}K_{{\rm Def}(X)})^{-1}$ 
is trivial.
\end{lem}

\begin{pf}
By the same argument as in 
\cite[p.68 l.25-l.28]{Yoshikawa04}, we see that 
$K_{\frak X}|_{X_{s}}\cong{\cal O}_{X_{s}}$ for
all $s\in{\rm Def}(X)$. 
Since ${\rm Def}(X)\cong\varDelta^{N+1}$, we get
the triviality of $K_{\frak X}$ by 
the same argument as in
\cite[p.68 l.29-l.33]{Yoshikawa04}.
\end{pf}

\par
Recall that the Kodaira-Spencer isomorphism
$$
\rho_{{\rm Def}(X)\setminus{\frak D}}\colon
\Theta_{{\rm Def}(X)\setminus{\frak D}}\to
R^{1}{\frak p}_{*}\Theta_{{\frak X}/{\rm Def}(X)}
|_{{\rm Def}(X)\setminus{\frak D}}
$$
was defined in Subsection 4.2.
By considering the dual of 
$\rho_{{\rm Def}(X)\setminus{\frak D}}$, 
the relative Serre duality induces an isomorphism of
${\cal O}_{{\rm Def}(X)}$-modules on
${\rm Def}(X)\setminus{\frak D}$:
$$
\rho_{{\rm Def}(X)\setminus{\frak D}}^{\lor}\colon
R^{n-1}{\frak p}_{*}(\Omega^{1}_{{\frak X}/{\rm Def}(X)}
\otimes K_{{\frak X}/{\rm Def}(X)})
|_{{\rm Def}(X)\setminus{\frak D}}
\cong
\Omega^{1}_{{\rm Def}(X)}
|_{{\rm Def}(X)\setminus{\frak D}}.
$$

\begin{thm}
The isomomorphism 
$\rho_{{\rm Def}(X)\setminus{\frak D}}^{\lor}$ 
extends to an isomorphism 
$$
\rho_{{\rm Def}(X)}^{\lor}\colon
R^{n-1}{\frak p}_{*}(\Omega^{1}_{{\frak X}/{\rm Def}(X)}
\otimes K_{{\frak X}/{\rm Def}(X)})
\cong
\Omega^{1}_{{\rm Def}(X)}.
$$
of ${\cal O}_{{\rm Def}(X)}$-modules over ${\rm Def}(X)$.
\end{thm}

The isomorphism $\rho_{{\rm Def}(X)}^{\lor}$ is again
called the Kodaira-Spencer isomorphism.
Before proving Theorem 6.2, we first prove 
an intermediate result in the next subsection.

\subsection
{}{\bf Blowing-up and the regularity of differential forms}
\par
Set
$$
\widetilde{\varDelta^{n+1}}=
\{(z,[\zeta])\in\varDelta^{n+1}\times{\Bbb P}^{n};\,
z_{i}\zeta_{j}-z_{j}\zeta_{i}=0
\quad 0\leq i,j\leq n\},
\qquad
q:={\rm pr}_{1}.
$$
Then 
$q\colon\widetilde{\varDelta^{n+1}}\to\varDelta^{n+1}$ 
is the blowing-up at the origin. 
Set $E:=q^{-1}(0)$ and
\begin{align*}
U_{i}:
&=
\{(z,[\zeta])\in\widetilde{\varDelta^{n+1}};\,
\zeta_{i}\not=0\},
\qquad
O_{i}:=\{z\in\varDelta^{n+1};\,
z_{i}\not=0\},
\\
W_{i}:
&=
\{(\zeta_{0},\ldots,\zeta_{i-1},z_{i},
\zeta_{i+1},\ldots,\zeta_{n})\in{\Bbb C}^{n+1};\,
|z_{i}|<1,\quad
|z_{i}\zeta_{j}|<1\quad(j\not=i)\}.
\end{align*}
Then $U_{i}\cong W_{i}\subset{\Bbb C}^{n+1}$ 
via the map
$$
\begin{aligned}
\,&
W_{i}\ni(\zeta_{0},\ldots,\zeta_{i-1},z_{i},
\zeta_{i+1},\ldots,\zeta_{n})
\\
&\to
((z_{i}\zeta_{0},\ldots,z_{i}\zeta_{i-1},z_{i},
z_{i}\zeta_{i+1},\ldots,z_{i}\zeta_{n}),
[\zeta_{0}:\cdots:\zeta_{i-1}:1:\zeta_{i+1}:
\cdots:\zeta_{n}])
\in U_{i}.
\end{aligned}
$$
By construction, we have
$\widetilde{\varDelta}^{n+1}=\bigcup_{i=0}^{n}U_{i}$
and
$$
E\cap U_{i}\cong
\{(\zeta_{0},\ldots,\zeta_{i-1},z_{i},
\zeta_{i+1},\ldots,\zeta_{n})\in W_{i};
\,z_{i}=0\},
\qquad
q(U_{i})\supset O_{i}.
$$
Let $\omega_{ij}$ be the $C^{\infty}$ $(n,0)$-form 
on $O_{i}$ defined by
$$
\omega_{ij}:=
\frac{|z_{j}|^{2}}{|z_{0}|^{2}+\cdots+|z_{n}|^{2}}\cdot
\frac{dz_{0}\wedge\cdots\wedge dz_{i-1}\wedge 
dz_{i+1}\wedge\cdots\wedge dz_{n}}
{z_{i}^{n-2}z_{j}}.
$$

\begin{lem} 
For all $0\leq i,j\leq n$,
the $C^{\infty}$ $(n,0)$-form $q^{*}\omega_{ij}$ 
on $q^{-1}(O_{i})=U_{i}\setminus E$ extends to 
a $C^{\infty}$ $(n,0)$-form on $U_{i}$ and satisfies 
$q^{*}\omega_{ij}|_{E\cap U_{i}}=0$.
\end{lem}

\begin{pf}
Since
$
q|_{W_{i}}(\zeta_{0},\ldots,\zeta_{i-1},z_{i},
\zeta_{i+1},\ldots,\zeta_{n})
=
(z_{i}\zeta_{0},\ldots,z_{i}\zeta_{i-1},z_{i},
z_{i}\zeta_{i+1},\ldots,z_{i}\zeta_{n})
$
under the identification $U_{i}\cong W_{i}$, 
we get the following two formulas:
$$
q^{*}\left(
\frac{|z_{j}|^{2}}{|z_{0}|^{2}+\cdots+|z_{n}|^{2}}
\right)
=\left\{
\begin{array}{ll}
|\zeta_{j}|^{2}(1+\|\zeta\|^{2})^{-1}
&(j\not=i)
\\
(1+\|\zeta\|^{2})^{-1}
&(j=i),
\end{array}
\right.
$$
$$
\begin{aligned}
\,&
q^{*}\left(z_{i}^{-(n-1)}
dz_{0}\wedge\cdots\wedge dz_{i-1}\wedge dz_{i+1}
\wedge\cdots\wedge dz_{n}\right)
\\
&=
z_{i}^{-(n-1)}
d(z_{i}\zeta_{0})\wedge\cdots\wedge d(z_{i}\zeta_{i-1})
\wedge d(z_{i}\zeta_{i+1})\wedge\cdots\wedge 
d(z_{i}\zeta_{n})
\\
&=
z_{i}\,d\zeta_{0}\wedge\cdots\widehat{d\zeta_{i}}
\cdots\wedge d\zeta_{n}
+dz_{i}\wedge
\sum_{j<i}(-1)^{j-1}
d\zeta_{0}\wedge\cdots\widehat{d\zeta_{j}}
\cdots\widehat{d\zeta_{i}}\cdots\wedge d\zeta_{n}
\\
&\quad+
dz_{i}\wedge
\sum_{j>i}(-1)^{j}
d\zeta_{0}\wedge\cdots\widehat{d\zeta_{i}}
\cdots\widehat{d\zeta_{j}}\cdots\wedge d\zeta_{n}
\in A^{n,0}(U_{i}),
\end{aligned}
$$
which yields that 
$q^{*}\omega_{ii}\in A^{n,0}(U_{i})$
and $q^{*}\omega_{ii}|_{E\cap U_{i}}=0$.
Since
$q^{*}\omega_{ij}=
\frac{\bar{\zeta}_{j}}{1+\|\zeta\|^{2}}\,
q^{*}\omega_{ii}$ when $j\not=i$, 
the assertion for $q^{*}\omega_{ij}$ $(i\not=j)$ 
follows from the assertion for $q^{*}\omega_{ii}$.
\end{pf}

\subsection
{}
{\bf Proof of Theorem 6.2}
\par
For simplicity, we set
$$
{\cal X}:={\frak X}, 
\quad
S:={\rm Def}(X),
\quad 
\pi:={\frak p}, 
\quad
0:=[X], 
\quad
X_{0}:=X, 
\quad
N+1=\dim S.
$$
Hence $(S,0)\cong(\varDelta^{N+1},0)$ and
$\pi\colon({\cal X},X_{0})\to(S,0)$
is the Kuranishi family of $X_{0}$.
\par
Let $s=(s_{0},\ldots,s_{N})$ be a system of coordinates 
of $S$ such that ${\frak D}={\rm div}(s_{0})$.
We set $s'=(s_{1},\ldots,s_{N})$.
Then $\partial/\partial s_{\alpha}$ is a nowhere vanishing 
holomorphic vector field on $S$ for $0\leq\alpha\leq N$.
\newline
{\bf (Step 1)}
The Kodaira-Spencer isomorphism 
$\rho_{S\setminus{\frak D}}\colon
\Theta_{S\setminus{\frak D}}\to
R^{1}\pi_{*}\Theta_{{\cal X}/S}|_{S\setminus{\frak D}}$
yields holomorphic sections
$\rho(\partial/\partial s_{\alpha})\in 
H^{0}(S\setminus{\frak D},
R^{1}\pi_{*}\Theta_{{\cal X}/S})$.
Let $\langle\cdot,\rangle_{s}$ be the Yoneda product
between 
$H^{n-1}(X_{s},
\Omega^{1}_{X_{s}}\otimes K_{X_{s}})$
and
${\rm Ext}^{1}_{{\cal O}_{X_{s}}}
(\Omega^{1}_{X_{s}}\otimes K_{X_{s}},
K_{X_{s}})$.
\par
Since $h^{n-1}(X_{s},\Omega^{1}_{X_{s}})=N+1$, 
there exist 
$\phi_{0},\ldots,\phi_{N}\in 
H^{n-1}({\cal X},
\Omega^{1}_{{\cal X}/S}\otimes K_{{\cal X}/S})$ 
such that
\newline{(i) }
$\{\phi_{0},\ldots,\phi_{N}\}$ is a basis of
$R^{n-1}\pi_{*}
(\Omega^{1}_{{\cal X}/S}\otimes K_{{\cal X}/S})$
as a free ${\cal O}_{S}$-module; 
\newline{(ii) }
$\{\phi_{0}|_{X_{s}},\ldots,\phi_{N}|_{X_{s}}\}$ 
is a basis of
$H^{n-1}(X_{s},
\Omega^{1}_{X_{s}}\otimes K_{X_{s}})$ 
for all $s\in S$;
\newline{(iii) }
$\langle\phi_{\alpha}|_{X_{0}},
\rho_{0}(\partial/\partial s_{\beta})\rangle_{0}
=\delta_{\alpha\beta}$ for $0\leq\alpha,\beta\leq N$.
\par
Let $\rho_{s}^{\lor}\colon
H^{n-1}(X_{s},\Omega^{1}_{X_{s}}\otimes K_{X_{s}})
\to\Omega^{1}_{S,s}$ be the dual of 
the Kodaira-Spencer map. For $s\in S$, set
$$
g_{\alpha\beta}(s):=
\langle\phi_{\alpha}|_{X_{s}},
\rho_{s}(\partial/\partial s_{\beta})\rangle_{s}=
\langle\langle\rho_{s}^{\lor}(\phi_{\alpha}|_{X_{s}}),
\partial/\partial s_{\beta}\rangle\rangle,
$$
where $\langle\langle\cdot,\cdot\rangle\rangle\colon
\Omega^{1}_{S,s}\times TS_{s}\to{\Bbb C}$ is 
the natural pairing.
Then $g_{\alpha\beta}$ is a function on $S$,
which is holomorphic on $S\setminus{\frak D}$ but
which may not be continuous on $S$, 
such that 
$$
g_{\alpha\beta}(0)=\delta_{\alpha\beta}.
$$ 
It suffices to prove $g_{\alpha\beta}\in C^{0}(S)$; 
if it is the case, $(g_{\alpha\beta}(s))$ is a family
of invertible matrices depending holomorphically 
on $s\in S$, so that
$R^{n-1}\pi_{*}
(\Omega^{1}_{{\cal X}/S}\otimes K_{{\cal X}/S})$
is the holomorphic dual bundle of $\Theta_{S}$ via 
the extension of $\rho_{S\setminus{\frak D}}^{\lor}$.
\newline
{\bf (Step 2)}
Let ${\cal A}_{{\cal X}}$ be the sheaf of germs of 
$C^{\infty}$ functions on $\cal X$, 
and let ${\cal A}^{p,q}_{{\cal X}}$ be the sheaf 
of germs of $C^{\infty}$ $(p,q)$-forms on $\cal X$.
Set
$$
A^{p,q}({\cal X},
\Omega^{1}_{{\cal X}/S}\otimes K_{{\cal X}/S})
:=
\Gamma({\cal X},
{\cal A}^{p,q}_{{\cal X}}\otimes_{{\cal O}_{{\cal X}}}
\Omega^{1}_{{\cal X}/S}\otimes K_{{\cal X}/S}).
$$
Then $A^{p,q}({\cal X},
\Omega^{1}_{{\cal X}/S}\otimes K_{{\cal X}/S})$
is the vector space of $C^{\infty}$ $(p,q)$-forms 
on $\cal X$ with values in 
$\Omega^{1}_{{\cal X}/S}\otimes K_{{\cal X}/S}$.
By Malgrange 
\cite[pp.\,88, Cor.\,1.12]{Malgrange66}, 
${\cal O}_{{\cal X}}$ is a flat 
${\cal A}_{\cal X}$-module. 
Hence we have the Dolbeault isomorphism 
\cite[Chap.\,VII, Prop.\,4.5]{Banica76}
$$
\begin{aligned}
\,
&
H^{n-1}({\cal X},
\Omega^{1}_{{\cal X}/S}\otimes K_{{\cal X}/S})
\\
&
=
\frac{\ker\{\bar{\partial}\colon
A^{0,n-1}({\cal X},
\Omega^{1}_{{\cal X}/S}\otimes K_{{\cal X}/S})
\to
A^{0,n}({\cal X},
\Omega^{1}_{{\cal X}/S}\otimes K_{{\cal X}/S})\}}
{{\rm Im}\{\bar{\partial}\colon
A^{0,n-2}({\cal X},
\Omega^{1}_{{\cal X}/S}\otimes K_{{\cal X}/S})
\to
A^{0,n-1}({\cal X},
\Omega^{1}_{{\cal X}/S}\otimes K_{{\cal X}/S})\}}.
\end{aligned}
$$
\par
Let
$\Phi_{\alpha}\in A^{0,n-1}({\cal X},
\Omega^{1}_{{\cal X}/S}\otimes K_{{\cal X}/S})$
be a $\bar{\partial}$-closed differential form 
representing $\phi_{\alpha}$, i.e., 
$\phi_{\alpha}=[\Phi_{\alpha}]$.
\newline
{\bf (Step 3)}
To study the behavior of $g_{\alpha\beta}(s)$ near 
$\frak D$,
we compute a representative of the Kodaira-Spencer
classes $\rho(\partial/\partial s_{\alpha})$
in the Dolbeault cohomology.
\par
Near the critical locus $\Sigma_{\pi}\subset{\cal X}$,
there is a neighborhood 
$V\cong\varDelta^{n+1}\times\varDelta^{N}$ 
of $\Sigma_{\pi}$ in $\cal X$ such that 
$\pi(z_{0},\ldots,z_{n},s')=
(z_{0}^{2}+\cdots+z_{n}^{2},s_{1},\ldots,s_{N})$.
Hence we have 
$\Sigma_{\pi}\cap V=\{0\}\times\varDelta^{N}$.
For $i=0,1,\ldots,n$, we set
$$
V_{i}:=
\varDelta^{i-1}\times\varDelta^{*}\times\varDelta^{n-i}
\times\varDelta^{N}
=
\{(z,s')\in\varDelta^{n+1}\times\varDelta^{N};
\,z_{i}\not=0\}.
$$
Then $\{V_{i}\}_{i}$ is an open covering of
$V\setminus\Sigma_{\pi}$, i.e.,
$V\setminus\Sigma_{\pi}=
\bigcup_{i=0}^{n}V_{i}$.
Let $\{V_{\lambda}\}_{\lambda\in\Lambda}$ be 
an open covering of $\overline{{\cal X}\setminus V}$
such that 
$V_{\lambda}\cong\varDelta^{n}\times\varDelta^{N+1}$
and $\pi|_{V_{\lambda}}={\rm pr}_{2}$. 
Then
${\frak V}:=\{V_{i}\}_{i}\cup\{V_{\lambda}\}_{\lambda}$ 
is an open covering of ${\cal X}\setminus\Sigma_{\pi}$.
\par
First let us construct a representative 
of the Kodaira-Spencer class 
$\rho(\partial/\partial s_{\alpha})$ in the Cech 
cohomology with respect to the covering $\frak V$.
\par
On $V_{i}$, set
$$
v^{(i)}_{0}:=\frac{1}{2z_{i}}
\frac{\partial}{\partial z_{i}},
\qquad
v^{(i)}_{\alpha}=\frac{\partial}{\partial s_{\alpha}}
\quad
(\alpha=1,\ldots,N).
$$
Then
$v^{(i)}_{0},\ldots v^{(i)}_{N}\in
H^{0}(V_{i},\Theta_{\cal X})$ and
$\pi_{*}(v^{(i)}_{\alpha})=
\frac{\partial}{\partial s_{\alpha}}$
$(\alpha=0,\ldots,N)$.
We also fix a holomorphic vector field 
$v^{(\lambda)}_{\alpha}$ such that
$v^{(\lambda)}_{\alpha}=\partial/\partial s_{\alpha}$ 
on every $V_{\lambda}$. We get in Cech cohomology
$$
\rho\left(\frac{\partial}{\partial s_{\alpha}}\right)
=
\{(v^{(\mu)}_{\alpha}-v^{(\nu)}_{\alpha})
|_{V_{\mu}\cap V_{\nu}}\}_{V_{\mu},V_{\nu}\in{\frak V}}
\in
H^{1}({\cal X}\setminus\Sigma_{\pi},
\Theta_{{\cal X}/S};{\frak V}).
$$
\par
Let 
$\{\chi_{i}\}_{i}\cup\{\chi_{\lambda}\}_{\lambda}$
be a partition of unity of 
${\cal X}\setminus\Sigma_{\pi}$ subject to 
the covering $\frak V$ such that on $V_{i}$,
$$
\chi_{i}(z)=
\frac{|z_{i}|^{2}}{|z_{0}|^{2}+\cdots+|z_{n}|^{2}},
\qquad
i=0,\ldots,n.
$$
Then the following differential form
$\xi_{\alpha}\in A^{0,1}({\cal X}\setminus\Sigma_{\pi},
\Theta_{{\cal X}/S})$ represents 
$\rho(\partial/\partial s_{\alpha})$:
$$
\xi_{\alpha}|_{V_{\nu}}:=
\sum_{\mu}\bar{\partial}\chi_{\mu}\otimes 
(v^{(\mu)}_{\alpha}-v^{(\nu)}_{\alpha}),
\qquad
\rho_{s}\left(\frac{\partial}{\partial s_{\alpha}}\right)
=[\xi_{\alpha}|_{X_{s}}]
\quad(s\in S\setminus{\frak D}).
$$
In particular, we get on $V\setminus\Sigma_{\pi}$
$$
\xi_{0}|_{V\setminus\Sigma_{\pi}}=
\sum_{i=0}^{n}\bar{\partial}\chi_{i}\otimes
\frac{1}{2z_{i}}\frac{\partial}{\partial z_{i}},
\qquad
\xi_{\alpha}|_{V\setminus\Sigma_{\pi}}=0
\quad(\alpha=1,\ldots,N).
$$
\newline
{\bf (Step 4)}
Let us study the behavior of 
$g_{\alpha\beta}|_{S\setminus{\frak D}}(s)$
as $s\to{\frak D}$. 
Let
$\varrho(z)\in C^{\infty}_{0}(\varDelta^{n+1})$ 
be a cut-off function with $\varrho\equiv1$ near 
$0\in\varDelta^{n+1}$.
Recall that $\iota(\cdot)$ denotes the interior product.
There exists 
$h_{\alpha\beta}(s)\in C^{\infty}(S)$ such that
for $s\in S\setminus{\frak D}$,
$$
g_{\alpha\beta}(s)=
\langle\phi_{\alpha}|_{X_{s}},
\rho_{s}(\partial/\partial s_{\beta})\rangle_{s}
=
\int_{X_{s}}
\iota(\xi_{\beta})\Phi_{\alpha}
=
\int_{X_{s}\cap V}
\varrho(z)\cdot\iota(\xi_{\beta})\Phi_{\alpha}+
h_{\alpha\beta}(s).
$$
Since $\xi_{\beta}\equiv0$ on $V\setminus\Sigma_{\pi}$
for $\beta\not=0$,
$g_{\alpha\beta}|_{S\setminus{\frak D}}(s)$ extends to 
a $C^{\infty}$ function on $S$ if $\beta\not=0$. 
Let us prove that $g_{\alpha0}|_{S\setminus{\frak D}}$ 
extends to a continuous function on $S$.
\par
Since $\Phi_{\alpha}$ is a $(0,n-1)$-form on $\cal X$
with values in 
$\Omega^{1}_{{\cal X}/S}\otimes K_{{\cal X}/S}$, 
we can write
$$
\Phi_{\alpha}|_{V}=
\sum_{i=0}^{n}\theta_{\alpha}^{i}(z,s)\,[dz_{i}]
\otimes\eta,
$$
with 
$[dz_{i}]=dz_{i}\mod(\pi^{*}ds_{0},\ldots,\pi^{*}ds_{N})$, 
$\theta_{\alpha}^{i}\in A^{0,n-1}(V)$, and
$$
\eta|_{V_{i}}
:=
(-1)^{i-1}
\left.\frac{dz_{0}\wedge\cdots\wedge 
dz_{i-1}\wedge dz_{i+1}
\wedge\cdots\wedge dz_{n}}{2z_{i}}\right|_{V_{i}}
=
\left.
{\rm Res}\left(\frac{dz_{0}\wedge\cdots\wedge dz_{n}}
{z_{0}^{2}+\cdots+z_{n}^{2}}\right)\right|_{V_{i}}.
$$
Hence we have the following formula on $V_{i}$
\begin{equation}
\begin{aligned}
\iota(\xi_{0})\Phi_{\alpha}|_{V_{i}}
&=
\iota(
\sum_{j=0}^{n}\bar{\partial}\chi_{j}\otimes
\frac{1}{2z_{j}}\frac{\partial}{\partial z_{j}})\,
\sum_{k=0}^{n}
\theta^{k}_{\alpha}\,[dz_{k}]\otimes\eta|_{V_{i}}
\\
&=
\frac{1}{2}\sum_{j=0}^{n}
\frac{\bar{\partial}\chi_{j}\wedge
\theta^{j}_{\alpha}}{z_{j}}
\wedge\eta|_{V_{i}}
=
\frac{1}{4}\sum_{i=0}^{n}(-1)^{n+i}\,
z_{i}^{n-3}\,\theta^{j}_{\alpha}
\wedge\bar{\partial}\omega_{ij},
\end{aligned}
\end{equation}
where we used the following relations 
to get the second equality:
$$
\iota(
\sum_{j=0}^{n}\bar{\partial}\chi_{j}\otimes
\frac{1}{2z_{j}}\frac{\partial}{\partial z_{j}})\,
\pi^{*}ds_{k}=0,
\qquad
k=0,\ldots,N.
$$
\par
Let $q\colon\widetilde{\cal X}\to{\cal X}$ be 
the blowing-up along the submanifold 
$\Sigma_{\pi}\subset V$ with exceptional divisor 
$E:=q^{-1}(\Sigma_{\pi})=
{\Bbb P}(N_{\Sigma_{\pi}/{\cal X}})$.
Then 
$q|_{E}\colon{\Bbb P}(N_{\Sigma_{\pi}/V})\to\Sigma_{\pi}$
is the standard projection. 
Since $n\geq3$ and 
since $\{U_{i}\times\varDelta^{N}\}_{i}$
is an open covering of $\widetilde{V}:=q^{-1}(V)$, 
we deduce from Lemma 6.3 and (6.1) that
$q^{*}(\iota(\xi_{0})\Phi_{\alpha})\in
A^{(n,n)}(\widetilde{\cal X})$. 
\par
Set $\widetilde{\pi}=\pi\circ q$. 
By King \cite[Th.\,3.3.2]{King71}, we have
$\widetilde{\pi}_{*}q^{*}
(\iota(\xi_{\alpha})\Phi_{\alpha})\in C^{0}(S)$. 
Since
$$
g_{\alpha0}|_{S\setminus{\frak D}}
=
\pi_{*}(\iota(\xi_{0})\Phi_{\alpha}) 
=
\widetilde{\pi}_{*}q^{*}(\iota(\xi_{0})\Phi_{\alpha}),
$$
$g_{\alpha0}|_{S\setminus{\frak D}}$ extends to
a continuous function on $S$. 
\newline
{\bf (Step 5)}
Let $s_{0}\in{\frak D}$. We must prove
$\lim_{S\setminus{\frak D}\ni s\to s_{0}}
g_{\alpha\beta}|_{S\setminus{\frak D}}(s)=
g_{\alpha\beta}(s_{0})$.
Let $Y_{s_{0}}$ be the proper transform of $X_{s_{0}}$. 
Since
$q^{-1}(X_{s_{0}})=Y_{s_{0}}\cup E$
and since $g_{\alpha\beta}|_{S\setminus{\frak D}}$ 
extends to a continuous function on $S$, we get
$$
\lim_{s\to s_{0}}
g_{\alpha\beta}|_{S\setminus{\frak D}}(s)
=
\int_{q^{-1}(X_{s_{0}})}
q^{*}(\iota(\xi_{\beta})\Phi_{\alpha})
=
\int_{Y_{s_{0}}}q^{*}(\iota(\xi_{\beta})\Phi_{\alpha})
+
\int_{E}q^{*}(\iota(\xi_{\beta})\Phi_{\alpha}).
$$
Since
$q^{*}(\iota(\xi_{\beta})\Phi_{\alpha})|_{E}=0$ 
by Lemma 6.3 and (6.1), we get
$$
\lim_{s\to s_{0}}g_{\alpha\beta}
|_{S\setminus{\frak D}}(s)
=
\int_{Y_{s_{0}}}q^{*}(\iota(\xi_{\beta})\Phi_{\alpha})
=
\int_{(X_{s_{0}})_{\rm reg}}\iota(\xi_{\beta})
\Phi_{\alpha}
=
\langle\phi_{\alpha}|_{X_{s_{0}}},
\rho_{s_{0}}(\frac{\partial}{\partial s_{\beta}})
\rangle_{s_{0}}
=g_{\alpha\beta}(s_{0}),
$$
where we used Lemma 2.9 to get the third equality.
This proves $g_{\alpha\beta}(s)\in C^{0}(S)$.
This completes the proof of Theorem 6.2.
\qed

%%%%%%%%%%%%%%%%%%%%%%%%%%%%%%%%%%%%%%%%%%%%%%%%%%%%%%%%%%%
%
%                 Section 7
%
%%%%%%%%%%%%%%%%%%%%%%%%%%%%%%%%%%%%%%%%%%%%%%%%%%%%%%%%%%%

\section
{\bf Behaviors of the Weil-Petersson metric and 
the Hodge metric}
\par
In this section, we study the boundary behavior
of the Weil-Petersson metric and the Hodge metric
for one-parameter families of Calabi-Yau threefolds
that shall be used later. 
We first recall some basic notions
about positive $(1,1)$-current and give two lemmas
on harmonic functions on $\varDelta^{*}$.
\par

\subsection
{}{\bf Positive $(1,1)$-currents and 
their trivial extensions}
\par
Let $u$ be a $(1,1)$-current on $\varDelta$. 
Then $u$ is {\it positive }
if $u$ is real and if the inequality $u(\varphi)\geq0$
holds for all non-negative function 
$\varphi\in C_{0}^{\infty}(\varDelta)$.
For real $(1,1)$-currents $u$, $v$ on $\varDelta$, 
$u\geq v$ if $u-v$ is a positive $(1,1)$-current 
on $\varDelta$.
For a divisor $H$ on $\varDelta$, let $\delta_{H}$ be 
the current of integration over $H$. 
A real-valued function 
$f\in L^{1}_{\rm loc}(\varDelta)$ is
{\it subharmonic} if $f$ is upper semi-continuous 
and if $dd^{c}f\geq0$ as currents on $\varDelta$. 
\par
Let $\omega_{\varDelta^{*}}$
be the K\"ahler form of the Poincar\'e metric 
on $\varDelta^{*}$:
$$
\omega_{\varDelta^{*}}
:=
\frac{\sqrt{-1}dt\wedge d\bar{t}}
{|t|^{2}(-\log|t|^{2})^{2}}
=
-dd^{c}\log(-\log|t|^{2}).
$$ 
\par
A $C^{\infty}$ real $(1,1)$-form $T$ on 
$\varDelta^{*}$ has Poincar\'e growth if there exists
$C>0$ with
\begin{equation}
-C\,\omega_{\varDelta^{*}}
\leq 
T
\leq 
C\,\omega_{\varDelta^{*}}.
\end{equation}
In that case, the coefficient of $T$ lies in 
$L^{1}_{\rm loc}(\varDelta)$.
The $(1,1)$-current on $\varDelta$ defined by
$$
\widetilde{T}(\psi)
:=
\int_{\varDelta}\psi\,T,
\qquad
\psi\in C_{0}^{\infty}(\varDelta)
$$
is called the {\it trivial extension} of $T$ 
from $\varDelta^{*}$ to $\varDelta$. We have
$\widetilde{\omega_{\varDelta^{*}}}=
-dd^{c}\log(-\log|t|^{2})$ as currents on $\varDelta$.

\subsection
{}{\bf Two lemmas on harmonic functions 
on $\varDelta^{*}$}

\begin{lem}
Let $H(t)$ be a real-valued harmonic function 
on $\varDelta^{*}$. 
\newline{$(1)$}
There exist $c\in{\Bbb R}$ and 
$F(t)\in{\cal O}(\varDelta^{*})$
with $H(t)=c\,\log|t|^{2}+2{\rm Re}\,F(t)$.
\newline{$(2)$}
If there exist $\gamma\in{\Bbb R}$ such that
$|t|^{\gamma}e^{H(t)}\in 
L^{1}_{\rm loc}(\varDelta)$,
then $F(t)\in{\cal O}(\varDelta)$.
\newline{$(3)$}
If $H(t)=O(\log(-\log|t|))$ as $t\to0$,
then $H(t)$ extends to a harmonic function
on $\varDelta$.
\end{lem}

\begin{pf}
{\bf (1)}
Since $H(t)$ is harmonic on $\varDelta^{*}$,
there exists $f(t)\in{\cal O}(\varDelta^{*})$
with $\partial H(t)=f(t)\,dt$. 
Let $f(t)=\sum_{n\in{\Bbb Z}}a_{n}\,t^{n}$
be the Laurent expansion of $f(t)$
and define the meromorphic function $F(t)$
on $\varDelta^{*}$ by
$F(t):=\sum_{n\not=-1}\frac{a_{n}}{n+1}\,t^{n+1}$.
By the reality of $H(t)$, we get
$dH(t)=
a_{-1}\frac{dt}{t}+\overline{a}_{-1}\frac{dt}{t}
+dF(t)+d\bar{F}(t)$.
Integrating the both hand sides 
over the circle $|t|=1/2$,
we get $a_{-1}\in{\Bbb R}$ by the Stokes theorem, 
so that
$dH(t)=a_{-1}\,d\log|t|^{2}
+2d\{{\rm Re}\,F(t)\}$.
This proves (1).
\newline{\bf (2)}
By assumption, we get
\begin{equation}
\int_{|t|<1/2}
|t|^{\gamma+2c}\,|e^{F(t)}|^{2}\,
\sqrt{-1}\,dt\wedge d\bar{t}
<+\infty.
\end{equation}
Since $e^{F(t)}$ is holomorphic 
on $\varDelta^{*}$, we deduce from (7.2)
that $e^{F(t)}$ is a meromorphic function
on $\varDelta$. 
There exist $\nu\in{\Bbb Z}$ and 
a nowhere vanishing holomorphic function 
$\epsilon(t)\in{\cal O}(\varDelta)$ with
$e^{F(t)}=t^{\nu}\,\epsilon(t)$.
Then $F'(t)=\nu\,t^{-1}+\epsilon'(t)\epsilon(t)^{-1}$.
Since $F(t)$ is a meromorphic function
on $\varDelta^{*}$, the residue of $F'(t)$
must vanish, i.e., $\nu=0$.
Thus we have proved that $F(t)=\log\epsilon(t)$ 
is holomorphic on $\varDelta$.
\newline{\bf (3)}
Since $e^{H(t)}\in L^{1}_{\rm loc}(\varDelta)$,
$H(t)-c\,\log|t|^{2}$ is a harmonic function
on $\varDelta$ by (1), (2). Since
$H(t)=O(\log(-\log|t|))$ as $t\to0$, we get $c=0$.
This completes the proof.
\end{pf}

\begin{lem}
Let $\lambda(t)$ be a positive, 
locally $L^{m}$-integrable
function on $\varDelta$ for some $m>0$.
Let $\chi(t)$ be a function on $\varDelta^{*}$
satisfying $\chi(t)\leq C\,(-\log|t|+2)$, 
where $C\in{\Bbb R}$ is a constant.
If $\log\lambda(t)+\chi(t)$ is harmonic 
on $\varDelta^{*}$, then there exists 
$c\in{\Bbb R}$ such that
$$
\log\lambda(t)
=
c\,\log|t|^{2}+O(|\chi(t)|+1)
\qquad
(t\to0).
$$
\end{lem}

\begin{pf}
Set $H(t):=\log\lambda(t)+\chi(t)$.
Since $\chi(t)\leq C\,(-\log|t|+2)$, we get
\begin{equation}
\log\lambda(t)
=
H(t)-\chi(t)
\geq
H(t)-C\,(-\log|t|+2).
\end{equation}
Since $\lambda(t)\in L^{m}(\varDelta(1/2))$, we get
\begin{equation}
e^{-2Cm}\int_{\varDelta(1/2)}
|t|^{Cm}e^{m\,H(t)}\,\sqrt{-1}\,dt\wedge d\bar{t}
\leq
\int_{\varDelta(1/2)}\lambda(t)^{m}\,
\sqrt{-1}\,dt\wedge d\bar{t}
<+\infty.
\end{equation}
By (7.4) and Lemma 7.1 (1), (2), 
there exists $c\in{\Bbb R}$ and 
$F(t)\in{\cal O}(\varDelta)$ with
\begin{equation}
H(t)=c\,\log|t|^{2}+2\,{\rm Re}\,F(t).
\end{equation}
Since $\log\lambda(t)=H(t)-\chi(t)$,
the result follows from (7.5).
\end{pf}

\subsection
{}
{\bf The boundary behaviors}
\par
In Subsect.\,7.3, we fix the following notation.
Let $\cal X$ be a (possibly) singular complex fourfold 
and let $\pi\colon{\cal X}\to\varDelta$ be
a proper surjective holomorphic function.
Assume that $X_{t}:=\pi^{-1}(t)$ is a smooth
Calabi-Yau threefold
for $t\in\varDelta^{*}$.
We do {\it not} assume that the central fiber
$X_{0}$ has only ODP's as its singular set.
Recall that the Weil-Petersson form 
$\omega_{{\rm WP},{\cal X}/\varDelta}$ and
the Hodge form
$\omega_{{\rm H},{\cal X}/\varDelta}$ for
$\pi\colon{\cal X}\to\varDelta$ were defined 
in Sects.\,4.3 and 4.4.3, respectively.

\begin{prop}
There exists a positive constant $C$ such that
\begin{equation}
0
\leq
\omega_{{\rm WP},{\cal X}/\varDelta}
\leq 
C\,\omega_{\varDelta^{*}},
\qquad
0
\leq
\omega_{{\rm H},{\cal X}/\varDelta}
\leq 
C\,\omega_{\varDelta^{*}}.
\end{equation}
In particular, the positive $(1,1)$-forms
$\omega_{{\rm WP},{\cal X}/\varDelta}$
and $\omega_{{\rm H},{\cal X}/\varDelta}$
on $\varDelta^{*}$ extend trivially to
closed positive $(1,1)$-currents on $\varDelta$.
\end{prop}

\begin{pf}
We follow \cite[Proof of Th.\,5.1]{LuSun04}.
Since (7.6) is obvious 
when $\omega_{{\rm H},{\cal X}/\varDelta}=0$,
we assume that $\omega_{{\rm H},{\cal X}/\varDelta}$ 
does not vanish identically on $\varDelta^{*}$.
Shrinking $\varDelta$ if necessary, we may assume 
that $\omega_{{\rm H},{\cal X}/\varDelta}$ is 
strictly positive on $\varDelta^{*}$.
Let $b\in\varDelta^{*}$.
Since $\omega_{{\rm H},{\cal X}/\varDelta}$ is 
non-degenerate at $b$, the deformation germ
$\pi\colon({\cal X},X_{b})\to(\varDelta,b)$
is induced from the Kuranishi family 
by an immersion of germs
$(\varDelta,b)\hookrightarrow({\rm Def}(X_{b}),[X_{b}])$.
Let $\omega_{\rm H}$ be the Hodge form 
on ${\rm Def}(X_{b})$.
By \cite[Th.\,1.1.2]{Lu01}, the holomorphic sectional
curvature of $({\rm Def}(X_{b}),\omega_{\rm H})$ is
bounded from above by $\alpha:=-(5+2\sqrt{3})^{-1}$.
Since $b\in\varDelta^{*}$ is an arbitrary point, 
the holomorphic sectional curvature of
$(\varDelta^{*},\omega_{{\rm H},{\cal X}/\varDelta})$
is bounded from above by $\alpha$ 
(cf. e.g. \cite[Prop.\,2.3.9]{Kobayashi98}).
The second inequality of (7.6) follows from 
the Schwarz lemma \cite[Th.\,2.3.5]{Kobayashi98}.
The first inequality of (7.6) follows from
the second one because 
$2\,\omega_{{\rm WP},{\cal X}/\varDelta}\leq
\omega_{{\rm H},{\cal X}/\varDelta}$ by 
\cite[p.107, l.17]{Lu01}.
\par
Since $(\varDelta(r)^{*},\omega_{\varDelta^{*}})$ 
has finite volume when $r<1$, 
the positive $(1,1)$-forms
$\omega_{{\rm WP},{\cal X}/\varDelta}$
and $\omega_{{\rm H},{\cal X}/\varDelta}$ 
extend trivially to closed positive $(1,1)$-currents 
on $\varDelta$.
\end{pf}

\begin{defn}
Define $\Omega_{{\rm WP},{\cal X}/\varDelta}$ and
$\Omega_{{\rm H},{\cal X}/\varDelta}$ 
as the trivial extensions of
$\omega_{{\rm WP},{\cal X}/\varDelta}$ and
$\omega_{{\rm H},{\cal X}/\varDelta}$ 
from $\varDelta^{*}$ to $\varDelta$, respectively.
\end{defn}

\begin{lem}
Let $A,B\in{\Bbb R}$.
Let $\lambda(t)$ be a positive, 
locally $L^{m}$-integrable $C^{\infty}$ function
on $\varDelta^{*}$ for some $m>0$ such that
$-dd^{c}\log\lambda=
A\,\omega_{{\rm H},{\cal X}/\varDelta}
+B\,\omega_{{\rm WP},{\cal X}/\varDelta}$.
\newline{$(1)$}
There exists $c\in{\Bbb R}$
such that as $t\to0$,
$$
\log\lambda(t)=c\,\log|t|^{2}+O(\log(-\log|t|)).
$$
\newline{$(2)$}
With the same constant $c$ as above,
the following equation of currents on $\varDelta$
holds:
$$
-dd^{c}\log\lambda
=
A\,\Omega_{{\rm H},{\cal X}/\varDelta}
+B\,\Omega_{{\rm WP},{\cal X}/\varDelta}
-c\,\delta_{0}.
$$
\end{lem}

\begin{pf}
We follow \cite[Prop.\,3.11]{Yoshikawa04}.
By \cite[Proof of Lemma 5.4]{Siu74},
there exist subharmonic functions
$\varphi$ and $\theta$ on $\varDelta$
such that the following equations of currents
on $\varDelta$ hold:
\begin{equation}
\Omega_{{\rm WP},{\cal X}/\varDelta}=dd^{c}\varphi,
\qquad
\Omega_{{\rm H},{\cal X}/\varDelta}=dd^{c}\theta.
\end{equation}
Since $\varphi$ and $\theta$ are subharmonic,
there exists $C_{0}\in{\Bbb R}$ with
\begin{equation}
\varphi(t)\leq C_{0},
\qquad
\theta(t)\leq C_{0},
\qquad
t\in\varDelta(1/2).
\end{equation}
Since 
$\widetilde{\omega_{\varDelta^{*}}}=
-dd^{c}\log(-\log|t|)$ as a current on $\varDelta$, 
we deduce from (7.6) that
$$
dd^{c}\left\{-C\log(-\log|t|)-\varphi\right\}
=C\,\widetilde{\omega_{\varDelta^{*}}}-
\Omega_{{\rm WP},{\cal X}/\varDelta}
\geq0,
$$
$$
dd^{c}\left\{-C\log(-\log|t|)-\theta\right\}
=C\,\widetilde{\omega_{\varDelta^{*}}}-
\Omega_{{\rm H},{\cal X}/\varDelta}
\geq0.
$$
Hence $-C\log(-\log|t|)-\varphi$ and
$-C\log(-\log|t|)-\theta$ are subharmonic
functions on $\varDelta$, so that there 
exists $C_{1}\in{\Bbb R}$ with
\begin{equation}
-C\log(-\log|t|)-\varphi(t)\leq C_{1},
\qquad
-C\log(-\log|t|)-\theta(t)\leq C_{1},
\qquad
\forall\,t\in\varDelta(1/2).
\end{equation}
By (7.8) and (7.9), there exists 
$C_{2}\in{\Bbb R}$ such that 
for all $t\in\varDelta(1/2)$,
\begin{equation}
-C\,\log(-\log|t|)-C_{1}
\leq
\varphi(t)
\leq
C_{0},
\qquad
-C\,\log(-\log|t|)-C_{1}
\leq
\theta(t)
\leq
C_{0}.
\end{equation}
\par
Set
\begin{equation}
H(t):=
\log\lambda(t)+A\,\theta(t)+B\,\varphi(t).
\end{equation}
Since $dd^{c}H=0$,
$H(t)$ is a harmonic function 
on $\varDelta^{*}$.
Since $\lambda(t)$ is locally $L^{m}$-integrable 
on $\varDelta$,
the first assertion follows from (7.10) and Lemma 7.2
by setting $\chi(t)=A\,\theta(t)+B\,\varphi(t)$.
The second assertion follows from (7.5), (7.7), (7.11).
\end{pf}

Let $g_{{\rm WP},{\cal X}/\varDelta}$ be
the K\"ahler metric on $\varDelta^{*}$ whose
K\"ahler form is 
$\omega_{{\rm WP},{\cal X}/\varDelta}$.

\begin{prop}
Assume that $h^{1,2}(X_{t})=1$ for all 
$t\in\varDelta^{*}$.
\newline{$(1)$} 
There exists $\alpha\in{\Bbb R}$ such that
as $t\to0$:
$$
\log g_{{\rm WP},{\cal X}/\varDelta}
\left(
\frac{\partial}{\partial t},
\frac{\partial}{\partial\bar{t}}
\right)
=
\alpha\,\log|t|^{2}+O(\log(-\log|t|)).
$$
\newline{$(2)$ }
With the same constant $\alpha$ as above,
the following equation of currents 
on $\varDelta$ holds:
$$
dd^{c}\log g_{{\rm WP},{\cal X}/\varDelta}
\left(
\frac{\partial}{\partial t},
\frac{\partial}{\partial\bar{t}}
\right)
=
\alpha\,\delta_{0}
-\Omega_{{\rm H},{\cal X}/\varDelta}
+4\,\Omega_{{\rm WP},{\cal X}/\varDelta}.
$$
\newline{$(3)$ }
If $X_{0}$ is a Calabi-Yau threefold with at most
one ODP and if 
$\pi\colon{\cal X}\to\varDelta$ is the Kuranishi
family of $X_{0}$, then $\alpha=0$.
\end{prop}

\begin{pf}
{\bf (1)} 
Set
$\lambda(t):=
g_{{\rm WP},{\cal X}/\varDelta}
(\frac{\partial}{\partial t},
\frac{\partial}{\partial\bar{t}})$
and $A=1$, $B=-4$ in Lemma 7.5. 
By the definition of Hodge form, we have 
$-dd^{c}\log\lambda=\omega_{{\rm H},{\cal X}/\varDelta}
-4\,\omega_{{\rm WP},{\cal X}/\varDelta}$
on $\varDelta^{*}$.
Since $\lambda(t)\in L^{1}_{\rm loc}(\varDelta)$ 
by Proposition 7.3,
the result follows from Lemma 7.5 (1). 
\newline{\bf (2)}
The result follows from Lemma 7.5 (2).
\newline{\bf (3)}
The result follows from \cite[Cor.\,5.1]{Tian92}.
This completes the proof.
\end{pf}

If ${\cal X}$ is smooth, $\pi_{*}K_{\cal X}$ is
locally free by \cite[p.391, Th.\,V]{Takegoshi95}.
Since $K_{\varDelta}$ is trivial and since
$h^{0}(X_{t},K_{\cal X}|X_{t})=1$ for $t\in\varDelta^{*}$,
$\pi_{*}K_{{\cal X}/\varDelta}=
\pi_{*}(K_{\cal X}\otimes\pi^{*}K_{\varDelta}^{-1})
\cong\pi_{*}K_{\cal X}$ is an invertible sheaf 
on $\varDelta$ in that case.

\begin{lem}
Assume that $X_{t}$ is Calabi-Yau for
all $t\in\varDelta^{*}$.
If ${\cal X}$ is smooth,
there exists $\xi\in H^{0}({\cal X},K_{\cal X})$
such that ${\rm div}(\xi)\subset X_{0}$.
\end{lem}

\begin{pf}
Since $\pi_{*}K_{\cal X}$ is an invertible sheaf
on $\varDelta$,
there exists $\xi\in H^{0}({\cal X},K_{\cal X})=
H^{0}(\varDelta,\pi_{*}K_{\cal X})$ that generates
$\pi_{*}K_{\cal X}$ as an ${\cal O}_{\varDelta}$-module, 
i.e.,
$\pi_{*}K_{\cal X}={\cal O}_{\varDelta}\cdot\xi$.
Since $H^{0}(X_{t},K_{\cal X}|_{X_{t}})\cong
H^{0}(X_{t},K_{X_{t}})\cong{\Bbb C}$ 
for all $t\in\varDelta^{*}$, we get
$H^{0}(X_{t},K_{\cal X}|_{X_{t}})={\Bbb C}\,\xi|_{X_{t}}$
in that case by 
\cite[Chap.\,3, Th.\,4.12 (ii)]{Banica76}. 
Since 
$K_{\cal X}|_{X_{t}}\cong K_{X_{t}}\cong{\cal O}_{X_{t}}$
for $t\in\varDelta^{*}$, $\xi|_{X_{t}}$ is nowhere vanishing
on $X_{t}$, $t\in\varDelta^{*}$.
This proves the lemma.
\end{pf}

If ${\cal X}$ is smooth, there exists
$\xi\in H^{0}({\cal X},K_{\cal X})$ by Lemma 7.7
such that ${\rm div}(\xi)\subset X_{0}$.
In that case, we define a section
$\eta_{{\cal X}/\varDelta}\in 
H^{0}({\cal X},K_{{\cal X}/\varDelta})$ by
$\eta_{{\cal X}/\varDelta}
:=\xi\otimes(\pi^{*}dt)^{-1}$.
We identify 
$\eta_{{\cal X}/\varDelta}|_{X_{t}}$ 
with the Poincar\'e residue
$\eta_{t}:={\rm Res}_{X_{t}}\xi/(\pi-t)
\in H^{0}(X_{t},K_{X_{t}})$ for $t\in\varDelta^{*}$.
Then
\begin{equation}
\xi|_{X_{t}}=\eta_{t}\otimes d\pi,
\end{equation}
and $\eta_{{\cal X}/\varDelta}$ is regarded 
as a family of holomorphic $3$-forms.
We also regard $\eta_{{\cal X}/\varDelta}$
as the corresponding element of 
$H^{0}(\varDelta,\pi_{*}K_{{\cal X}/\varDelta})$.

\begin{prop}
Assume that ${\cal X}$ is smooth.
Let $\eta_{{\cal X}/\varDelta}$ be a nowhere vanishing 
holomorphic section of $\pi_{*}K_{{\cal X}/\varDelta}$.
\newline{$(1)$} 
There exists $\beta\in{\Bbb R}$ such that
as $t\to0$:
$$
\log \|\eta_{{\cal X}/\varDelta}(t)\|_{L^{2}}^{2}
=\beta\,\log|t|^{2}+O(\log(-\log|t|)).
$$
\newline{$(2)$ }
With the same constant $\beta$ as above,
the following equation of currents 
on $\varDelta$ holds:
$$
dd^{c}\log \|\eta_{{\cal X}/\varDelta}(t)\|_{L^{2}}^{2}
=
\beta\,\delta_{0}
-\Omega_{{\rm WP},{\cal X}/\varDelta}.
$$
\newline{$(3)$ }
If $X_{0}$ is a Calabi-Yau threefold with at most
one ODP and if $\xi$ is nowhere vanishing on ${\cal X}$, 
then 
$\log \|\eta_{{\cal X}/\varDelta}(t)\|_{L^{2}}^{2}$ 
extends to a continuous function on $\varDelta$.
In particular, $\beta=0$.
\end{prop}

\begin{pf}
{\bf (1) }
Set 
$\lambda(t):=
\|\eta_{{\cal X}/\varDelta}(t)\|_{L^{2}}^{2}$ 
and $A=0$, $B=1$ in Lemma 7.5. Since
$$
\begin{aligned}
\int_{\varDelta(1/2)}
\lambda(t)\sqrt{-1}dt\wedge d\bar{t}
&
=
\int_{\varDelta(1/2)}
\pi_{*}(\sqrt{-1}\eta_{{\cal X}/\varDelta}\wedge
\overline{\eta_{{\cal X}/\varDelta}})
\sqrt{-1}dt\wedge d\bar{t}
\\
&
=
\int_{\pi^{-1}(\varDelta(1/2))}\xi\wedge\bar{\xi}
<+\infty
\end{aligned}
$$
by (7.12), we get
$\lambda(t)\in L^{1}_{\rm loc}(\varDelta)$.
Since $-dd^{c}\log\lambda
=\omega_{{\rm WP},{\cal X}/\varDelta}$ 
by the definition of the Weil-Petersson form,
the result follows from Lemma 7.5 (1).
\newline{\bf (2)}
The result follows from Lemma 7.5 (2).
\newline{\bf (3)}
The result follows from e.g.
\cite[Proof of Th.\,8.1]{Yoshikawa98}.
This completes the proof.
\end{pf}

\subsection
{}
{\bf The boundary behavior of the anomaly term}
\par
In Subsection 7.4, we fix the following notation.
Let $\pi\colon{\cal X}\to\varDelta$ be a proper surjective
holomorphic function on a smooth K\"ahler fourfold
with critical locus $\Sigma_{\pi}$, so that
$\pi$ has relative dimension $3$.
Assume that $\Sigma_{\pi}\subset X_{0}$ and that
$X_{t}$ is a smooth Calabi-Yau threefold for all 
$t\in\varDelta^{*}$.
\par
Let $g_{\cal X}$ be a K\"ahler metric on ${\cal X}$.
Let $\gamma_{\cal X}$ be the K\"ahler form 
of $g_{\cal X}$ and set 
$\gamma_{t}:=\gamma_{\cal X}|_{X_{t}}$.
Recall that the anomaly term ${\cal A}(X_{t},\gamma_{t})$ 
was defined in Definition 4.1. 
The following result is a generalization of
\cite[(6.17), (6.19)]{Yoshikawa04}.

\begin{prop}
$(1)$ There exists $c\in{\Bbb R}$ such that as $t\to0$:
$$
\log{\cal A}(X_{t},\gamma_{t})
=
c\,\log|t|^{2}+O(\log(-\log|t|)).
$$
\newline{$(2)$}
If $\Sigma_{\pi}$ consists of a unique ODP and 
if $X_{0}$ is Calabi-Yau, then as $t\to0$
$$
\log{\cal A}(X_{t},\gamma_{t})
=
-\frac{1}{12}\,\log|t|^{2}+O(1).
$$
\end{prop}

\begin{pf}
{\bf (1)}
Let $g_{{\cal X}/\varDelta}$ be the Hermitian metric
on $T{\cal X}/\varDelta$ induced from $g_{\cal X}$,
and let $\gamma_{{\cal X}/\varDelta}$ be the 
corresponding $(1,1)$-form on $T{\cal X}/\varDelta$.
Then we may identify $\gamma_{{\cal X}/\varDelta}$
with the family of K\"ahler forms
$\{\gamma_{t}\}_{t\in\varDelta}$.
Let $N_{X_{t}/{\cal X}}^{*}$ be the conormal bundle
of $X_{t}$ in ${\cal X}$ for $t\in\varDelta^{*}$.
Then 
$d\pi=\pi^{*}dt\in H^{0}(X_{t},N^{*}_{X_{t}/{\cal X}})$
generates $N^{*}_{X_{t}/{\cal X}}$ for 
$t\in\varDelta^{*}$, so that 
$N^{*}_{X_{t}/{\cal X}}$ is trivial in that case.
Since the Hermitian metric
on $\Omega^{1}_{X_{t}}$ is induced from $g_{\cal X}$
via the $C^{\infty}$ identification
$\Omega^{1}_{X_{t}}\cong(N_{X_{t}/{\cal X}}^{*})^{\perp}$
and since 
$(\gamma_{{\cal X}/\varDelta}^{3}/3!)|_{X_{t}}$ 
is the volume form on $\Omega^{1}_{X_{t}}$,
we get
\begin{equation}
\frac{\gamma_{\cal X}^{4}}{4!}
=
\frac{\gamma_{{\cal X}/\varDelta}^{3}}{3!}\wedge
\left(\sqrt{-1}\,\frac{d\pi}{\|d\pi\|}\wedge
\frac{\overline{d\pi}}{\|d\pi\|}\right).
\end{equation}
\par
By Lemma 7.7, there exists
$\xi\in H^{0}({\cal X},K_{\cal X})$
such that ${\rm div}(\xi)\subset X_{0}$.
As before, define $\eta_{{\cal X}/\varDelta}\in 
H^{0}({\cal X},K_{{\cal X}/\varDelta})$ by
$\eta_{{\cal X}/\varDelta}
:=\xi\otimes(\pi^{*}dt)^{-1}$,
and identify 
$\eta_{{\cal X}/\varDelta}|_{X_{t}}$ 
with the Poincar\'e residue
$\eta_{t}:={\rm Res}_{X_{t}}\xi/(\pi-t)
\in H^{0}(X_{t},K_{X_{t}})$ for $t\in\varDelta^{*}$.
Then
$\eta_{{\cal X}/\varDelta}$ is regarded 
as a family of holomorphic $3$-forms $\{\eta_{t}\}$.
By (7.12) and (7.13), we get
\begin{equation}
\frac{\sqrt{-1}\,\eta_{{\cal X}/\varDelta}\wedge
\overline{\eta_{{\cal X}/\varDelta}}}
{\gamma_{{\cal X}/\varDelta}^{3}/3!}
=
\frac{(-1)^{3}\sqrt{-1}\,\xi\wedge\overline{\xi}}
{(\gamma_{{\cal X}/\varDelta}^{3}/3!)\wedge 
d\pi\wedge\overline{d\pi}}
=
\frac{\xi\wedge\overline{\xi}}
{\gamma_{{\cal X}/\varDelta}^{4}/4!}
\cdot
\frac{1}{\|d\pi\|^{2}}
=
\frac{\|\xi\|^{2}}{\|d\pi\|^{2}}.
\end{equation}
\par
Let $X$ denote a general fiber of 
$\pi\colon{\cal X}\to\varDelta$.
Let ${\cal A}({\cal X}/\varDelta)$ be the function 
on $\varDelta^{*}$ defined by 
${\cal A}({\cal X}/\varDelta)(t):=
{\cal A}(X_{t},\gamma_{t})$.
Then
\begin{equation}
\begin{aligned}
\log{\cal A}({\cal X}/\varDelta)
&=-\frac{1}{12}
\pi_{*}\left[\log\left(
\frac{\sqrt{-1}\,\eta_{{\cal X}/\varDelta}\wedge
\overline{\eta_{{\cal X}/\varDelta}}}
{\gamma_{{\cal X}/\varDelta}^{3}/3!}\right)\,
c_{3}(T{\cal X}/\varDelta,g_{{\cal X}/\varDelta})
\right]
\\
&\quad
+\frac{\chi(X)}{12}\log
\|\eta_{{\cal X}/\varDelta}\|_{L^{2}}^{2}.
\end{aligned}
\end{equation}
\par
We use the notation in Subsection 5.3.
Hence $q\colon\widetilde{\cal X}\to{\cal X}$
is the resolution of the Gauss maps $\mu$ and $\nu$.
Substituting (7.14) into (7.15) and using (5.1), 
we get
\begin{equation}
\begin{aligned}
\log{\cal A}({\cal X}/\varDelta)
&
=
-\frac{1}{12}\pi_{*}\left[\log
\left(\frac{\|\xi\|^{2}}{\|d\pi\|^{2}}\right)\,
c_{3}(T{\cal X}/\varDelta,g_{{\cal X}/\varDelta})\right]
+\frac{\chi(X)}{12}
\log\|\eta_{{\cal X}/\varDelta}\|_{L^{2}}^{2}
\\
&
=
-\frac{1}{12}\widetilde{\pi}_{*}\left[\log
q^{*}\left(\frac{\|\xi\|^{2}}{\|d\pi\|^{2}}\right)\,
\widetilde{\mu}^{*}c_{3}(U,g_{U})\right]
+\frac{\chi(X)}{12}
\log\|\eta_{{\cal X}/\varDelta}\|_{L^{2}}^{2}.
\end{aligned}
\end{equation}
Since 
${\rm div}(q^{*}\xi)\subset\widetilde{\pi}^{-1}(0)$
by the condition ${\rm div}(\xi)\subset X_{0}$, 
the assertion follows from Lemma 5.8 
and Proposition 7.8 (1) applied to the second line
of (7.16). 
\newline
{\bf (2)}
Assume that $\Sigma_{\pi}$ consists of a unique ODP 
and that $X_{0}$ is Calabi-Yau. We use the notation
in Subsect.\,5.9.
We may assume by Lemma 6.1
that $\xi$ is nowhere vanishing on ${\cal X}$.
Hence ${\rm div}(q^{*}\xi)=\emptyset$,
and
$\widetilde{\pi}_{*}\{q^{*}\log\|\xi\|^{2}\,
\widetilde{\mu}^{*}c_{3}(U,g_{U})\}$ and
$\log\|\eta_{{\cal X}/\varDelta}\|_{L^{2}}^{2}$ 
are bounded as $t\to0$ 
by the first equation of Lemma 5.8 and 
by Proposition 7.8 (3).
We deduce from (7.16) that
\begin{equation}
\log{\cal A}({\cal X}/\varDelta)
=
\frac{1}{12}\widetilde{\pi}_{*}\{
q^{*}(\log\|d\pi\|^{2})\,
\widetilde{\mu}^{*}c_{3}(U,g_{U})\}
+O(1).
\end{equation}
Since $E={\Bbb P}^{3}$ and 
$\varphi=(-1)^{3}\widetilde{\mu}^{*}c_{3}(U)$ 
in the second equation of Lemma 5.8, we get
\begin{equation}
\log{\cal A}({\cal X}/\varDelta)(t)
=
\left(\frac{1}{12}\int_{{\Bbb P}^{3}}c_{3}(U)\right)\,
\log|t|^{2}+O(1)
=
\frac{(-1)^{3}}{12}\log|t|^{2}+O(1).
\end{equation}
This proves (2).
\end{pf}

\subsection
{}
{\bf The Weil-Petersson and Hodge metrics 
on the Kuranishi space}
\par
In Subsect.\,7.5, we fix the following notation.
Let $X$ be a smoothable Calabi-Yau threefold 
with only one ODP as its singular set,
and let 
${\frak p}\colon({\frak X},X)\to({\rm Def}(X),[X])$
be the Kuranishi family with discriminant locus 
$\frak D$. Assume that
$\dim{\rm Def}(X)=h^{1,2}(X)=1$.
\par
By Lemma 6.1, there exists a nowhere vanishing 
holomorphic $4$-form $\xi$ on $\frak X$.
Then $\eta_{{\frak X}/{\rm Def}(X)}=
\xi\otimes\pi^{*}(ds)^{-1}$ 
is a nowhere vanishing holomorphic section of 
${\frak p}_{*}K_{{\frak X}/{\rm Def}(X)}$.
Set $\eta_{s}:=\eta_{{\frak X}/{\rm Def}(X)}|_{X_{s}}$.
We identify $\eta_{s}$ with the corresponding
holomorphic $3$-form on $(X_{s})_{\rm reg}$ such that
$\eta_{s}\otimes(ds)=\xi|_{X_{s}}$ 
under the canonical isomorphism
$K_{X_{s}}\otimes
{\frak p}^{*}K_{{\rm Def}(X)}|_{X_{s}}=
K_{\frak X}|_{X_{s}}$.
Then $\{\eta_{s}\}_{s\in S}$ is regarded as
a holomorphic family of nowhere vanishing 
holomorphic $3$-forms.
\par
For $p=0,1$ and $q\geq0$, the direct image sheaves
$R^{q}{\frak p}_{*}\Omega^{p}_{{\frak X}/{\rm Def}(X)}$
are locally free by Definition 2.1 (ii) and 
Theorem 2.11. For $p=0,1$, let $\sigma_{p}$ be 
a nowhere vanishing holomorphic section of
$\lambda(\Omega^{p}_{{\frak X}/{\rm Def}(X)})$.
\par
By Proposition 2.8, there exists 
a K\"ahler metric $g_{\frak X}$ on $\frak X$.
Let $g_{{\frak X}/{\rm Def}(X)}$ be the Hermitian
metric on $T{\frak X}/{\rm Def}(X)
|_{{\frak X}\setminus\Sigma_{\frak p}}$ induced from
$g_{\frak X}$. Set $g_{s}:=g_{\frak X}|_{X_{s}}$
for $s\in{\rm Def}(X)$.

\begin{thm}
The following formula holds for $p=0,1$:
$$
\log\|\sigma_{p}(s)\|^{2}_{\lambda
(\Omega^{p}_{{\frak X}/{\rm Def}(X)}),L^{2},
g_{{\frak X}/{\rm Def}(X)}}
=O(\log(-\log|s|)).
$$
\end{thm}

\begin{pf}
Let $p=0$. Let $1$ be the section of
${\frak p}_{*}{\cal O}_{\frak X}$ such that
$1_{s}=1\in H^{0}(X_{s},{\cal O}_{X_{s}})$.
Regard $\eta_{{\frak X}/{\rm Def}(X)}$
as a nowhere vanishing holomorphic section of
$(R^{3}{\frak p}_{*}{\cal O}_{\frak X})^{\lor}$
by the relative Serre duality. Set
$\sigma_{0}:=1\otimes\eta_{{\frak X}/{\rm Def}(X)}$.
Since 
$$
\log\|\sigma_{0}(s)\|_{L^{2},g_{s}}^{2}=
\log{\rm Vol}(X_{s},g_{s})+\log\|\eta_{s}\|_{L^{2}}^{2}
=\log\|\eta_{s}\|_{L^{2}}^{2}+O(1),
$$
the assertion for $p=0$ follows from 
Proposition 7.8 (3).
\par
Let $p=1$.
Let ${\bf e}_{1},\ldots,{\bf e}_{b_{2}(X)}$ be 
a $\Bbb Z$-basis of $H^{2}(X,{\Bbb Z})/{\rm Torsion}$.
There exist holomorphic line bundles 
${\cal L}_{1},\ldots,{\cal L}_{b_{2}(X)}$ on $\frak X$ 
by Lemma 2.16
such that $c_{1}({\cal L}_{i})|_{X}={\bf e}_{i}$
for $1\leq i\leq b_{2}(X)$, and such that 
the Dolbeault cohomology classes of their Chern forms
${\frak C}_{1}({\cal L}_{1}),\ldots,
{\frak C}_{1}({\cal L}_{b_{2}(X)})$ 
form a local basis of 
$R^{1}\pi_{*}\Omega^{1}_{{\frak X}/{\rm Def}(X)}$
as a ${\cal O}_{{\rm Def}(X)}$-module.
\par
By Theorem 6.2, 
$(\rho_{s}^{\lor})^{-1}(ds)\otimes\eta_{s}^{-1}$ 
is a local basis of
$R^{2}\pi_{*}\Omega^{1}_{{\frak X}/{\rm Def}(X)}$
as an ${\cal O}_{{\rm Def}(X)}$-module.
For $s\in{\rm Def}(X)$, set
$$
\sigma_{1}(s):=
({\frak C}_{1}({\cal L}_{1})\wedge\cdots\wedge 
{\frak C}_{1}({\cal L}_{b_{2}(X)}))^{-1}\otimes
((\rho_{s}^{\lor})^{-1}(ds)\otimes\eta_{s}^{-1}).
$$
Then $\sigma_{1}$ is a nowhere vanishing holomorphic
section of 
$\lambda(\Omega^{1}_{{\frak X}/{\rm Def}(X)})$.
\par
Let $\gamma_{s}$ be the K\"ahler form of
$g_{\frak X}|_{X_{s}}$. Since $g_{\frak X}$ is
a K\"ahler metric on $\frak X$, 
the section ${\rm Def}(X)\ni s\to[\gamma_{s}]
\in H^{2}(X_{s},{\Bbb R})$ of
$R^{2}{\frak p}_{*}{\Bbb R}$ is constant.
Let $[\gamma]\in H^{2}(X,{\Bbb R})$ be the element
corresponding to $[\gamma_{s}]$. By Lemma 4.12,
$$
\|{\frak C}_{1}({\cal L}_{1})\wedge\cdots\wedge 
{\frak C}_{1}({\cal L}_{b_{2}(X)})
\|_{L^{2},g_{s}}^{2}(s)
=
{\rm Vol}_{L^{2}}(H^{2}(X,{\Bbb Z}),[\gamma])
\not=0
$$
is a constant function on ${\rm Def}(X)$.
Hence we get
$$
\begin{aligned}
\log\|\sigma_{1}(s)\|_{L^{2},g_{s}}^{2}
&=
-\log{\rm Vol}_{L^{2}}(H^{2}(X,{\Bbb Z}),[\gamma])
-\log g_{\rm WP}
(\frac{\partial}{\partial s},\frac{\partial}{\partial s})
-h^{1,2}(X)\,
\log\|\eta_{s}\|_{L^{2}}^{2}
\\
&=
O(\log(-\log|s|))
\end{aligned}
$$
by Propositions 4.4, 7.6 (3) and 7.8 (3). 
This proves the theorem.
\end{pf}

%%%%%%%%%%%%%%%%%%%%%%%%%%%%%%%%%%%%%%%%%%%%%%%%%%%%%%%%%%%
%
%                 Section 8
%
%%%%%%%%%%%%%%%%%%%%%%%%%%%%%%%%%%%%%%%%%%%%%%%%%%%%%%%%%%%

\section
{\bf The singularity of the BCOV invariant I -- 
the case of ODP}
\par
In Sect.\,8, we fix the following notation.
Let $\pi\colon{\cal X}\to S$ be a proper, surjective, 
flat holomorphic map from a compact, connected smooth 
K\"ahler fourfold to a compact Riemann surface. 
Let $\cal D$ be the discriminant locus and 
let $0\in{\cal D}$.
We assume that $X:=X_{0}$ is a Calabi-Yau threefold
with a unique ODP as its singular set 
satisfying $h^{2}(\Omega^{1}_{X})=1$. 
The deformation germ
$\pi\colon({\cal X},X)\to(S,0)$
is a smoothing of $X$, and a general fiber of $\pi$
is a smooth Calabi-Yau threefold.
We set $o:={\rm Sing}\,X$.
\par
Let 
${\frak p}\colon({\frak X},X)\to({\rm Def}(X),[X])$
be the Kuranishi family of $X$ with discriminant locus
$\frak D=[X]$. 
Since $h^{2}(\Omega^{1}_{X})=1$,
we have $\dim{\rm Def}(X)=1$.
By Proposition 2.8, 
${\frak X}$ is K\"ahler.
Let $g_{\frak X}$ be a K\"ahler metric on $\frak X$, 
and set 
$g_{{\frak X}/{\rm Def}(X)}:=
g_{\frak X}|_{T{\frak X}/{\rm Def}(X)}$.
\par
Let $\mu\colon(S,0)\to({\rm Def}(X),[X])$ be 
the holomorphic map that induces the family 
$\pi\colon({\cal X},X)\to(S,0)$
from the Kuranishi family.
By the local description (2.2), we have
${\cal O}_{{\cal X},o}\cong
{\Bbb C}\{z_{0},z_{1},z_{2},z_{3}\}/
(z_{0}^{2}+\cdots+z_{3}^{2}-\mu(t))$.
Since ${\cal X}$ is smooth,
$\frak D=\mu(0)$ is not a critical value of $\mu$,
and the morphism of germs 
$\mu\colon(S,0)\to({\rm Def}(X),[X])$
is an isomorphism. Hence there exist a neighborhood
$\cal U$ of $0\in S$ and an isomorphism of families
$f\colon
{\cal X}|_{\cal U}\cong{\frak X}|_{\mu({\cal U})}$. 
\par
Let $g_{\pi^{-1}(\cal U)}$ be the K\"ahler metric 
on $\pi^{-1}(\cal U)$ defined as
$$
g_{\pi^{-1}({\cal U})}=
f^{*}g_{\frak X}.
$$
Let
$g_{{\cal X}/S}$ be the Hermitian metric on 
$T{\cal X}/S|_{\pi^{-1}({\cal U})\setminus\Sigma_{\pi}}$ 
induced from $g_{\pi^{-1}(\cal U)}$. Then
$$
g_{{\cal X}/S}=f^{*}g_{{\frak X}/{\rm Def}(X)}.
$$
\par
Let
$\|\cdot\|^{2}_{\lambda
({\cal E}^{p}_{{\cal X}/S}),L^{2},g_{{\cal X}/S}}$
be the $L^{2}$-metric on the K\"ahler extension
$\lambda({\cal E}^{p}_{{\cal X}/S})|_{\cal U}$ 
with respect to $g_{{\cal X}/S}$. 
Since ${\cal F}_{{\cal X}/S}^{p}$ is acyclic
on ${\cal X}$ for $p=0,1$, we have the following
isomorphisms for $p=0,1$:
\begin{equation}
\lambda({\cal E}^{p}_{{\cal X}/S})|_{\cal U}
\cong
\mu^{*}\lambda(\Omega^{p}_{{\frak X}/{\rm Def}(X)}),
\qquad
\|\cdot\|_{L^{2},g_{{\cal X}/S}}=
\mu^{*}\|\cdot\|_{L^{2},g_{{\frak X}/{\rm Def}(X)}}.
\end{equation}
\par
Let $t$ be a local coordinate of $S$ centered at $0$.
Let $\sigma_{p}$ be a nowhere vanishing 
holomorphic section of the K\"ahler extension
$\lambda({\cal E}^{p}_{{\cal X}/S})$
near $0\in{\cal D}$.

\begin{thm}
The following formula holds as $t\to0$:
$$
(-1)^{p}\log\|\sigma_{p}(t)\|^{2}_{\lambda
({\cal E}^{p}_{{\cal X}/S}),L^{2},g_{{\cal X}/S}}
=
\left\{
\begin{array}{ll}
O(\log(-\log|t|))
&(p=0,1)
\\
-\log|t|^{2}+O(\log(-\log|t|))
&(p=2,3).
\end{array}
\right.
$$
\end{thm}

\begin{pf}
Let $p=0,1$.
Since $\mu\colon(S,0)\to({\rm Def}(X),[X])$
is an isomorphism, 
the assertion follows from Theorem 7.10 and (8.1).
\par
Let $p=2,3$. Recall that the canonical element 
${\bf 1}_{p,3-p}(X_{t})\in
\lambda(\Omega_{X_{t}}^{p})\otimes
\lambda(\Omega_{X_{t}}^{3-p})^{\lor}$ 
was defined in Subsection 3.3.
Let ${\bf 1}_{p,3-p,S^{o}}$ be the nowhere vanishing
holomorphic section of
$\lambda(\Omega_{{\cal X}^{o}/S^{o}}^{p})\otimes
\lambda(\Omega_{{\cal X}^{o}/S^{o}}^{3-p})^{\lor}$
defined by
$$
{\bf 1}_{p,3-p,S^{o}}(t):={\bf 1}_{p,3-p}(X_{t})
\in
\lambda(\Omega_{X_{t}}^{p})\otimes
\lambda(\Omega_{X_{t}}^{3-p})^{\lor},
\qquad
t\in S^{o}.
$$ 
Then
\begin{equation}
\|{\bf 1}_{p,3-p,S^{o}}(t)\|_{L^{2},g_{{\cal X}/S}}
=
\|{\bf 1}_{p,3-p,S^{o}}(t)\|_{Q,g_{{\cal X}/S}}
=
1,
\qquad
t\in S^{o}.
\end{equation}
by Proposition 3.4.
\par
By Theorem 5.11, we get
\begin{equation}
\begin{aligned}
\,&
\log
\|\sigma_{p}(t)\otimes\sigma_{3-p}(t)^{-1}
\|^{2}_{\lambda({\cal E}^{p}_{{\cal X}/S})\otimes
\lambda({\cal E}^{3-p}_{{\cal X}/S})^{\lor},
Q,g_{{\cal X}/S}}
\\
&=
(-1)^{3-p}\delta(3,p)\log|t|^{2}
+(-1)^{3}\cdot(-1)^{3-(3-p)}\delta(3,3-p)\log|t|^{2}
+O(1)
\\
&=
(-1)^{3-p}\log|t|^{2}+O(1),
\end{aligned}
\end{equation}
where we used the first identity of Lemma 5.12 
to get the last equality of (8.3).
\par
Set
$$
f_{p}(t):=
\frac{\sigma_{p}(t)\otimes\sigma_{3-p}(t)^{-1}}
{{\bf 1}_{p,3-p}(t)}\in{\cal O}(S^{o}).
$$
By (8.2), we get
\begin{equation}
\begin{aligned}
\,&
\|\sigma_{p}(t)\otimes\sigma_{3-p}(t)^{-1}
\|^{2}_{\lambda({\cal E}^{p}_{{\cal X}/S})\otimes
\lambda({\cal E}^{3-p}_{{\cal X}/S})^{\lor},Q,
g_{{\cal X}/S}}
\\
&
=
|f_{p}(t)|^{2}\cdot
\|{\bf 1}_{p,3-p}(t)
\|^{2}_{\lambda({\cal E}^{p}_{{\cal X}/S})\otimes
\lambda({\cal E}^{3-p}_{{\cal X}/S})^{\lor},Q,
g_{{\cal X}/S}}
\\
&
=
|f_{p}(t)|^{2}\cdot
\|{\bf 1}_{p,3-p}(t)
\|^{2}_{\lambda({\cal E}^{p}_{{\cal X}/S})\otimes
\lambda({\cal E}^{3-p}_{{\cal X}/S})^{\lor},L^{2},
g_{{\cal X}/S}}
\\
&
=
\|\sigma_{p}(t)\otimes\sigma_{3-p}(t)^{-1}
\|^{2}_{\lambda({\cal E}^{p}_{{\cal X}/S})\otimes
\lambda({\cal E}^{3-p}_{{\cal X}/S})^{\lor},
L^{2},g_{{\cal X}/S}},
\end{aligned}
\end{equation}
which, together with (8.3), yields that
\begin{equation}
\log
\|\sigma_{p}(t)\otimes\sigma_{3-p}(t)^{-1}
\|^{2}_{\lambda({\cal E}^{p}_{{\cal X}/S})\otimes
\lambda({\cal E}^{3-p}_{{\cal X}/S})^{\lor},
L^{2},g_{{\cal X}/S}}
=(-1)^{3-p}\,\log|t|^{2}+O(1).
\end{equation}
By Theorem 8.1 for $p=0,1$ and (8.4), 
we get
$$
\begin{aligned}
\,&
(-1)^{p}\log\|\sigma_{p}(t)\|^{2}_{
\lambda({\cal E}^{p}_{{\cal X}/S}),
L^{2},g_{{\cal X}/S}}
\\
&
=
(-1)^{p}\,\log
\|\sigma_{p}(t)\otimes\sigma_{3-p}(t)^{-1}
\|^{2}_{\lambda({\cal E}^{p}_{{\cal X}/S})\otimes
\lambda({\cal E}^{3-p}_{{\cal X}/S})^{\lor},
L^{2},g_{{\cal X}/S}}
\\
&
\quad+
(-1)^{p}\log\|\sigma_{3-p}(t)\|^{2}_{
\lambda({\cal E}^{3-p}_{{\cal X}/S}),
L^{2},g_{{\cal X}/S}}
\\
&
=
-\log|t|^{2}+O(\log(-\log|t|)).
\end{aligned}
$$
This proves the theorem for $p=2,3$.
\end{pf}

Let $\gamma_{t}$ be the K\"ahler form of
$g_{{\cal X}/S}|_{X_{t}}$.

\begin{thm}
The following formula holds as $t\to0$:
$$
\log\tau_{\rm BCOV}(X_{t})
=\frac{1}{6}\log|t|^{2}+O(\log(-\log|t|)).
$$
\end{thm}

\begin{pf}
By the definition of the BCOV torsion of 
$(X_{t},\gamma_{t})$, we have
\begin{equation}
\begin{aligned}
\log{\cal T}_{\rm BCOV}(X_{t},\gamma_{t})
&=
\sum_{p\geq0}(-1)^{p}p\,
\log\|\sigma_{p}(t)\|^{2}_{\lambda
({\cal E}^{p}_{{\cal X}/S}),Q,g_{{\cal X}/S}}
\\
&\quad
-\sum_{p\geq0}(-1)^{p}p\,
\log\|\sigma_{p}(t)\|^{2}_{\lambda
({\cal E}^{p}_{{\cal X}/S}),L^{2},g_{{\cal X}/S}}.
\\
&=
-\frac{19}{4}\log |t|^{2}
+\sum_{p=2}^{3}p\log|t|^{2}+O(\log(-\log|t|^{2}))
\\
&
=
\frac{1}{4}\log|t|^{2}+O(\log(-\log|t|^{2})),
\end{aligned}
\end{equation}
where we used Theorems 5.13 and 8.1 to get
the second equality.
Since
$$
\log{\rm Vol}(X_{t},\gamma_{t})=O(1),
\qquad
\log{\rm Vol}_{L^{2}}(H^{2}(X_{t},{\Bbb Z}),[\gamma_{t}])
=O(1),
$$
we deduce from Proposition 7.9 (2) and (8.6) that
$$
\begin{aligned}
\log\tau_{\rm BCOV}(X_{t})
&=
\log{\cal A}(X_{t},\gamma_{t})+
\log{\cal T}_{\rm BCOV}(X_{t},\gamma_{t})+O(1)
\\
&=
\frac{1}{6}\log|t|^{2}+O(\log(-\log|t|^{2})).
\end{aligned}
$$
This proves the theorem.
\end{pf}

%%%%%%%%%%%%%%%%%%%%%%%%%%%%%%%%%%%%%%%%%%%%%%%%%%%%%%%%%%%
%
%                 Section 9
%
%%%%%%%%%%%%%%%%%%%%%%%%%%%%%%%%%%%%%%%%%%%%%%%%%%%%%%%%%%%

\section
{\bf The singularity of the BCOV invariant II -- 
general degenerations}
\par
In Section 9, we fix the following notation:
Let $\cal X$ be an irreducible projective algebraic 
fourfold and let $S$ be a compact Riemann surface.
Let $\pi\colon{\cal X}\to S$ be a surjective, flat 
holomorphic map. Let ${\cal D}\subset S$ be
a reduced divisor and set
${\cal X}^{o}:={\cal X}\setminus\pi^{-1}({\cal D})$,
$S^{o}:=S\setminus{\cal D}$,
$\pi^{o}:=\pi|_{{\cal X}^{o}}$.
Let $0\in{\cal D}$, and let $(U,t)$ be 
a coordinate neighborhood of $S$ centered at $0$ 
such that
$U\setminus\{0\}\cong\varDelta^{*}$.
\par
In Section 9, we shall prove a generalization of 
Theorem 8.2.

\begin{thm}
If $\pi^{o}\colon{\cal X}^{o}\to S^{o}$ is a smooth
morphism whose fibers are Calabi-Yau threefolds,
then there exists $\alpha\in{\Bbb R}$ such that
as $t\to0$,
$$
\log\tau_{\rm BCOV}(X_{t})=
\alpha\,\log|t|^{2}+O(\log(-\log|t|^{2})).
$$
\end{thm}

First, we shall prove Theorem 9.1 when 
$\pi\colon{\cal X}\to S$ is a semi-stable family.
Then we shall reduce the general case 
to this particular case 
by the semi-stable reduction theorem of Mumford
\cite{Mumford73}. We set $D:=X_{0}$ in this section.

\subsection
{}{\bf The singularity of $L^{2}$ metrics for
semi-stable degenerations}
\par
In Subsections 9.1 and 9.2, we assume that $\cal X$
is smooth and that $D=X_{0}$ is a reduced divisor of
normal crossing, i.e., 
for every $x\in D$, there exist integers 
$\epsilon_{0},\epsilon_{1},\epsilon_{2},\epsilon_{3}
\in\{0,1\}$ and a coordinate neighborhood
$({\cal U},(z_{0},z_{1},z_{2},z_{3}))$ of $\cal X$
centered at $x$ such that
$$
\pi(z)=
z_{0}^{\epsilon_{0}}
z_{1}^{\epsilon_{1}}
z_{2}^{\epsilon_{2}}
z_{3}^{\epsilon_{3}},
\qquad
z\in{\cal U}.
$$
\par
Let $\Omega^{1}_{{\cal X}/S}(\log D)$ be the sheaf
of meromorphic $1$-forms on $\cal X$ with logarithmic
pole along $D$. Then 
$\Omega^{1}_{\cal X}(\log D)|_{{\cal X}\setminus D}
=\Omega^{1}_{\cal X}|_{{\cal X}\setminus D}$, 
and $\Omega^{1}_{\cal X}(\log D)|_{\cal U}$
is a free ${\cal O}_{\cal U}$-module generated by
$dz_{0}/z_{0}^{\epsilon_{0}},\,
dz_{1}/z_{1}^{\epsilon_{1}},\,
dz_{2}/z_{2}^{\epsilon_{2}},\,
dz_{3}/z_{3}^{\epsilon_{3}}$.
\par
Let $\Omega^{1}_{S}(\log0)$ be
the sheaf of meromorphic $1$-forms on $S$ with
logarithmic pole at $0$. Then
$\Omega^{1}_{S}(\log0)_{0}={\cal O}_{S,0}\,dt/t$.
We set 
$$
\Omega^{1}_{{\cal X}/S}(\log D):=
\Omega^{1}_{\cal X}(\log D)/\pi^{*}\Omega^{1}_{S}(\log0).
$$
See e.g. \cite[Sect.\,2]{Steenbrink77}, 
\cite[Chap.\,3, Sect.\,2]{Voisin99} for more details
about $\Omega^{1}_{{\cal X}/S}(\log D)$.
\par
Let $g_{\cal X}$ be a K\"ahler metric on $\cal X$ 
whose K\"ahler class is integral.
Let $\kappa\in H^{2}({\cal X},{\Bbb Z})$ be
the K\"ahler class of $g_{\cal X}$. 
We set
$g_{{\cal X}/S}:=g_{\cal X}|_{T_{\cal X}/S}$.

\subsubsection
{The canonical extension of the Hodge bundles}
\par
For the proof of Theorem 9.1, let us recall
some results of Schmid \cite{Schmid73} and
Steenbrink \cite{Steenbrink77}.
Set $U^{o}:=U\setminus\{0\}$. We fix $b\in U^{o}$
and set $W:=H^{m}(X_{b},{\Bbb C})$ and $l:=\dim W$.
\par
Let
${}^{o}{\bf H}^{m}:=R^{m}\pi_{*}{\Bbb C}\otimes_{\Bbb C}
{\cal O}_{U^{o}}$ and consider the Gauss-Manin
connection on ${}^{o}{\bf H}^{m}$. 
The canonical extension ${\bf H}^{m}$
of ${}^{o}{\bf H}^{m}$ from $U^{o}$ to $U$
is defined as follows:
Let $\{v_{1},\ldots,v_{l}\}$ be a basis of $W$,
and let $\gamma\in GL(W)$ be the Picard-Lefschetz
transformation.
There exists a Nilpotent operator $N\in{\rm End}(W)$ 
with $\gamma=\exp(N)$.
\par
Let 
$\psi\colon\widetilde{U^{o}}\ni z\to
\exp(2\pi\sqrt{-1}z)\in U^{o}$ be the universal covering.
Since ${}^{o}{\bf H}^{m}$ is flat, the vectors
$v_{i}$ extend to flat holomorphic sections 
${\bf v}_{i}\in\Gamma(\widetilde{U^{o}},
\psi^{*}({}^{o}{\bf H}^{m}))$, 
which induce an isomorphism
$\psi^{*}({}^{o}{\bf H}^{m})\cong
{\cal O}_{\widetilde{U^{o}}}\otimes_{\Bbb C}W$ of
flat bundles. Under this trivialization of 
$\psi^{*}({}^{o}{\bf H}^{m})$, we have
${\bf v}_{i}(z+1)=\gamma\cdot{\bf v}_{i}(z)$
for all $i$.
After Schmid \cite[pp.234-236]{Schmid73}, 
we define holomorphic frame fields of
$\psi^{*}({}^{o}{\bf H}^{m})$ by
\begin{equation}
{\bf s}_{i}(\exp 2\pi\sqrt{-1}z)
:=
\exp\left(-z\,N\right)\,{\bf v}_{i}(z)
=\sum_{k\geq0}\frac{1}{k!}(-z\,N)^{k}{\bf v}_{i}(z).
\end{equation}
Since
${\bf s}_{1},\ldots,{\bf s}_{l}\in
\Gamma(\widetilde{U^{o}},\psi^{*}({}^{o}{\bf H}^{m}))$
are invariant under the translation $z\to z+1$,
they descend to single-valued holomorphic 
frame fields of ${}^{o}{\bf H}^{m}$.
Then ${\bf H}^{m}$ 
is a locally free sheaf on $U$ defined as
${\bf H}^{m}:=
{\cal O}_{U}\,{\bf s}_{1}\oplus\cdots
\oplus{\cal O}_{U}\,{\bf s}_{l}$.
\par 
By Hodge theory, 
${}^{o}{\bf H}^{m}$ carries the {\it Hodge filtration}
$0\subset{}^{o}{\bf F}^{m}\subset\cdots\subset
{}^{o}{\bf F}^{1}\subset{}^{o}{\bf H}^{m}$
such that ${}^{o}{\bf F}^{p}$ is a holomorphic subbundle 
of ${}^{o}{\bf H}^{m}$ with
${}^{o}{\bf F}^{p}/{}^{o}{\bf F}^{p+1}\cong
R^{m-p}\pi_{*}\Omega^{p}_{{\cal X}/S}|_{U^{o}}$.
For $t\in U^{o}$, we have the natural identification
${}^{o}{\bf F}^{p}_{t}=
\bigoplus_{i\geq p}H^{m-i}(X_{t},\Omega^{i}_{X_{t}})$.
\par
By \cite[p.235]{Schmid73}, 
\cite[Th.\,2.11]{Steenbrink77},
\cite[pp.130 Cor.]{Zucker84}, 
the filtration $\{{}^{o}{\bf F}^{p}\}$ extends to 
a filtration $\{{\bf F}^{p}\}$ of
${\bf H}^{m}$ such that
${\bf F}^{p}/{\bf F}^{p+1}
\cong
R^{m-p}\pi_{*}\Omega^{p}_{{\cal X}/S}(\log D)|_{U}$.
Under this isomorphism, we have an identification
of holomorphic line bundles on $U$:
\begin{equation}
i_{p}:
(\det{\bf F}^{p})\otimes(\det{\bf F}^{p+1})^{-1}
\cong
\det R^{m-p}\pi_{*}\Omega^{p}_{{\cal X}/S}(\log D)|_{U}.
\end{equation}
\par
Since ${}^{o}{\bf H}^{m}_{t}=H^{m}(X_{t},{\Bbb C})$ 
for $t\in U^{o}$,
${}^{o}{\bf H}^{m}$ is equipped with the $L^{2}$-metric
$h_{R^{m}\pi_{*}{\Bbb C}}$ with respect to 
$g_{{\cal X}/S}$. 
Recall that the $C^{\infty}$ vector bundles
${\cal K}^{p,q}({\cal X}^{o}/U^{o})$ on $U^{o}$
were defined in Subsect.\,3.5.
Let $h_{{\bf F}^{p}}$ be the $L^{2}$-metric 
on ${}^{o}{\bf F}^{p}$ induced from 
$h_{R^{m}\pi_{*}{\Bbb C}}$ by the $C^{\infty}$ 
isomorphism
${}^{o}{\bf F}^{p}\cong\bigoplus_{i\geq p}
{\cal K}^{i,m-i}({\cal X}^{o}/U^{o})$.
By the definition of $L^{2}$-metrics, 
the isomorphism $i_{p}|_{U^{o}}$ 
induces an isometry of Hermitian line bundles
on $U^{o}$:
\begin{equation}
\left(
(\det{}^{o}{\bf F}^{p})\otimes
(\det{}^{o}{\bf F}^{p+1})^{-1},
\det h_{{\bf F}^{p}}\otimes(\det h_{{\bf F}^{p+1}})^{-1}
\right)
\cong
(\det R^{m-p}\Omega^{p}_{{\cal X}/S},
\|\cdot\|_{L^{2}}).
\end{equation}
\par
Recall that the K\"ahler operator 
$L\colon H^{m}(X_{t},{\Bbb C})\to
H^{m+2}(X_{t},{\Bbb C})$
with respect to $\kappa|_{X_{t}}$ 
was defined in Subsect.\,4.4.1.
Then $L$ induces a homomorphism of 
${\cal O}_{U}$-modules
$L\colon{\bf H}^{m}\to{\bf H}^{m+2}$. 
The primitive part of ${\bf H}^{m}$ is 
the holomorphic flat subbundle of ${\bf H}^{m}$
defined as
${\bf P}^{m}:={\bf H}^{m}\cap\ker L^{4-m}$.
The Picard-Lefschetz transformation $\gamma$
preserves ${\bf P}^{m}$.
If ${\bf s}_{i}\in\Gamma(U,{\bf P}^{m})$, 
there exists $k\in{\Bbb Z}$, $C\in{\Bbb R}$
by \cite[p.252 Th.\,6.6']{Schmid73} such that
\begin{equation}
\|{\bf s}_{i}(t)\|_{L^{2}}^{2}
\leq C\,(-\log|t|)^{k},
\qquad
t\in U^{o}.
\end{equation}

\subsubsection
{Singularities of the $L^{2}$-metrics: 
the case of canonical extension}
\par

\begin{lem}
Let $m=3$.
Let ${\bf f}_{p}$ be a nowhere vanishing holomorphic
section of $\det{\bf F}^{p}$ defined
on $U$.
Then there exists $c_{p}\in{\Bbb R}$ such that
as $t\to0$,
$$
\log\|{\bf f}_{p}(t)\|_{L^{2}}^{2}
=
c_{p}\,\log|t|^{2}+O(\log(-\log|t|)).
$$
\end{lem}

\begin{pf}
Since $m=3$, we have ${\bf H}^{3}={\bf P}^{3}$,
i.e., the groups $H^{3}(X_{t},{\Bbb C})$ are
primitive.
By (9.4), there exists a constant $C>0$ and
$l\in{\Bbb Z}$ such that
\begin{equation}
\lambda_{p}(t):=\|{\bf f}_{p}(t)\|_{L^{2}}^{2}
\leq C\,(-\log|t|)^{l},
\qquad
t\in U^{o}.
\end{equation}
We set $\lambda_{4}(t)=1$.
By Proposition 4.6 and (9.3), 
we get the following on $U^{o}$:
\begin{equation}
-dd^{c}(\log\lambda_{p}-\log\lambda_{p+1})
=
\left\{
\begin{array}{ll}
-\omega_{{\rm WP},{\cal X}^{o}/U^{o}}
& (p=0)
\\
-\omega_{{\rm H},{\cal X}^{o}/U^{o}}
+3\,\omega_{{\rm WP},{\cal X}^{o}/U^{o}}
&(p=1)
\\
\omega_{{\rm H},{\cal X}^{o}/U^{o}}
-3\,\omega_{{\rm WP},{\cal X}^{o}/U^{o}}
&(p=2)
\\
\omega_{{\rm WP},{\cal X}^{o}/U^{o}}
&(p=3).
\end{array}
\right.
\end{equation}
Since $\lambda_{p}\in L^{1}_{\rm loc}(U)$
by (9.5), the result follows from
Lemma 7.5 (1) and (9.6).
\end{pf}

\par
Let $\sigma_{p}$ be a nowhere vanishing holomorphic
section of $\lambda({\cal E}^{p}_{{\cal X}/S})$
near $0$.

\begin{prop}
There exists $\beta_{0}\in{\Bbb R}$ such that 
as $t\to0$:
$$
\log\|\sigma_{0}(t)\|^{2}_{\lambda
({\cal O}_{\cal X}),L^{2},g_{{\cal X}/S}}
=
\beta_{0}\,\log|t|^{2}+O(\log(-\log|t|)).
$$
\end{prop}

\begin{pf}
We may assume that 
$\sigma_{0}={\bf f}_{0}\otimes{\bf f}_{1}^{-1}$
under the isomorphism (9.2). 
Since (9.2) induces the isometry (9.3), 
the result follows from Lemma 9.2.
\end{pf}

\par
By \cite[Th.\,2.11]{Steenbrink77}, 
$R^{q}\pi_{*}\Omega^{1}_{{\cal X}/S}(\log D)$ is
locally free.
Set $r:={\rm rk}\,R^{q}\pi_{*}
\Omega^{1}_{{\cal X}/S}(\log D)$.
Let $e_{1}(t),\ldots,e_{r}(t)$ be a basis
of $R^{q}\pi_{*}\Omega^{1}_{{\cal X}/S}(\log D)$
as a free ${\cal O}_{U}$-module.

\begin{prop}
For $0\leq q\leq3$, there exists 
$\delta_{q}\in{\Bbb R}$ such that as $t\to0$,
$$
\log \|e_{1}(t)\wedge\cdots\wedge e_{r}(t)\|^{2}_{
\det R^{q}\pi_{*}\Omega^{1}_{{\cal X}/S}(\log D),
L^{2},g_{{\cal X}/S}}
=
\delta_{q}\,\log|t|^{2}+O(\log(-\log|t|)).
$$
\end{prop}

\begin{pf}
Since $r=0$ when $q=0,3$, it suffices to prove
the cases $q=1,2$.
\par{\bf (Case 1)}
Let $q=2$. There exists a nowhere vanishing
holomorphic function $h(t)$ on $U$ such that
$e_{1}(t)\wedge\cdots\wedge e_{r}(t)=
h(t)\,{\bf f}_{1}(t)\otimes{\bf f}_{2}(t)^{-1}$
under the isomorphism (9.2). Since (9.2) induces
the isometry (9.3), the result follows 
from Lemma 9.2.
\par{\bf (Case 2)}
Let $q=1$. When $m=2$, we have
${\bf H}^{2}={\bf F}^{1}$.
Hence $r=l$.
Identify the integral K\"ahler class $\kappa$
on ${\cal X}$ with the corresponding flat section
of ${\bf H}^{2}$. Then ${\bf P}^{m}$
and ${\cal O}_{U}\,\kappa$ are holomorphic flat
subbundles of ${\bf H}^{m}$ preserved by
the Picard-Lefschetz transformation $\gamma$.
Hence we have a decomposition
${\bf H}^{2}={\bf P}^{2}\oplus{\cal O}_{U}\,\kappa$
of $\gamma$-invariant flat bundles on $U$.
Choose $v_{1}=\kappa_{b}$ and
$v_{2},\ldots,v_{l}\in{\bf P}^{2}_{b}\cap 
H^{2}(X_{b},{\Bbb Z})/{\rm Torsion}$ 
in Subsect.\,9.1.1. Then
${\bf s}_{1}=\kappa$ and
${\bf P}^{m}={\cal O}_{U}\,{\bf s}_{2}
\oplus\cdots\oplus{\cal O}_{U}\,{\bf s}_{l}$.
Since ${\bf v}_{1}(z),\ldots,{\bf v}_{l}(z)$ are
identified with $v_{1},\ldots,v_{l}$ via the
$C^{\infty}$ trivialization
${\cal X}^{o}\times_{U^{o}}\widetilde{U}^{o}
\cong X_{b}\times\widetilde{U}^{o}$, we get
by Definition 4.11 and Lemma 4.12
\begin{equation}
\|{\bf v}_{1}(z)\wedge\cdots\wedge{\bf v}_{l}(z)
\|_{L^{2},\kappa}^{2}=
{\rm Vol}_{L^{2}}(H^{2}(X_{b},{\Bbb Z}),\kappa_{b}),
\qquad
\forall\,z\in\widetilde{U}^{o}.
\end{equation}
Since $N$ is nilpotent and hence
$\det\exp(-zN)=1$ for all $z\in\widetilde{U}^{o}$, 
we get
\begin{equation}
\begin{aligned}
{\bf s}_{1}(e^{2\pi\sqrt{-1}z})\wedge\cdots\wedge
{\bf s}_{l}(e^{2\pi\sqrt{-1}z})
&
=
\exp(-z\,N)\,{\bf v}_{1}(z)\wedge\cdots\wedge
\exp(-z\,N)\,{\bf v}_{l}(z)
\\
&
=
\det\exp(-zN)\cdot{\bf v}_{1}(z)\wedge\cdots\wedge
{\bf v}_{l}(z)
\\
&
=
{\bf v}_{1}(z)\wedge\cdots\wedge{\bf v}_{l}(z).
\end{aligned}
\end{equation}
By (9.7), (9.8), we get for all $t\in U^{o}$:
\begin{equation}
\|{\bf s}_{1}(t)\wedge\cdots\wedge{\bf s}_{l}(t)
\|_{L^{2},\kappa}^{2}
=
{\rm Vol}_{L^{2}}(H^{2}(X_{b},{\Bbb Z}),\kappa_{b}).
\end{equation}
Since $\{{\bf s}_{1}(t),\ldots,{\bf s}_{l}(t)\}$
is a basis of 
$R^{1}\pi_{*}\Omega^{1}_{{\cal X}/S}(\log D)$
as a free ${\cal O}_{S}$-module, the result follows
from (9.9).
This completes the proof.
\end{pf}

\subsubsection
{Comparison of the K\"ahler extension and
the canonical extension}

\begin{prop}
There exists $\beta_{1}\in{\Bbb R}$ such that
$$
\log\|\sigma_{1}(t)\|^{2}_{\lambda
({\Omega}^{1}_{{\cal X}/S}),L^{2},g_{{\cal X}/S}}
=
\beta_{1}\,\log|t|^{2}+O(\log(-\log|t|))
\qquad
(t\to0).
$$
\end{prop}

\begin{pf}
Consider the natural injection
$0\to\Omega^{1}_{{\cal X}/S}\to
\Omega^{1}_{{\cal X}/S}(\log D)$, and set 
$Q:=\Omega^{1}_{{\cal X}/S}(\log D)/
\Omega^{1}_{{\cal X}/S}$. Then $Q$ is a torsion
sheaf on $\cal X$ whose support is contained in
${\rm Sing}(D)$. 
Consider the long exact sequence of direct image 
sheaves induced by the short exact sequence of sheaves 
$0\to\Omega^{1}_{\cal X/S}\to
\Omega^{1}_{{\cal X}/S}(\log D)\to Q\to0$ on $\cal X$:
$$
R^{q-1}\pi_{*}\Omega^{1}_{{\cal X}/S}(\log D)\to 
R^{q-1}\pi_{*}Q\to 
R^{q}\pi_{*}\Omega^{1}_{{\cal X}/S}\to
R^{q}\pi_{*}\Omega^{1}_{{\cal X}/S}(\log D)
\to R^{q}\pi_{*}Q.
$$
Since $R^{q}\pi_{*}Q$ is a torsion sheaf on $U$
supported at $\{0\}$ for all $q$, there exist
torsion sheaves $M_{q}$, $N_{q}$ on $U$ 
supported at $\{0\}$ and an exact sequence 
of coherent sheaves on $U$:
\begin{equation}
\begin{CD}
0\to M_{q}\to
R^{q}\pi_{*}\Omega^{1}_{{\cal X}/S}
@>j>>
R^{q}\pi_{*}\Omega^{1}_{{\cal X}/S}(\log D)
\to N_{q}\to 0.
\end{CD}
\end{equation}
Since $U\cong\varDelta$ and hence
${\cal O}_{U,t}$ is a discrete valuation ring 
for all $t\in U$, the image 
$j(R^{q}\pi_{*}\Omega^{1}_{{\cal X}/S})$ is
a locally free submodule of
$R^{q}\pi_{*}\Omega^{1}_{{\cal X}/S}(\log D)$.
Hence 
$(R^{q}\pi_{*}\Omega^{1}_{{\cal X}/S})_{\rm tor}$,
the torsion part of
$R^{q}\pi_{*}\Omega^{1}_{{\cal X}/S}$, is contained
in $\ker j$. Since $M_{q}\subset
(R^{q}\pi_{*}\Omega^{1}_{{\cal X}/S})_{\rm tor}$,
we have
\begin{equation}
M_{q}=(R^{q}\pi_{*}\Omega^{1}_{{\cal X}/S})_{\rm tor}.
\end{equation}
\par
Since 
$N_{q}=R^{q}\pi_{*}\Omega^{1}_{{\cal X}/S}(\log D)/
j(R^{q}\pi_{*}\Omega^{1}_{{\cal X}/S})$ is
a torsion sheaf, there exist integers
$\nu_{1},\ldots,\nu_{r}\geq0$ such that
$N_{q}\cong
{\Bbb C}\{t\}/(t^{\nu_{1}})
\,\oplus\cdots\oplus\,
{\Bbb C}\{t\}/(t^{\nu_{r}})$ and
$j(R^{q}\pi_{*}\Omega^{1}_{{\cal X}/S})=
{\cal O}_{U}\,t^{\nu_{1}}e_{1}(t)\oplus\cdots\oplus
{\cal O}_{U}\,t^{\nu_{r}}e_{r}(t)$. 
Hence
\begin{equation}
\det\,j(R^{q}\pi_{*}\Omega^{1}_{{\cal X}/S})=
{\cal O}_{U}\cdot
t^{\nu_{1}}e_{1}(t)\wedge\cdots\wedge
t^{\nu_{r}}e_{r}(t).
\end{equation}
\par
By \cite[p.110 3. Proof of the theorem]{Banica76},
there exists a complex of locally free sheaves
of finite rank on $U$
$$
\CD
E_{\bullet}
\colon
0
\to 
E_{-1}
@>v_{-1}>> 
E_{0}
@>v_{0}>>
\cdots 
@>v_{k-1}>> 
E_{k}
\to
0
\endCD
$$
such that $R^{q}\pi_{*}\Omega^{1}_{{\cal X}/S}$
is the $q$-th cohomology sheaf of $E_{\bullet}$, i.e.,
$R^{q}\pi_{*}\Omega^{1}_{{\cal X}/S}\cong
H^{q}(E_{\bullet})$ for all $q\geq0$.
Since $U\cong\varDelta$, $\ker v_{q}\subset E_{q}$ 
and ${\rm Im}\,v_{q}\subset E_{q+1}$ are
locally free sheaves on $U$ for all $q\geq-1$. 
Let $\xi_{q}$ be the inverse image of 
$(R^{q}\pi_{*}\Omega^{1}_{{\cal X}/S})_{\rm tor}$
by the natural surjection 
$\ker v_{q}\to R^{q}\pi_{*}\Omega^{1}_{{\cal X}/S}$,
and set $\eta_{q}:={\rm Im}\,v_{q-1}$.
There exists an exact sequence of coherent sheaves
on $U$
$$
\CD
0
\to
\eta_{q}
@>\varphi_{q}>>
\xi_{q}
\to
(R^{q}\pi_{*}\Omega^{1}_{{\cal X}/S})_{\rm tor}
\to
0
\endCD
$$
such that $\eta_{q}$, $\xi_{q}$ are locally free
with equal rank.
Under the canonical isomorphism
$\det(R^{q}\pi_{*}\Omega^{1}_{{\cal X}/S})_{\rm tor}
\cong\det\xi_{q}\otimes(\det\eta_{q})^{-1}$,
the canonical section
$\det\varphi_{q}\in
H^{0}(U,\det\xi_{q}\otimes(\det\eta_{q})^{-1})$
induces the trivialization
$\det(R^{q}\pi_{*}\Omega^{1}_{{\cal X}/S})_{\rm tor}
\cong{\cal O}_{U}$ on $U^{o}$ by 
\cite[pp.118, Proof of Lemma 1, First Case]{Soule92}:
\begin{equation}
\det(R^{q}\pi_{*}\Omega^{1}_{{\cal X}/S})_{\rm tor}
\ni\det\varphi_{q}\to 1\in{\cal O}_{U}.
\end{equation}
\par
Since
$\det R^{q}\pi_{*}\Omega^{1}_{{\cal X}/S}\cong
\det j(R^{q}\pi_{*}\Omega^{1}_{{\cal X}/S})\otimes
\det(R^{q}\pi_{*}\Omega^{1}_{{\cal X}/S})_{\rm tor}$
by (9.10) and (9.11),
we deduce from (9.12), (9.13) that
the following expression $s_{1,q}$ is a holomorphic 
section of 
$\det R^{q}\pi_{*}\Omega^{1}_{{\cal X}/S}$:
$$
s_{1,q}(t):=
(t^{\nu_{1}}e_{1}(t)\wedge\cdots\wedge
t^{\nu_{r}}e_{r}(t))\otimes\det\varphi_{q}(t).
$$
Since $s_{1,q}(t)|_{U^{o}}$ 
is identified with the section 
$t^{\nu_{1}}e_{1}(t)\wedge\cdots\wedge
t^{\nu_{r}}e_{r}(t)|_{U^{o}}$ 
under the identification
$\det R^{q}\pi_{*}\Omega^{1}_{{\cal X}/S}|_{U^{o}}
\cong
\det j(R^{q}\pi_{*}\Omega^{1}_{{\cal X}/S})|_{U^{o}}$ 
induced by (9.13), we deduce from 
Proposition 9.4 that for $t\in U^{o}$,
\begin{equation}
\begin{aligned}
\log\|s_{1,q}(t)\|^{2}_{L^{2},g_{{\cal X}/S}}
&
=
\log\|t^{\nu_{1}}e_{1}(t)\wedge\cdots\wedge
t^{\nu_{r}}e_{r}(t)\|^{2}_{L^{2},g_{{\cal X}/S}}
\\
&
=
\dim_{\Bbb C}N_{q}\log|t|^{2}+
\log\|e_{1}(t)\wedge\cdots\wedge
e_{r}(t)\|^{2}_{L^{2},g_{{\cal X}/S}}
\\
&
=
(\dim_{\Bbb C}N_{q}+\delta_{q})\,\log|t|^{2}+
O(\log(-\log|t|)).
\end{aligned}
\end{equation}
Since  
$\det\varphi_{q}$ vanishes at $t=0$ with multiplicity 
$\dim_{\Bbb C}M_{q}$, 
$\sigma_{1,q}(t):=t^{-\dim_{\Bbb C}M_{q}}\,s_{1,q}(t)$
is a nowhere vanishing holomorphic section of
$\det R^{q}\pi_{*}\Omega^{1}_{{\cal X}/S}$. 
By (9.14), we get
\begin{equation}
\log\|\sigma_{1,q}(t)\|^{2}_{L^{2},g_{{\cal X}/S}}
=
(\dim_{\Bbb C}N_{q}+\delta_{q}-\dim_{\Bbb C}M_{q})\,
\log|t|^{2}+O(\log(-\log|t|)).
\end{equation}
The result follows from (9.15).
This completes the proof of Proposition 9.5.
\end{pf}

\begin{prop}
Let $p=2,3$.
There exists $\beta_{p}\in{\Bbb R}$ 
such that as $t\to0$,
$$
\log\|\sigma_{p}(t)\|^{2}_{\lambda
({\cal E}^{p}_{{\cal X}/S}),L^{2},g_{{\cal X}/S}}
=
\beta_{p}\log|t|^{2}+O(\log(-\log|t|)).
$$
\end{prop}

\begin{pf}
We keep the notation in Section 8, Proof of 
Theorem 8.1. 
By Theorem 5.4, there exists $a_{p}\in{\Bbb Q}$
such that
\begin{equation}
\log
\|\sigma_{p}(t)\otimes\sigma_{3-p}(t)^{-1}
\|^{2}_{\lambda({\cal E}^{p}_{{\cal X}/S})\otimes
\lambda({\cal E}^{3-p}_{{\cal X}/S})^{\lor},
Q,g_{{\cal X}/S}}
=
a_{p}\log|t|^{2}+O(1).
\end{equation}
By the same argument as in the proof 
of Theorem 8.1 (8.4) using (9.16) 
in stead of (8.3), we get
$$
\log\|\sigma_{p}(t)\otimes\sigma_{3-p}(t)^{-1}
\|^{2}_{\lambda({\cal E}^{p}_{{\cal X}/S})\otimes
\lambda({\cal E}^{3-p}_{{\cal X}/S})^{\lor},L^{2},
g_{{\cal X}/S}}
=
a_{p}\,\log|t|^{2}+O(1),
$$
which, together with Propositions 9.3 and 9.5,
yields the existence of $\beta_{p}\in{\Bbb R}$
such that
$$
\log\|\sigma_{p}(t)\|^{2}_{
\lambda({\cal E}^{p}_{{\cal X}/S}),L^{2},
g_{{\cal X}/S}} 
=
\beta_{p}\log|t|^{2}+O(\log(-\log|t|)).
$$
This proves the proposition.
\end{pf}

\subsection
{}
{\bf Proof of Theorem 9.1: 
the case of semi-stable degenerations}
\par
Let $\gamma_{t}$ be the K\"ahler form of
$g_{{\cal X}/S}|_{X_{t}}$.
By the definition of the BCOV torsion of 
$(X_{t},\gamma_{t})$, we have
$$
\log{\cal T}_{\rm BCOV}(X_{t},\gamma_{t})
=\sum_{p}(-1)^{p}p\,\{
\log\|\sigma_{p}(t)\|^{2}_{\lambda
({\cal E}^{p}_{{\cal X}/S}),Q,g_{{\cal X}/S}}-
\log\|\sigma_{p}(t)\|^{2}_{\lambda
({\cal E}^{p}_{{\cal X}/S}),L^{2},g_{{\cal X}/S}}\}.
$$
By Theorem 5.4 and Propositions 9.3, 9.5, 9.6, 
there exists $a\in{\Bbb R}$ such that 
\begin{equation}
\log{\cal T}_{\rm BCOV}(X_{t},\gamma_{t})
=
a\,\log|t|^{2}+O(\log(-\log|t|^{2})).
\end{equation}
Since the K\"ahler class of $g_{\cal X}$ is integral,
there exist positive constants $A,B\in{\Bbb Q}$
by Lemma 4.12 such that for all $t\in U^{o} $,
\begin{equation}
\log{\rm Vol}(X_{t},\gamma_{t})=A,
\qquad
\log{\rm Vol}_{L^{2}}
(H^{2}(X_{t},{\Bbb Z}),[\gamma_{t}])=B.
\end{equation}
By Proposition 7.9 (1), 
there exists $\epsilon\in{\Bbb R}$
such that
\begin{equation}
\log{\cal A}(X_{t},\gamma_{t})
=
\epsilon\,\log|t|^{2}+O(\log(-\log|t|^{2})).
\end{equation}
By (9.17), (9.18), (9.19), we get
$$
\begin{aligned}
\log\tau_{\rm BCOV}(X_{t})
&=
\log{\cal A}(X_{t},\gamma_{t})+
\log{\cal T}_{\rm BCOV}(X_{t},\gamma_{t})+O(1)
\\
&=
(a+\epsilon)\,\log|t|^{2}
+O(\log(-\log|t|^{2})).
\end{aligned}
$$
This proves the theorem.
\qed

\subsection
{}{\bf Proof of Theorem 9.1: general cases}
\par
In Subsection 9.3, we only assume that 
$\pi^{o}\colon{\cal X}^{o}\to S^{o}$ is a smooth
morphism whose fibers are Calabi-Yau threefolds.
\par
By the semi-stable reduction theorem
\cite[Chap.\,II]{Mumford73}, 
there exist a pointed projective curve $(B,o)$,
a finite surjective holomorphic map 
$f\colon(B,o)\to(S,0)$, and a holomorphic surjection 
$p\colon{\cal Y}\to B$ from a projective fourfold 
$\cal Y$ to $B$ satisfying the following conditions:
\newline{(i) }
Let $V$ be the component of $f^{-1}(U)$ containing $o$.
Then $f\colon V\setminus\{o\}\to U\setminus\{0\}$ 
is an isomorphism;
\newline{(ii) }
Set $U^{*}=U\setminus\{0\}$ and $V^{*}=V\setminus\{o\}$.
Then $p|_{V^{*}}\colon{\cal Y}|_{V^{*}}\to V^{*}$ is
induced from 
$\pi|_{U^{*}}\colon{\cal X}|_{U^{*}}\to U^{*}$
by $f|_{V^{*}}$;
\newline{(iii) }
$\cal Y$ is smooth, and
$Y_{o}$ is a reduced divisor of normal
crossing.
\par
Let $b$ be the coordinate on $V$ centered at $o$.
By condition (i), we may assume that there exists 
$\nu\in{\Bbb N}$ such that $f^{*}t=b^{\nu}$. 
Let $\tau_{U^{*}}$ and $\tau_{V^{*}}$
be the functions on $U^{*}$ and $V^{*}$ defined by
$$
\tau_{U^{*}}(t):=\tau_{\rm BCOV}(X_{t}),
\qquad
\tau_{V^{*}}(b):=\tau_{\rm BCOV}(Y_{b})
$$
for $t\in U^{*}$ and $b\in V^{*}$, respectively.
By condition (ii) and Theorem 4.16, we get
\begin{equation}
\tau_{V^{*}}=f^{*}\tau_{U^{*}}
\end{equation}
We can apply Theorem 9.1
to the family $p|_{V}\colon{\cal Y}|_{V}\to V$
by condition (iii), so that there exists
$\alpha\in{\Bbb R}$ such that as $b\to0$,
\begin{equation}
\log\tau_{V^{*}}(b)=\alpha\,\log|b|^{2}+
O(\log(-\log|b|)).
\end{equation}
Since $b=t^{\nu}$, the desired formula follows from
(9.20) and (9.21). 
This completes the proof of Theorem 9.1.
\qed

%%%%%%%%%%%%%%%%%%%%%%%%%%%%%%%%%%%%%%%%%%%%%%%%%%%%%%%%%%%
%
%                 Section 10
%
%%%%%%%%%%%%%%%%%%%%%%%%%%%%%%%%%%%%%%%%%%%%%%%%%%%%%%%%%%%

\section
{\bf The curvature current of the BCOV invariant}
\par
Following \cite[Sect.\,7]{Yoshikawa04}, 
we extend Theorem 4.14 to the Kuranishi space 
of Calabi-Yau threefold with a unique ODP
as its singular set.

\subsection
{}
{\bf  The curvature current of $\tau_{\rm BCOV}$:
general cases}
\par
In Subsection 10.1, we fix the following notation.
Let $\cal X$ be an irreducible projective algebraic 
fourfold and let $S$ be a compact Riemann surface.
Let $\pi\colon{\cal X}\to S$ be a surjective, flat 
holomorphic map. Let ${\cal D}\subset S$ be
a reduced divisor and set
${\cal X}^{o}:={\cal X}\setminus\pi^{-1}({\cal D})$,
$S^{o}:=S\setminus{\cal D}$,
$\pi^{o}:=\pi|_{{\cal X}^{o}}$.
We assume that the fibers of 
$\pi^{o}\colon{\cal X}^{o}\to S^{o}$ are Calabi-Yau
threefolds with $h^{2}(\Omega^{1}_{X_{s}})=1$
for $s\in S^{o}$.
Let $\chi(X)$ denote
the topological Euler number of $X_{s}$, $s\in S^{o}$.
\par
Let $\Omega_{{\rm WP},{\cal X}/S}$ and
$\Omega_{{\rm H},{\cal X}/S}$ be the trivial extensions 
of the Weil-Petersson form and the Hodge form from 
$S^{o}$ to $S$ (cf. Proposition 7.3 and
Definition 7.4). Then the $(1,1)$-currents
$\Omega_{{\rm WP},{\cal X}/S}$ and
$\Omega_{{\rm H},{\cal X}/S}$ are positive.
\par
Let $0\in{\cal D}$ and let $(U,t)$ be a coordinate
neighborhood of $S$ centered at $0$.
By Eq.\,(7.7), there exist subharmonic functions
$\varphi$ and $\theta$ on $U$ satisfying
the following equations of currents on $U$:
\begin{equation}
dd^{c}\varphi=\Omega_{{\rm WP},{\cal X}/S}|_{U},
\qquad
dd^{c}\theta=\Omega_{{\rm H},{\cal X}/S}|_{U}.
\end{equation}
\par
As in Subsection 4.4.2, we define a function on $S$
by
$$
\tau_{\rm BCOV}({\cal X}/S)(t)
:=
\tau_{\rm BCOV}(X_{t}),
\qquad
t\in S.
$$
By Theorems 4.14 and 9.1, 
$\log\tau_{\rm BCOV}({\cal X}/S)
\in C^{\infty}(S^{o})\cap L^{1}(S)$.

\begin{thm}
Set
$$
a
:=
\lim_{t\to0}
\frac{\log\tau_{\rm BCOV}({\cal X}/S)|_{U}(t)}
{\log|t|^{2}}
\in{\Bbb R}.
$$
Then the following equation of currents on $U$ holds:
$$
dd^{c}\log\tau_{\rm BCOV}({\cal X}/S)
=
-\frac{\chi(X)}{12}\,
\Omega_{{\rm WP},{\cal X}/S}
-\Omega_{{\rm H},{\cal X}/S}
+a\,\delta_{0}.
$$
\end{thm}

\begin{pf}
Identify $U$ with $\varDelta$ in what follows.
By Theorem 9.1, there exists a positive constant $K$ 
such that
\begin{equation}
\left|
\log\tau_{\rm BCOV}({\cal X}/S)(t)-a\,\log|t|^{2}
\right|
\leq 
K\,\log(-\log|t|),
\qquad
t\in\varDelta(1/2)^{*}.
\end{equation}
\par
For $t\in\varDelta(1/2)^{*}$, set 
$$
P(t):=
\left(
\log\tau_{\rm BCOV}({\cal X}/S)(t)-a\,\log|t|^{2}
\right)
+\frac{\chi(X)}{12}\,\varphi(t)
+\theta(t).
$$
Then 
$P(t)\in C^{\infty}(\varDelta(1/2)^{*})$.
By (7.10) and (10.2), 
there exists a positive constant $L$ such that
\begin{equation}
|P(t)|\leq L\,\log(-\log|t|^{2}),
\qquad
t\in\varDelta(1/2)^{*}.
\end{equation}
Since $P$ is harmonic on $\varDelta(1/2)^{*}$ by
Theorem 4.14 and (10.1), 
we deduce from Lemma 7.1 (3)
that $P$ extends to a harmonic function 
on $\varDelta(1/2)$.
Since $P$ is harmonic on $\varDelta(1/2)$, 
it follows from (7.10) that
\begin{equation}
\log\tau_{\rm BCOV}({\cal X}/S)
=
a\,\log|t|^{2}
-\frac{\chi(X)}{12}\,\varphi
-\theta+P
\in L^{1}_{\rm loc}(\varDelta(1/2)).
\end{equation}
Since $dd^{c}P=0$ on $\varDelta$, 
Eq.\,(10.4), together with (10.1), yields the assertion.
\end{pf}

\subsection
{}{\bf The curvature current of $\tau_{\rm BCOV}$: 
the case of Kuranishi families}

\par
In Subsection 10.2, we fix the following notation:
Let $X$ be a smoothable Calabi-Yau threefold with
only one ODP as its singular set. Let ${\rm Def}(X)$
be the Kuranishi space of $X$ with discriminant locus
$\frak D$, and let 
${\frak p}\colon({\frak X},X)\to({\rm Def}(X),[X])$
be the Kuranishi family of $X$. 
Assume that
$\dim{\rm Def}(X)=h^{2}(\Omega^{1}_{X})=1$.
Let $s$ be a coordinate on ${\rm Def}(X)$ 
such that ${\frak D}={\rm div}(s)$.
We identify ${\rm Def}(X)$ with the disc
$\varDelta$ equipped with the coordinate $s$. 
Then ${\rm Def}(X)\setminus{\frak D}\cong
\varDelta^{*}$.
\par
Let $\Omega_{\rm WP}$ and $\Omega_{\rm H}$ be
the trivial extensions of the Weil-Petersson form
and the Hodge form from 
${\rm Def}(X)\setminus{\frak D}$ to
${\rm Def}(X)$.
Let $\chi(X_{\rm gen})$ denote
the topological Euler number of a general fiber
of the Kuranishi family.

\begin{thm}
The function $\log\tau_{\rm BCOV}$ is 
locally integrable on ${\rm Def}(X)$, and
the following equation of currents on ${\rm Def}(X)$
holds:
$$
dd^{c}\log\tau_{\rm BCOV}
=
-\frac{\chi(X_{\rm gen})}{12}\,\Omega_{\rm WP}
-\Omega_{\rm H}
+\frac{1}{6}\,\delta_{\frak D}.
$$
\end{thm}

\begin{pf}
By Proposition 2.8, 
there exist a pointed projective curve $(B,0)$, 
a projective fourfold $\frak Z$, 
and a surjective, proper, flat holomorphic map
$f\colon{\frak Z}\to B$ such that the deformation germ
$f\colon({\frak Z},f^{-1}(0))\to(B,0)$ is isomorphic to
the Kuranishi family
${\frak p}\colon({\frak X},X)\to({\rm Def}(X),[X])$.
Since ${\rm Def}(X)$ is smooth at $[X]$,
so is $B$ at $0$.
By Theorem 9.1, we get $\log\tau_{\rm BCOV}
\in L^{1}_{\rm loc}({\rm Def}(X))$.
Let $\gamma:=\lim_{t\to0}\log\tau_{\rm BCOV}(X_{t})/
\log|t|^{2}$.
Since $\gamma=\frac{1}{6}$ by Theorem 8.2, 
the result follows from Theorem 10.1.
\end{pf}

\subsection
{}{\bf The curvature current of $\tau_{\rm BCOV}$:
the case of induced families }
\par
We keep the notation in Subsection 10.2.
Let $\mu\colon(\varDelta,0)\to({\rm Def}(X),[X])$ be
a holomorphic map and 
let $\pi\colon{\cal X}\to\varDelta$ be the family of 
Calabi-Yau threefolds induced from the Kuranishi family 
${\frak p}\colon({\frak X},X)\to({\rm Def}(X),[X])$ 
by $\mu$. Notice that ${\cal X}$ is singular
if $0$ is a critical point of $\mu$.

\begin{thm}
The function $\log\tau_{\rm BCOV}({\cal X}/\varDelta)$ 
lies in $L^{1}_{\rm loc}(\varDelta)$, and
the following equation of currents on $\varDelta$ holds:
$$
dd^{c}\log\tau_{\rm BCOV}({\cal X}/\varDelta)
=
-\frac{\chi(X_{\rm gen})}{12}\,
\Omega_{{\rm WP},{\cal X}/\varDelta}
-\Omega_{{\rm H},{\cal X}/\varDelta}
+\frac{1}{6}\,\delta_{\mu^{*}{\frak D}}.
$$
\end{thm}

\begin{pf}
Let $f\in{\cal O}_{{\rm Def}(X),[X]}$ be 
such that $\frak D={\rm div}(f)$.
Let $\Omega_{\rm WP}$ and $\Omega_{\rm H}$
be the trivial extensions of the Weil-Petersson
and the Hodge forms on ${\rm Def}(X)$, respectively.
As in Eq.\,(7.7),
let $\varphi$ and $\theta$ be the subharmonic
functions on ${\rm Def}(X)$ 
with $\Omega_{\rm WP}=dd^{c}\varphi$ and 
$\Omega_{\rm H}=dd^{c}\theta$.
Then $\mu^{*}\varphi$ and $\mu^{*}\theta$ are
subharmonic functions on $\varDelta$ with
\begin{equation}
dd^{c}(\mu^{*}\varphi)|_{\varDelta^{*}}
=\omega_{{\rm WP},{\cal X}/\varDelta},
\qquad
dd^{c}(\mu^{*}\theta)|_{\varDelta^{*}}
=\omega_{{\rm H},{\cal X}/\varDelta}.
\end{equation}
\par
After shrinking ${\rm Def}(X)$ if necessary, 
we may assume by (7.10) the existence
of constants $C_{0},C_{1}>0$ with
\begin{equation}
-C_{0}\,\log(-\log|f|^{2})
\leq
\varphi|_{{\rm Def}(X)\setminus{\frak D}}
\leq
C_{1},
\qquad
-C_{0}\,\log(-\log|f|^{2})
\leq
\theta|_{{\rm Def}(X)\setminus{\frak D}}
\leq
C_{1}.
\end{equation}
Since 
$\mu^{-1}({\frak D})\cap\varDelta=\{0\}$,
there exist a positive integer $k$ and
a nowhere vanishing holomorphic function 
$\varepsilon(s)\in{\cal O}(\varDelta)$ with
\begin{equation}
\mu^{*}f(s)=s^{k}\,\varepsilon(s).
\end{equation}
After shrinking $\varDelta$ if necessary,
the following inequality holds by (10.6)
\begin{equation}
-C_{2}\,\log(-\log|s|^{2})
\leq
\mu^{*}\varphi|_{\varDelta^{*}}
\leq
C_{1},
\qquad
-C_{2}\,\log(-\log|s|^{2})
\leq
\mu^{*}\theta|_{\varDelta^{*}}
\leq
C_{1},
\end{equation}
where $C_{2}>0$ is a constant.
By (10.5), (10.8) and Lemma 7.5 (2), 
we get the following equations of currents 
on $\varDelta$:
\begin{equation}
\Omega_{{\rm WP},{\cal X}/\varDelta}
=dd^{c}(\mu^{*}\varphi),
\qquad
\Omega_{{\rm H},{\cal X}/\varDelta}
=dd^{c}(\mu^{*}\theta).
\end{equation}
\par
By (10.4) and Theorem 10.2, there exists 
a harmonic function $P$ on ${\rm Def}(X)$ such that
$$
\log\tau_{\rm BCOV}
=
\frac{1}{6}\,\log|f|^{2}
-\frac{\chi(X)}{12}\,\varphi-\theta+P.
$$
Since 
$\tau_{\rm BCOV}({\cal X}/\varDelta)
=\mu^{*}\tau_{\rm BCOV}$,
we get
\begin{equation}
\log\tau_{\rm BCOV}({\cal X}/\varDelta)
=
\frac{1}{6}\,\mu^{*}\log|f|^{2}
-\frac{\chi(X)}{12}\,\mu^{*}\varphi
-\mu^{*}\theta
+\mu^{*}P.
\end{equation}
By (10.8), (10.10), we get 
$\log\tau_{\rm BCOV}({\cal X}/\varDelta)
\in L^{1}_{\rm loc}(\varDelta)$.
By (10.9), (10.10), we get the desired equation
of currents.
This complete the proof.
\end{pf}

%%%%%%%%%%%%%%%%%%%%%%%%%%%%%%%%%%%%%%%%%%%%%%%%%%%%%%%%%%%
%
%                 Section 11
%
%%%%%%%%%%%%%%%%%%%%%%%%%%%%%%%%%%%%%%%%%%%%%%%%%%%%%%%%%%%

\section
{\bf The BCOV invariant of Calabi-Yau threefolds with
$h^{1,2}=1$}
\par
In Section 11, we fix the following notation.
Let $\cal X$ be a possibly singular irreducible
projective fourfold and let $S$ be a compact 
Riemann surface.
Let $\pi\colon{\cal X}\to S$ be
a proper, surjective, flat morphism
with discriminant locus 
${\cal D}:=
\{s\in S;\,{\rm Sing}\,X_{s}\not=\emptyset\}$. 
We set
$$
S^{o}:=S\setminus{\cal D},
\qquad 
{\cal X}^{o}:=\pi^{-1}(S^{o}),
$$
$$
{\cal D}^{*}:=\{s\in{\cal D};\,
{\rm Sing}\,X_{s}\hbox{ consists of a unique ODP}\},
$$
and
$$
S^{*}:=S^{o}\cup{\cal D}^{*},
\qquad
{\cal X}^{*}:=\pi^{-1}(S^{*}).
$$
\par
In Section 11, we make the following:
\newline
\newline{\bf Assumption}
(i) 
$X_{s}$ is a Calabi-Yau threefold with
$h^{2}(\Omega^{1}_{X_{s}})=1$ 
for all $s\in S^{*}$;
\newline{(ii)}
${\cal D}^{*}$ is a non-empty finite set,
and ${\cal D}\setminus{\cal D}^{*}$ consists of
a unique point $\infty\in S$;
\newline{(iii)}
${\rm Sing}({\cal X})\cap X_{\infty}=\emptyset$
and $X_{\infty}$ is a divisor of normal crossing.

\begin{lem}
Let $p\in{\cal D}^{*}$. Then $X_{p}$ is smoothable 
in the sense of Definition $2.2$. 
\end{lem}

\begin{pf}
To see this, let $o={\rm Sing}\,X_{p}$, and
let $f\colon\widetilde{X}_{p}\to X_{p}$ be a small
resolution such that $C:=f^{-1}(o)\cong{\Bbb P}^{1}$ and 
$\widetilde{X}_{p}\setminus C\cong X_{p}\setminus\{o\}$. 
Let $[C]\in H_{2}(\widetilde{X}_{p},{\Bbb Z})$ be
the homology class of $C$.
Since $X_{p}$ is smoothable by a flat deformation 
by Assumption (ii), we get $[C]=0$ by
\cite[Th.\,2.5 (2)$\Rightarrow$(3)]{Namikawa02b}.
Hence the map $\gamma'$ in \cite[p.16, l.28]{Namikawa02}
is zero. By the commutative diagram
\cite[p.16 (14)]{Namikawa02}, the natural map
${\rm Ext}^{1}(\Omega^{1}_{X_{p}},{\cal O}_{X_{p}})\to
H^{0}(X,{\cal Ext}(\Omega^{1}_{X_{p}},{\cal O}_{X_{p}}))$ 
is not zero. Let ${\rm Def}(X_{p},o)$ be the Kuranishi
space of the ODP $(X_{p},o)$ and let
$\phi\colon({\rm Def}_{[X_{p}]}(X_{p}),[X_{p}])
\to({\rm Def}(X_{p},o),o)$ be the map of germs 
induced from the Kuranishi family of $X_{p}$. Since 
${\rm Ext}^{1}(\Omega^{1}_{X_{p}},{\cal O}_{X_{p}})=
T_{[X_{p}]}{\rm Def}(X_{p})$ and
$H^{0}(X,{\cal Ext}(\Omega^{1}_{X_{p}},{\cal O}_{X_{p}}))=
T_{p}{\rm Def}(X_{p},o)$ via the Kodaira-Spencer map
and since the natural map
${\rm Ext}^{1}(\Omega^{1}_{X_{p}},{\cal O}_{X_{p}})\to
H^{0}(X,{\cal Ext}(\Omega^{1}_{X_{p}},{\cal O}_{X_{p}}))$
is identified with the differential of $\phi$
at $[X_{p}]$, we get $(d\phi)_{[X_{p}]}\not=0$.
Since $\dim T_{[X_{p}]}{\rm Def}(X_{p})=
\dim T_{o}{\rm Def}(X_{p},o)=1$ by Assumption (i),
$(d\phi)_{[X_{p}]}$ is an isomorphism.
By \cite[Prop.\,5.3]{Namikawa02} and the smoothness
of ${\rm Def}(X_{p},o)$, $\phi$ is an isomorphism
of germs. This implies the smoothness of the total space 
of the Kuranishi family of $X_{p}$.
\end{pf}

\par
The {\it ramification divisor} of the family
$\pi\colon{\cal X}\to S$ is defined as follows.
For $s\in S^{*}$, let 
$\mu_{s}\colon(S,s)\to({\rm Def}(X_{s}),[X_{s}])$
be the map of germs of analytic sets defined by
$$
\mu_{s}(t):=[X_{t}]\in{\rm Def}(X_{s}).
$$
By Lemmas 2.7 and 11.1, $\mu_{p}$ is not a constant
map for $p\in{\cal D}^{*}$. 
Since ${\cal D}^{*}\not=\emptyset$ by Assumption (ii),
$\mu_{s}$ is not constant for all $s\in S^{*}$.
Since $\dim{\rm Def}(X_{s})=1$, we may identify
$({\rm Def}(X_{s}),[X_{s}])$ with $({\Bbb C},0)$.
Let $z$ be the coordinate of $\Bbb C$, 
so that $z\circ\mu_{s}(t)\in{\cal O}_{S,s}$. 
We define the ramification index of
$\pi\colon{\cal X}\to S$ at $s\in S$ by
$$
r_{{\cal X}/S}(s)
:=
{\rm ord}_{t=s}z\circ\mu_{s}(t)\in{\Bbb N}.
$$
Let $\{R_{j}\}_{j\in J}$ be the set of points
of $S$ whose ramification index is $>1$.
The ramification divisor is then defined as
$$
{\cal R}:=\sum_{j\in J}(r_{j}-1)\,R_{j},
\qquad
r_{j}:=r_{{\cal X}/S}(R_{j}).
$$
\par
Let $p\in{\cal D}^{*}$ and
${\rm Sing}(X_{p})=\{o\}$.
By the local description (2.2), 
we have an isomorphism of local rings 
\begin{equation}
{\cal O}_{{\cal X},o}
\cong
{\Bbb C}\{x,y,z,w,t\}/
(x^{2}+y^{2}+z^{2}+w^{2}+t^{r_{{\cal X}/S}(p)}).
\end{equation}
Write ${\cal D}^{*}=\{D_{k}\}_{k\in K}$.
As a divisor of $S$, we define
$$
{\cal D}^{*}
:=
\sum_{k\in K}r_{k}\,D_{k},
\qquad
r_{k}:=r_{{\cal X}/S}(D_{k}).
$$
\par
Since ${\rm Sing}\,{\cal X}\subset
\cup_{s\in{\cal D}^{*}}{\rm Sing}\,X_{s}$,
${\cal X}$ has at most isolated hypersurface
singularities as its singular points by (11.1).
Hence $K_{\cal X}$ and 
$K_{{\cal X}/S}:=K_{\cal X}\otimes\pi^{*}K_{S}^{-1}$
are invertible sheaves on ${\cal X}$.

\begin{lem}
The sheaf $\pi_{*}K_{{\cal X}/S}$ is 
an invertible sheaf on $S$.
\end{lem}

\begin{pf}
Since $\pi^{-1}(S\setminus{\cal D}^{*})$ is smooth,
$\pi_{*}K_{{\cal X}/S}$ is an invertible sheaf 
on $S\setminus{\cal D}^{*}$ by Assumption (i) and 
\cite[p.391, Th.\,V]{Takegoshi95}. 
Let $s\in{\cal D}^{*}$.
Since the conormal bundle of $(X_{s})_{\rm reg}$
in ${\cal X}_{\rm reg}$ is trivial,
we have $K_{{\cal X}/S}|_{(X_{s})_{\rm reg}}\cong
K_{(X_{s})_{\rm reg}}$. Since $K_{{\cal X}/S}|_{X_{s}}$
and $K_{X_{s}}$ are invertible sheaves on $X_{s}$,
we get $K_{{\cal X}/S}|_{X_{s}}\cong K_{X_{s}}$
by the normality of $X_{s}$. Since $X_{s}$ is
Calabi-Yau, we have
$h^{0}(K_{{\cal X}/S}|_{X_{s}})=h^{0}(K_{X_{s}})=1$.
By \cite[Th.\,4.12 (ii)]{Banica76}, 
$\pi_{*}K_{{\cal X}/S}$ is an invertible sheaf
near $s\in{\cal D}^{*}$.
This proves the lemma.
\end{pf}

Let $\chi$ be the topological Euler number
of a general fiber $X_{s}$, $s\in S^{o}$.
Let $\|\cdot\|$ be the Hermitian metric on
$(\pi_{*}K_{{\cal X}/S})^{\otimes(48+\chi)}
\otimes(TS)^{\otimes12}|_{S^{o}}$ induced from
the $L^{2}$-metric on $\pi_{*}K_{{\cal X}/S}$
and from the Weil-Petersson metric 
$g_{{\rm WP},{\cal X}/S}$ on $S^{o}$.
The following is the main result of this paper.

\begin{Mthm}
Let $\varXi$ be a meromorphic section of
$\pi_{*}K_{{\cal X}/S}$ on $S$ with
$$
{\rm div}(\varXi)
=
\sum_{i\in I}m_{i}\,P_{i}+m_{\infty}P_{\infty},
\qquad
P_{i}\not=P_{\infty}\,(i\in I),
$$ 
and let $V$ be a meromorphic vector field on $S$.
Then the following hold:
\newline{$(1)$ }
There exists a locally integrable function
$F_{\varXi,V}$ on $S$ with
$$
\begin{aligned}
dd^{c}F_{\varXi,V}
&
=
\left\{
(24+\frac{\chi}{2})\,\deg\pi_{*}K_{{\cal X}/S}
+6\,\chi(S)+6\,\deg{\cal R}
-\deg{\cal D}^{*}
\right\}\,\delta_{\infty}
\\
&
\quad
+\delta_{{\cal D}^{*}}
-(24+\frac{\chi}{2})\,\delta_{{\rm div}(\varXi)}
-6\,\delta_{{\rm div}(V)}-6\,\delta_{\cal R}
\end{aligned}
$$ 
such that
$$
\tau_{\rm BCOV}({\cal X}/S)
=
\left\|
e^{F_{\varXi,V}}\,\varXi^{48+\chi}
\otimes V^{12}
\right\|^{\frac{1}{6}}.
$$
\newline{$(2)$ }
When $S={\Bbb P}^{1}$,
let $\psi$ be the inhomogeneous coordinate
of ${\Bbb P}^{1}$ with $\psi(\infty)=\infty$.
Identify the points $P_{i},R_{j},D_{k}$ 
with their coordinates 
$\psi(P_{i}),\psi(R_{j}),\psi(D_{k})$,
respectively.
Then there exists a constant $C\not=0$ 
such that
$$
\tau_{\rm BCOV}(X_{\psi})
=
C\,\left\|
\prod_{i\in I,j\in J,k\in K}
\frac{(\psi-D_{k})^{2r_{k}}}
{(\psi-P_{i})^{(48+\chi)m_{i}}
(\psi-R_{j})^{12(r_{j}-1)}}
\,\varXi_{\psi}^{48+\chi}
\otimes
\left(\frac{\partial}{\partial\psi}\right)^{12}
\right\|^{\frac{1}{6}}.
$$
\end{Mthm}

In the rest of this section, we shall prove
Theorem 11.3.
For $p\in{\cal D}$,
let $(U_{p},t)$ be a coordinate neighborhood of $S$
centered at $p$ with
$U_{p}\cap{\cal D}=\{p\}$ 
and $U_{p}\setminus\{p\}\cong\varDelta^{*}$.
\par
By Proposition 7.3, the positive $(1,1)$-forms
$\omega_{{\rm WP},{\cal X}/S}$ and 
$\omega_{{\rm H},{\cal X}/S}$
on $S^{o}$ extend trivially to
closed positive $(1,1)$-currents on $S$.

\begin{defn}
Let
$\Omega_{{\rm WP},{\cal X}/S}$ and 
$\Omega_{{\rm H},{\cal X}/S}$
be the trivial extensions of
$\omega_{{\rm WP},{\cal X}/S}$ and 
$\omega_{{\rm H},{\cal X}/S}$
from $S^{o}$ to $S$, respectively.
\end{defn}

\begin{prop}
$(1)$
There exists $a(p)\in{\Bbb R}$ such that the 
following equation of currents on $U_{p}$ holds:
$$
dd^{c}\log\Omega_{{\rm WP},{\cal X}/S}|_{U_{p}}
\left(
\frac{\partial}{\partial t},
\frac{\partial}{\partial\bar{t}}
\right)
=
a(p)\,\delta_{p}
-\Omega_{{\rm H},{\cal X}/S}+
4\,\Omega_{{\rm WP},{\cal X}/S}.
$$
\newline{$(2)$}
For $D_{k}\in{\cal D}^{*}$, one has
$a(D_{k})=r_{k}-1$.
\end{prop}

\begin{pf}
We get (1) by Proposition 7.6 (2).
Let $p=D_{k}$. Under the identification of 
the Kuranishi space $({\rm Def}(X_{p}),[X_{p}])$
with $({\Bbb C},0)$, we may assume 
by the definition of the ramification index
$r_{{\cal X}/S}$ that 
$\pi|_{U_{p}}\colon{\cal X}|_{U_{p}}\to U_{p}$ is
induced from the Kuranishi family of $X_{p}$ by
the map $\mu(t)=t^{r_{k}}$. 
Let $\omega_{\rm WP}$ be the Weil-Petersson form
on ${\rm Def}(X_{p})$. 
Since 
$\Omega_{{\rm WP},{\cal X}/S}|_{U_{p}\setminus\{p\}}
=\mu^{*}\omega_{\rm WP}$,
we deduce from Proposition 7.6 (1), (3) that
as $t\to0$,
\begin{equation}
\begin{aligned}
\log\Omega_{{\rm WP},{\cal X}/S}|_{U_{p}}
\left(
\frac{\partial}{\partial t},
\frac{\partial}{\partial\bar{t}}
\right)
&
=
\log\omega_{\rm WP}
\left(
\mu_{*}\frac{\partial}{\partial t},
\mu_{*}\frac{\partial}{\partial\bar{t}}
\right)
\\
&
=
(r_{k}-1)\,\log|t|^{2}+O(\log(-\log|t|)).
\end{aligned}
\end{equation}
By (11.2), we get $a(p)=r_{k}-1$.
This completes the proof.
\end{pf}

\begin{prop}
There exists $b(\infty)\in{\Bbb R}$ such that
the following equation of currents on $S$ holds:
\begin{equation}
dd^{c}\log \|\varXi\|_{L^{2}}^{2}
=
b(\infty)\,\delta_{\infty}
+\delta_{{\rm div}(\varXi)}
-\Omega_{{\rm WP},{\cal X}/S}.
\end{equation}
\end{prop}

\begin{pf}
Let $s\in S$ be an arbitrary point. 
It suffices to prove Eq.\,(11.3) on a neighborhood
of $s$. By Proposition 7.8 (2), Eq.\,(11.3) holds
on a neighborhood of $\infty$. 
\par
Assume that
$s\in S^{*}$. Let 
${\frak p}\colon({\frak X},X_{s})\to
({\rm Def}(X_{s}),[X_{s}])$
be the Kuranishi family of $X_{s}$.
Since $\pi\colon({\cal X},X_{s})\to(S,s)$ is
induced from the Kuranishi family by the map
$\mu_{s}\colon(S,s)\to({\rm Def}(X_{s}),[X_{s}])$,
there exists a morphism of deformation germs 
$f_{\mu_{s}}\colon({\cal X},X_{s})\to({\frak X},X_{s})$
satisfying the commutative diagram:
$$
\begin{CD}
({\cal X},X_{s})@> f_{\mu_{s}} >> ({\frak X},X_{s})
\\
@V \pi VV @V {\frak p} VV
\\
(S,s)
@> \mu_{s} >> ({\rm Def}(X_{s}),[X_{s}]).
\end{CD}
$$
\par
Let $U_{s}\cong\varDelta$ be a neighborhood of $s$ 
in $S$ such that $\mu_{s}$ (resp. $f_{\mu_{s}}$) is
defined on $U_{s}$ (resp. $\pi^{-1}(U_{s})$)
and such that $\mu_{s}$ has no critical points
on $U_{s}^{o}:=U_{s}\setminus\{s\}$. Since
\begin{equation}
f_{\mu_{s}}^{*}K_{{\frak X}/{\rm Def}(X_{s})}
=
K_{{\cal X}/S}
\end{equation}
on $\pi^{-1}(U_{s})\setminus{\rm Sing}\,X_{s}$, 
the normality of ${\cal X}$ implies that
(11.4) holds on $\pi^{-1}(U_{s})$.
\par
By Lemma 6.1, $K_{{\frak X}/{\rm Def}(X_{s})}$
is trivial.
Let $\eta_{{\frak X}/{\rm Def}(X_{s})}$ be
a nowhere vanishing holomorphic section of
$K_{{\frak X}/{\rm Def}(X_{s})}$ defined on
${\rm Def}(X_{s})$.
We regard $\eta_{{\frak X}/{\rm Def}(X_{s})}$
as a family of holomorphic $3$-forms
$\{\eta_{{\frak X}/{\rm Def}(X_{s})}|_{{\frak X}_{b}}
\}_{b\in{\rm Def}(X_{s})}$.
Since $X_{s}$ has at most one ODP as its singular set,
$\log\|\eta_{{\frak X}/{\rm Def}(X_{s})}\|_{L^{2}}
\in C^{0}({\rm Def}(X_{s}))$ by Proposition 7.8 (3). 
\par
Since 
$f_{\mu_{s}}^{*}\eta_{{\frak X}/{\rm Def}(X_{s})}
\in H^{0}(\pi^{-1}(U_{s}),K_{{\cal X}/S})
=H^{0}(U_{s},\pi_{*}K_{{\cal X}/S})$
is nowhere vanishing, 
$f_{\mu_{s}}^{*}\eta_{{\frak X}/{\rm Def}(X_{s})}$
generates $\pi_{*}K_{{\cal X}/S}$ on $U_{s}$ as
an ${\cal O}_{U_{s}}$-module.
Since
$$
\|f_{\mu_{s}}^{*}\eta_{{\frak X}/{\rm Def}(X_{s})}
\|_{L^{2}}(t)
=
\|\eta_{{\frak X}/{\rm Def}(X_{s})}\|_{L^{2}}(\mu_{s}(t)),
\qquad
t\in U_{s}^{o}
$$
by (11.4) and since
$\log\|\eta_{{\frak X}/{\rm Def}(X_{s})}\|_{L^{2}}
\in C^{0}({\rm Def}(X_{s}))$, 
$\log\|f_{\mu_{s}}^{*}
\eta_{{\frak X}/{\rm Def}(X_{s})}\|_{L^{2}}$
is a continuous function on $U_{s}$.
Since
$-dd^{c}\log\|f_{\mu_{s}}^{*}
\eta_{{\frak X}/{\rm Def}(X_{s})}\|_{L^{2}}=
\Omega_{{\rm WP},{\cal X}/S}$ on $U_{s}^{o}$, 
we get the following equation of currents 
on $U_{s}$ by Lemma 7.5 (1), (2):
\begin{equation}
-dd^{c}\log\|f_{\mu_{s}}^{*}
\eta_{{\frak X}/{\rm Def}(X_{s})}\|_{L^{2}}
=\Omega_{{\rm WP},{\cal X}/S}.
\end{equation}
Since 
$f_{\mu_{s}}^{*}\eta_{{\frak X}/{\rm Def}(X_{s})}
\in H^{0}(U_{s},\pi_{*}K_{{\cal X}/S})$
is nowhere vanishing, there exists 
$h(t)\in{\cal O}(U_{s})$ such that
$\varXi=h\cdot
f_{\mu_{s}}^{*}\eta_{{\frak X}/{\rm Def}(X_{s})}$
on $U_{s}$. By (11.5), we get
\begin{equation}
-dd^{c}\log\|\varXi\|_{L^{2}}^{2}
=
\Omega_{{\rm WP},{\cal X}/S}-\delta_{{\rm div}(h)}
\end{equation}
as currents on $U_{s}$.
Eq.\,(11.3) on $U_{s}$ follows from (11.6).
\end{pf}

\begin{thm}
There exists $c(\infty)\in{\Bbb Q}$ 
such that
the following equation of currents on $S$
holds:
\begin{equation}
dd^{c}\log\tau_{\rm BCOV}({\cal X}/S)
=
-\frac{\chi}{12}\,\Omega_{{\rm WP},{\cal X}/S}
-\Omega_{{\rm H},{\cal X}/S}
+\frac{1}{6}\,\delta_{{\cal D}^{*}}
+c(\infty)\,\delta_{\infty}.
\end{equation}
\end{thm}

\begin{pf}
The result follows from Theorems 10.1 and 10.3.
\end{pf}

{\bf Proof of Theorem 11.3}
{\bf (1)}
By Proposition 11.5 and (4.1),
we get the following equation of currents on $S$:
\begin{equation}
dd^{c}\log\|V\|^{2}
=
a(\infty)\,\delta_{\infty}
+\delta_{\cal R}+\delta_{{\rm div}(V)}
-\Omega_{{\rm H},{\cal X}/S}+
4\,\Omega_{{\rm WP},{\cal X}/S}.
\end{equation}
By (11.3), (11.7), (11.8), we get
\begin{equation}
\begin{aligned}
dd^{c}\log
\|V^{12}\otimes\varXi^{48+\chi}\|^{2}
&
=
12(a(\infty)\,\delta_{\infty}
+\delta_{\cal R}+\delta_{{\rm div}(V)})
-12\,\Omega_{{\rm H},{\cal X}/S}
+48\,\Omega_{{\rm WP},{\cal X}/S}
\\
&
\quad
+(48+\chi)\,(b(\infty)\,\delta_{\infty}
+\delta_{{\rm div}(\varXi)})
-(48+\chi)\,\Omega_{{\rm WP},{\cal X}/S}
\\
&
=
12\,dd^{c}\log\tau_{\rm BCOV}({\cal X}/S)
\\
&
\quad
+\{12\,a(\infty)+(48+\chi)\,b(\infty)
-12\,c(\infty)\}\,\delta_{\infty}
\\
&
\quad
-2\,\delta_{{\cal D}^{*}}+12\,\delta_{\cal R}
+12\,\delta_{{\rm div}(V)}
+(48+\chi)\,\delta_{{\rm div}(\varXi)}.
\end{aligned}
\end{equation}
Integrating the both hand sides of (11.9)
over $S$, we get
\begin{equation}
\{12\,a(\infty)+(48+\chi)\,b(\infty)
-12\,c(\infty)\}
-2\,\deg{\cal D}^{*}+12\,\deg{\cal R}
+12\,\chi(S)
+(48+\chi)\,\deg\varXi
=0.
\end{equation}
By (11.9) and (11.10),
$$
F_{\varXi,V}
:=
\log\tau_{\rm BCOV}({\cal X}/S)^{6}
-\log\|V^{12}\otimes\varXi^{48+\chi}\|
$$
is a harmonic function on 
$S\setminus({\cal D}\cup{\cal R})$ satisfying 
Theorem 11.3 (1). This proves (1).
\newline{\bf (2)}
We set $V(\psi):=\partial/\partial\psi
\in H^{0}({\Bbb P}^{1},T{\Bbb P}^{1})$.
Then ${\rm div}(V)=2\,\infty$, so that
$F_{\varXi,V}$ satisfies the following equation 
of currents on ${\Bbb P}^{1}$ by (11.9), (11.10):
\begin{equation}
\begin{aligned}
dd^{c}F_{\varXi,V}
&
=
\left\{
(24+\frac{\chi}{2})\,\deg\pi_{*}K_{{\cal X}/S}
+6\,\deg{\cal R}-\deg{\cal D}^{*}
\right\}\,\delta_{\infty}
\\
&
\quad
+\delta_{{\cal D}^{*}}
-(24+\frac{\chi}{2})\,\delta_{{\rm div}(\varXi)}
-6\,\delta_{\cal R}.
\end{aligned}
\end{equation}
Up to a constant, the solution of Eq.\,(11.11)
is given by the following formula:
\begin{equation}
F_{\varXi,V}(\psi)
=
\log
\left|
\prod_{i\in I,j\in J,k\in K}
\frac{(\psi-D_{k})^{2r_{k}}}
{(\psi-P_{i})^{(48+\chi)m_{i}}
(\psi-R_{j})^{12(r_{j}-1)}}
\right|.
\end{equation}
The second assertion of Theorem 11.3 follows
from (11.12).
This completes the proof of Theorem 11.3.
\qed

%%%%%%%%%%%%%%%%%%%%%%%%%%%%%%%%%%%%%%%%%%%%%%%%%%%%%%%%%%%
%
%                 Section 12
%
%%%%%%%%%%%%%%%%%%%%%%%%%%%%%%%%%%%%%%%%%%%%%%%%%%%%%%%%%%%

\section
{\bf The BCOV invariant of quintic mirror threefolds}
\par

\subsection
{}
{\bf Quintic mirror threefolds}
\par
Let $p\colon{\cal X}\to{\Bbb P}^{1}$ be the pencil of 
quintic threefolds in ${\Bbb P}^{4}$ defined by
$$
{\cal X}:=
\{([z],\psi)\in{\Bbb P}^{4}\times{\Bbb P}^{1};\,
F_{\psi}(z)=0\},
\qquad
p={\rm pr}_{2},
$$
$$
F_{\psi}(z):=
z_{0}^{5}+z_{1}^{5}+z_{2}^{5}+z_{3}^{5}+z_{4}^{5}
-5\psi\,z_{0}z_{1}z_{2}z_{3}z_{4}.
$$
The parameter $\psi$ is regarded as the inhomogeneous
coordinate of ${\Bbb P}^{1}$.
Identify ${\Bbb Z}_{5}$ with the set of fifth roots
of unity:
${\Bbb Z}_{5}=\{\zeta\in{\Bbb C};\,\zeta^{5}=1\}$. 
We define
$$
G:=
\frac{
\{(a_{0},a_{1},a_{2},a_{3},a_{4})\in({\Bbb Z}_{5})^{5};\,
a_{0}a_{1}a_{2}a_{3}a_{4}=1\}}
{{\Bbb Z}_{5}(1,1,1,1,1)}
\cong
{\Bbb Z}_{5}^{3}.
$$
The group $G\times{\Bbb Z}_{5}$ acts
on ${\cal X}$ and ${\Bbb P}^{1}$ by
$$
(a,b)\cdot([z],\psi)
:=
([b^{-1}a_{0}z_{0}:a_{1}z_{1}:a_{2}z_{2}:a_{3}z_{3}:
a_{4}z_{4}],b\psi),
\qquad
(a,b)\cdot\psi:=b\,\psi.
$$
Then the projection $p\colon{\cal X}\to{\Bbb P}^{1}$
is $G\times{\Bbb Z}_{5}$-equivariant.
Since $G$ preserves the fibers of $p$, 
we have the induced family
$$
p\colon{\cal X}/G\to{\Bbb P}^{1}
$$
equipped with the induced ${\Bbb Z}_{5}$-action.
We set
$$
{\cal D}^{*}:=
\left\{
\exp\frac{2\pi\sqrt{-1}\,m}{5}\in{\Bbb P}^{1};\,
0\leq m\leq 4
\right\}
\subset{\Bbb P}^{1},
\qquad
{\cal D}:={\cal D}^{*}\cup\{\infty\}
\subset{\Bbb P}^{1}.
$$
Then ${\cal D}$ is the discriminant locus of
the family $p\colon{\cal X}\to{\Bbb P}^{1}$ 
by \cite[p.27]{CdlOGP93}.

\begin{prop}
There exists a resolution 
$f\colon{\cal W}\to{\cal X}/G$
satisfying the following conditions:
\newline{$(1)$ }
Set $f_{\psi}:=f|_{W_{\psi}}$.
Then $f_{\psi}\colon W_{\psi}\to X_{\psi}/G$ 
is a crepant resolution 
for $\psi\in{\Bbb P}^{1}\setminus{\cal D}$.
In particular, $W_{\psi}$ is a smooth Calabi-Yau
threefold for $\psi\in{\Bbb P}^{1}\setminus{\cal D}$;
\newline{$(2)$ }
${\rm Sing}\,W_{\psi}$ consists of a unique ODP
if $\psi^{5}=1$;
\newline{$(3)$ }
$W_{\infty}$ is a divisor of normal crossing.
\end{prop}

\begin{pf}
See \cite[Appendix B]{Morrison93}, \cite{Batyrev94},
\cite[Sects.\,2.2 and 2.4]{CoxKatz99} for (1)
and \cite[p.27]{CdlOGP93} for (2).
The last assertion follows from Hironaka's theorem.
\end{pf}

Notice that the choice of a resolution
$f\colon{\cal W}\to{\cal X}/G$ as above is 
{\it not} unique.

\begin{defn}
Set $\pi:=p\circ f$.
Any family $\pi\colon{\cal W}\to{\Bbb P}^{1}$ 
satisfying the conditions (1), (2), (3) as above 
is called a {\it family of quintic mirror threefolds}.
The induced family 
$\pi\colon{\cal W}/{\Bbb Z}_{5}\to
{\Bbb P}^{1}/{\Bbb Z}_{5}$
is also called a family of quintic mirror threefolds.
\end{defn}

\begin{lem}
If $\psi\in{\Bbb P}^{1}\setminus{\cal D}$,
then
$$
h^{1,2}(W_{\psi})=1,
\qquad
h^{1,1}(W_{\psi})=101,
\qquad
\chi(W_{\psi})=200.
$$
\end{lem}

\begin{pf}
Since $h^{1,1}(X_{\psi})=1$, $h^{1,2}(X_{\psi})=101$,
and $\chi(X_{\psi})=-200$, the result
follows from \cite{Batyrev94}, \cite[Th.\,4.1.5]{CoxKatz99},
\cite[Th.\,4.30]{Voisin99}.
\end{pf}

We refer to 
\cite{CdlOGP93}, \cite{CoxKatz99}, 
\cite{Morrison93}, \cite{Voisin99} for more details 
about quintic mirror threefolds.

\subsection
{}{\bf The mirror map}
\par

\begin{defn}
The {\it mirror map} is the holomorphic map 
from a neighborhood of $\infty\in{\Bbb P}^{1}$
to a neighborhood of $0\in\varDelta$
defined by the following formula:
$$
q:=(5\psi)^{-5}\exp\left(\frac{5}{y_{0}(\psi)}
\sum_{n=1}^{\infty}\frac{(5n)!}{(n!)^{5}}
\left\{\sum_{j=n+1}^{5n}\frac{1}{j}\right\}
\frac{1}{(5\psi)^{5n}}\right),
\qquad
|\psi|\gg1,
$$
where
$$
y_{0}(\psi)
:=
\sum_{n=1}^{\infty}\frac{(5n)!}{(n!)^{5}(5\psi)^{5n}},
\qquad
|\psi|>1.
$$
The inverse of the mirror map is denoted by
$\psi(q)$.
\end{defn}

For $\psi\in{\Bbb P}^{1}\setminus{\cal D}$,
we define a holomorphic $3$-form on $X_{\psi}$ by
$$
\varOmega_{\psi}
:=
\left(\frac{2\pi\sqrt{-1}}{5}\right)^{-3}
5\psi\frac{z_{4}\,dz_{0}\wedge dz_{1}\wedge dz_{2}}
{\partial F_{\psi}(z)/\partial z_{3}}.
$$
Since $\varOmega_{\psi}$ is $G$-invariant,
$\varOmega_{\psi}$ induces a holomorphic $3$-form
on $X_{\psi}/G$ in the sense of orbifolds.
We identify $\varOmega_{\psi}$ with the corresponding
holomorphic $3$-form on $X_{\psi}/G$.
Then $f_{\psi}^{*}\Omega_{\psi}$ is a holomorphic
$3$-form on $W_{\psi}$. Set
$$
\varXi_{\psi}:=f_{\psi}^{*}\varOmega_{\psi}
\in H^{0}(W_{\psi},K_{W_{\psi}}).
$$
By Lemma 12.3, we know 
${\rm rk}\,H_{3}(W_{\psi},{\Bbb Z})=4$.
There exists a symplectic basis 
$\{A^{1},A^{2},B_{1},B_{2}\}$ of 
$H_{3}(W_{\psi},{\Bbb Q})$, $\psi\not\in{\cal D}$,
such that $A^{a}\cap B_{b}=\delta_{an}$,
$A^{a}\cap A^{b}=B_{a}\cap B_{b}=0$.
By \cite{CdlOGP93}, \cite[p.245 l.13]{Morrison93},
the mirror map $q(\psi)$ is expressed as follows:
$$
q
=
\exp
\left(2\pi\sqrt{-1}\,
\frac{\int_{2B_{1}-A^{1}}\varXi_{\psi}}
{\int_{A^{2}}\varXi_{\psi}}
\right),
\qquad
y_{0}(\psi)
=
\int_{A^{2}}\varXi_{\psi}.
$$

We refer to 
\cite{CdlOGP93}, 
\cite[Sect.\,2.3, Sect.\,6.3]{CoxKatz99}, 
\cite{Morrison93}, \cite[Chap.\,3]{Voisin99}
for more details about the mirror map.

\subsection
{}{\bf Conjectures of Bershadsky-Cecotti-Ooguri-Vafa}
\par

\begin{defn}
Under the identification of the local parameters
$\psi^{5}$ and $q$ via the mirror map,
define a multi-valued analytic function 
near $\infty\in{\Bbb P}^{1}$ as
$$
F_{1,B}^{\rm top}(\psi)
:=
\left(\frac{\psi}{y_{0}(\psi)}\right)^{\frac{62}{3}}
(\psi^{5}-1)^{-\frac{1}{6}}\,q\,\frac{d\psi}{dq}
$$
and a power series in $q$ as
$$
F_{1,A}^{\rm top}(q)
:=
F_{1,B}^{\rm top}(\psi(q)).
$$
\end{defn}

Set 
$$
\eta(q):=\prod_{n=1}^{\infty}(1-q^{n}).
$$
In \cite[Eq.(16), (23), (24)]{BCOV93} and 
\cite[l.34]{BCOV94},
Bershadsky-Cecotti-Ooguri-Vafa conjectured 
the following:

\begin{conj}
{\bf (A)}
Let $N_{g}(d)$ be the genus-$g$ Gromov-Witten 
invariant of degree $d$ of a general quintic 
threefold in ${\Bbb P}^{4}$
(cf. \cite{LiZinger04}).
Then the following identity holds:
$$
q\,\frac{d}{dq}\log F_{1,A}^{\rm top}(q)
=
\frac{50}{12}
-
\sum_{n,d=1}^{\infty}N_{1}(d)\,
\frac{2nd\,q^{nd}}{1-q^{nd}}
-
\sum_{d=1}^{\infty}N_{0}(d)\,
\frac{2d\,q^{d}}{12(1-q^{d})},
$$
or equivalently
$$
F_{1,A}^{\rm top}(q)
=
{\rm Const.}\,
\left\{
q^{25/12}\,\prod_{d=1}^{\infty}
\eta(q^{d})^{N_{1}(d)}
(1-q^{d})^{N_{0}(d)/12}\right\}^{2}.
$$
\newline{\bf (B)}
Let $\|\cdot\|$ be the Hermitian metric 
on the line bundle
$(\pi_{*}K_{{\cal W}/{\Bbb P}^{1}})^{\otimes62}
\otimes(T{\Bbb P}^{1})^{\otimes3}
|_{{\Bbb P}^{1}\setminus{\cal D}}$
induced from the $L^{2}$-metric on 
$\pi_{*}K_{{\cal W}/{\Bbb P}^{1}}$ and from
the Weil-Petersson metric on $T{\Bbb P}^{1}$.
Then the following identity holds:
$$
\begin{aligned}
\tau_{\rm BCOV}(W_{\psi})
&
=
{\rm Const.}\,
\left\|
\psi^{-62}
(\psi^{5}-1)^{\frac{1}{2}}\,
(\varXi_{\psi})^{62}
\otimes\left(\frac{d}{d\psi}\right)^{3}
\right\|^{\frac{2}{3}}
\\
&
=
{\rm Const.}\,
\left\|\frac{1}{F_{1,B}^{\rm top}(\psi)^{3}}\,
\left(
\frac{\varXi_{\psi}}{y_{0}(\psi)}
\right)^{62}
\otimes\left(q\frac{d}{dq}\right)^{3}
\right\|^{\frac{2}{3}}.
\end{aligned}
$$
\end{conj}

\begin{rem}
Under Conjecture 12.6, 
the Gromov-Witten invariants 
$\{N_{g}(d)\}_{g\leq1,d\in{\Bbb N}}$
of a general quintic threefold in ${\Bbb P}^{4}$ 
and the BCOV invariant of the mirror
quintic threefolds satisfy the following relation:
$$
\tau_{\rm BCOV}(W_{\psi})
=
{\rm Const.}\,
\left\|
\left\{q^{\frac{25}{12}}\,
\prod_{d=1}^{\infty}\eta(q^{d})^{N_{1}(d)}
(1-q^{d})^{\frac{N_{0}(d)}{12}}\right\}^{6}
\,
\left(\frac{\varXi_{\psi}}{y_{0}(\psi)}\right)^{62}
\otimes\left(q\frac{d}{dq}\right)^{3}
\right\|^{\frac{2}{3}}.
$$
\end{rem}

In the rest of this section,
we prove Conjecture 12.6 (B)
as an application of Theorem 11.3.

\subsection
{}
{\bf Proof of Conjecture 12.6 (B)}
\par
Let $\pi\colon{\cal W}\to{\Bbb P}^{1}$ be
a family of quintic mirror threefolds.
Let $K(\psi)$ be the K\"ahler potential of 
the Weil-Petersson form $\Omega_{\rm WP}$ 
defined as
$$
K(\psi)
:=
-\log\left(
\sqrt{-1}
\int_{W_{\psi}}\varXi_{\psi}\wedge
\overline{\varXi}_{\psi}
\right).
$$
Define a function $G(\psi)$ by
$G(\psi)=g_{\rm WP}(\frac{\partial}{\partial\psi},
\frac{\partial}{\partial\bar{\psi}})$, so that
$$
\Omega_{\rm WP}(\psi)
=
\sqrt{-1}\,G(\psi)\,d\psi\wedge d\bar{\psi}
=
\frac{\sqrt{-1}}{2\pi}
\frac{\partial^{2}K(\psi)}
{\partial\psi\partial\bar{\psi}}\,
d\psi\wedge d\bar{\psi}.
$$

\begin{prop}
The following estimates hold
\begin{equation}
K(\psi)
=
\left\{
\begin{array}{ll}
\log|\psi|^{2}+O(1)
&(\psi\to0)
\\
O(1)
&(\psi^{5}\to1)
\\
O(\log\log|\psi|)
&(\psi\to\infty),
\end{array}
\right.
\end{equation}
\begin{equation}
\log G(\psi)
=
\left\{
\begin{array}{ll}
O(1)
&(\psi\to0)
\\
O(\log(-\log|\psi^{5}-1|))
&(\psi^{5}\to1)
\\
-\log|\psi|^{2}+O(\log\log|\psi|)
&(\psi\to\infty).
\end{array}
\right.
\end{equation}
In particular, 
${\cal R}\cap{\cal D}^{*}=\emptyset$
for any family of quintic mirror threefolds.
\end{prop}

\begin{pf}
See \cite[p.50 Table 2]{CdlOGP93}.
\end{pf}

\begin{prop}
The family of quintic mirror threefolds has
trivial ramification divisor, i.e.,
${\cal R}=0$ for the family 
$\pi\colon{\cal W}\to{\Bbb P}^{1}$.
\end{prop}

\begin{pf}
By (11.2) and Proposition 12.8, 
if suffices to prove that $G(\psi)>0$ on
${\Bbb P}^{1}\setminus{\cal D}$. Since
$$
K(\psi)
=
-\log\left(
\frac{\sqrt{-1}}{|G|}
\int_{X_{\psi}}\varOmega_{\psi}\wedge
\overline{\varOmega}_{\psi}
\right),
$$
$\Omega_{\rm WP}(\psi)$ coincides with
the Weil-Petersson form for $X_{\psi}$ by (4.1).
Thus $G(\psi)>0$ if and only if 
the Kodaira-Spencer map 
$\mu_{\psi}\colon T_{\psi}{\Bbb P}^{1}\to
H^{1}(X_{\psi},\Theta_{X_{\psi}})$
for $p\colon{\cal X}\to{\Bbb P}^{1}$ 
is non-degenerate at 
$\psi\in{\Bbb P}^{1}\setminus{\cal D}$.
By \cite[p.53 l.18-l.27]{Voisin99},
$\mu_{\psi}$ is non-degenerate for all
$\psi\in{\Bbb P}^{1}\setminus{\cal D}$.
This proves the proposition.
\end{pf}

\begin{thm}
Conjecture 12.6 (B) holds.
\end{thm}

\begin{pf}
For a point $z=(1:z)\in{\Bbb P}^{1}$, 
let $[z]=[(1:z)]$ denote the corresponding divisor.
By Proposition 12.1, we get
\begin{equation}
{\rm div}({\cal D}^{*})
=
\sum_{\zeta^{5}=1}[\zeta],
\end{equation}
which is a reduced divisor.
By (12.1), we have
\begin{equation}
{\rm div}(\varXi)=[0].
\end{equation}
Substituting (12.3), (12.4) and ${\cal R}=0$
into the formula for $\tau_{\rm BCOV}$ in Theorem 11.3 (2) 
and using $\chi(W_{\psi})=200$, we get
\begin{equation}
\begin{aligned}
\tau_{\rm BCOV}(W_{\psi})
&
=
{\rm Const.}\,\left\|
\frac{\prod_{\zeta^{5}=1}(\psi-\zeta)^{2}}
{\psi^{48+\chi}}\,
\varXi_{\psi}^{48+\chi}\otimes
\left(\frac{\partial}{\partial\psi}\right)^{12}
\right\|^{1/6}
\\
&
=
{\rm Const.}\,\left\|
\frac{(\psi^{5}-1)^{2}}
{\psi^{248}}\,
\varXi_{\psi}^{248}\otimes
\left(\frac{\partial}{\partial\psi}\right)^{12}
\right\|^{1/6}
\\
&
=
{\rm Const.}\,\left\|
\psi^{-62}(\psi^{5}-1)^{1/2}\,
\varXi_{\psi}^{62}\otimes
\left(\frac{\partial}{\partial\psi}\right)^{3}
\right\|^{2/3}.
\end{aligned}
\end{equation}
This proves Conjecture 12.6 (B).
\end{pf}

\begin{rem}
It seems that the families of Calabi-Yau
threefolds over ${\Bbb P}^{1}$ studied in
\cite[Eqs.\, (2.1), (2.2)]{KlemmTheisen93}
satisfy Assumption (i), (ii), (iii) of Sect.\,11.
(See \cite[p.157, last five lines]{KlemmTheisen93}.)
By the explicit formula for the Yukawa coupling
\cite[Eq.\,(4.6)]{KlemmTheisen93}, we get 
${\cal R}\cap({\Bbb P}^{1}\setminus{\cal D})=\emptyset$
for these examples. 
If ${\cal R}\cap{\cal D}^{*}=\emptyset$,
the conjectured formulas for the BCOV invariants 
of these families \cite[p.294]{BCOV93}
follow from Theorem 11.3 (2).
\end{rem}

%%%%%%%%%%%%%%%%%%%%%%%%%%%%%%%%%%%%%%%%%%%%%%%%%%%%%%%%%%%
%
%                 Section 13
%
%%%%%%%%%%%%%%%%%%%%%%%%%%%%%%%%%%%%%%%%%%%%%%%%%%%%%%%%%%%
\section
{\bf The BCOV invariant of FHSV threefolds}
\par

\subsection
{}
{\bf The threefolds of Ferrara-Harvey-Strominger-Vafa}
\par
A compact connected complex surface $S$ is
an {\it Enriques surface} if it satisfies
$H^{1}(S,{\cal O}_{S})=0$, 
$K_{S}\not\cong{\cal O}_{S}$, and 
$K_{S}^{2}\cong{\cal O}_{S}$. An Enriques surface $S$
is an algebraic surface with 
$\pi_{1}(S)\cong{\Bbb Z}_{2}$ whose universal covering
$\widetilde{S}$ is a {\it $K3$ surface}. 
For an Enriques surface $S$,
let $\iota_{S}\colon\widetilde{S}\to\widetilde{S}$
be the non-trivial covering transformation that
generates $\pi_{1}(S)$. 
Then $(\widetilde{S},\iota_{S})$ 
is a $2$-elementray $K3$ surface.
(See \cite[Sect.\,8.1]{Yoshikawa04}.)
\par
Let ${\Bbb H}\subset{\Bbb C}$ be the complex
upper-half plane.
For $\tau\in{\Bbb H}$, let $E_{\tau}$ denote
the elliptic curve ${\Bbb C}/{\Bbb Z}+\tau{\Bbb Z}$.
For an elliptic curve $T=E_{\tau}$, 
let $-1_{T}$ be the involution on $T$ defined as
$-1_{T}(z)=-z$ for 
$z\in{\Bbb C}/{\Bbb Z}+\tau{\Bbb Z}$.
\par
Let ${\Bbb Z}_{2}$ be a group of order $2$
with generator $\theta$. Then ${\Bbb Z}_{2}$
acts on the spaces $\widetilde{S}$, $T$, and 
$\widetilde{S}\times T$ by identifying
$\theta$ with $\iota_{S}$, $-1_{T}$ and
$\iota_{S}\times(-1_{T})$, respectively.

\begin{defn}
For an Enriques surface $S$ and 
an elliptic curve $T$, define 
$$
X_{(S,T)}:=\widetilde{S}\times T/{\Bbb Z}_{2}.
$$
\end{defn}

Since $\iota_{S}\times(-1_{T})$ has no fixed points,
$X_{(S,T)}$ is a smooth Calabi-Yau threefold.
Let 
$p_{1}\colon X_{(S,T)}\to S=
\widetilde{S}/{\Bbb Z}_{2}$ and let
$p_{2}\colon X_{(S,T)}\to{\Bbb P}^{1}/{\Bbb Z}_{2}$
be the natural projections. 
Then $p_{1}$ is an elliptic fibration with constant
fiber $T$, and $p_{2}$ is a $K3$ fibration with
constant fiber $\widetilde{S}$.
After Ferrara-Harvey-Strominger-Vafa \cite{FHSV95}, 
the Calabi-Yau threefold $X_{(S,T)}$ 
is called the {\it FHSV threefold} 
associated with $(S,T)$. We have
\begin{equation}
\chi(X_{(S,T)})
=
\frac{1}{2}\,\chi(\widetilde{S}\times T)
=
\frac{1}{2}\,\chi(\widetilde{S})\chi(T)
=0.
\end{equation}

\subsection
{}
{\bf The moduli space of FHSV threefolds}
\par
The period of an Enriques surface $S$ is defined
as the period of $(\widetilde{S},\iota_{S})$
and lies in the quotient space
$\Omega/\Gamma$, where
$\Omega$ is a symmetric bounded domain of type IV
of dimension $10$ and where $\Gamma$ is
an arithmetic subgroup of ${\rm Aut}(\Omega)$.
The period of $S$ is denoted by $[S]\in\Omega/\Gamma$.
There exists a $\Gamma$-invariant divisor 
$D\subset\Omega$, called the discriminant locus,
such that $(\Omega\setminus D)/\Gamma$
is a coarse moduli space of Enriques surfaces
via the period map. We refer to
e.g. \cite[Chap.\,8, Sects.\,19-21]{BPV84}
for more details about the moduli space of
Enriques surfaces.
\par
In \cite{Borcherds96},
Borcherds constructed an automorphic form $\Phi$
on $\Omega$ for $\Gamma$ of weight $4$ with
${\rm div}(\Phi)=D$. The automorphic form $\Phi$
is called the {\it Borcherds $\Phi$-function}.
Let $B_{\Omega}$ be the Bergman kernel function
of $\Omega$.
The Petersson norm of the Borcherds $\Phi$-function
is the $\Gamma$-invariant $C^{\infty}$ function 
on $\Omega$ defined 
as 
$$
\|\Phi\|^{2}:=B_{\Omega}^{4}|\Phi|^{2}.
$$
By the $\Gamma$-invariance of $\|\Phi\|^{2}$,
it descends to a function on $\Omega/\Gamma$,
denoted again by $\|\Phi\|^{2}$.
Then $\|\Phi([S])\|^{2}$ is the value of
the Petersson norm of the Borcherds $\Phi$-function
at the period point of an Enriques surface $S$.
We refer to \cite{Borcherds96}, \cite{Yoshikawa04}
for more details about the Borcherds $\Phi$-function.
\par
For an elliptic curve $T\cong E_{\tau}$,
the period of $T$ is defined as
the $SL_{2}({\Bbb Z})$-orbit of $\tau\in{\Bbb H}$
and is denoted by $[T]\in{\Bbb H}/SL_{2}({\Bbb Z})$. 
The quotient space ${\Bbb H}/SL_{2}({\Bbb Z})$ is
a coarse moduli space of elliptic curves 
via the period map.
Let 
$$
\Delta(\tau):=q\,\prod_{n=1}^{\infty}(1-q^{n})^{24},
\qquad
q:=\exp(2\pi\sqrt{-1}\tau)
$$
be the Jacobi $\Delta$-function, which is
a unique cusp form of weight $12$.
The Petersson norm
of the Jacobi $\Delta$-function is 
a $SL_{2}({\Bbb Z})$-invariant $C^{\infty}$ function
on ${\Bbb H}$ defined as
$$
\|\Delta(\tau)\|^{2}:=(\det{\rm Im}\,\tau)^{12}
|\Delta(\tau)|^{2}.
$$ 
By the $SL_{2}({\Bbb Z})$-invariance of
$\|\Delta\|^{2}$, it descends to a function 
on ${\Bbb H}/SL_{2}({\Bbb Z})$.
Then $\|\Delta([T])\|^{2}$ is the value of
the Petersson norm of the Jacobi $\Delta$-function
at the period point of an elliptic curve $T$.

\begin{thm}
The analytic space
$[(\Omega\setminus D)/\Gamma]\times
[{\Bbb H}/SL_{2}({\Bbb Z})]$
is a coarse moduli space of FHSV threefolds.
\end{thm}

\begin{pf}
Since $(\Omega\setminus D)/\Gamma$ is a coarse
moduli space of Enriques surfaces 
\cite[Chap.\,8, Ths.\,21.2 and 21.4]{BPV84} and 
since ${\Bbb H}/SL_{2}({\Bbb Z})$ is a coarse moduli
space of elliptic curves via the elliptic $j$-function, 
it suffices to prove
that $X_{(S,T)}\cong X_{(S',T')}$ if and only if
$S\cong S'$ and $T\cong T'$. This statement 
follows from \cite[Sect.\,3]{Beauville83}.
\end{pf}

\subsection
{}
{\bf A Conjecture of Harvey-Moore}
\par
Following \cite[Sect.\,V]{HarveyMoore98}
and \cite[Sect.\,8.1]{Yoshikawa04},
we interpret a result of the third-named
author \cite[Th.\,8.3]{Yoshikawa04}
in terms of the BCOV torsion of FHSV threefolds.
The following formula was conjectured by
Harvey-Moore \cite[Eq.\,(4.9)]{HarveyMoore98}.

\begin{thm}
There exists a constant $C$ such that
for every Enriques surface $S$ and
for every elliptic curve $T$,
$$
\tau_{\rm BCOV}(X_{(S,T)})
=
C\,\|\Phi([S])\|^{2}\,\|\Delta([T])\|^{2}.
$$
\end{thm}

\par
For the proof of Theorem 13.3, we need some
intermediary results.
Let $H^{2}_{+}(\widetilde{S},{\Bbb Z})$
be the invariant subspace of
$H^{2}(\widetilde{S},{\Bbb Z})$ 
with respect to the $\iota_{S}$-action.
Let $H\in H^{2}_{+}(\widetilde{S},{\Bbb Z})$ 
be an $\iota_{S}$-invariant K\"ahler class
on $\widetilde{S}$, and 
let ${\bf v}\in H^{2}(T,{\Bbb Z})$ 
be the generator with $\int_{T}{\bf v}=1$. 
Let $\pi\colon\widetilde{S}\times T\to X_{(S,T)}$
be the natural projection. 
We define $\kappa\in H^{2}(X_{(S,T)},{\Bbb Z})$ 
to be the K\"ahler class on $X_{(S,T)}$
such that $\pi^{*}\kappa=H+{\bf v}$.
By \cite{Yau78}, there exists a unique
Ricci-flat K\"ahler form $\gamma=\gamma_{\kappa}$ 
on $X_{(S,T)}$ with K\"ahler class $\kappa$.
By \cite{Yau78} again, 
there exist a unique Ricci-flat K\"ahler form 
$\gamma_{H}$ on $\widetilde{S}$ and 
a unique Ricci-flat K\"ahler form $\gamma_{T}$
on $T$ such that
$$
\pi^{*}\gamma_{\kappa}
=\gamma_{H}+\gamma_{T},
\qquad
[\gamma_{H}]=H,
\qquad
[\gamma_{T}]={\bf v}.
$$
\par
Let $\langle\cdot,\cdot\rangle$ denote
the cup-product pairing on 
$H^{2}(\widetilde{S},{\Bbb Z})$.
Since $\int_{T}{\bf v}=1$ and
$\langle a,b\rangle=\int_{\widetilde{S}}a\wedge b$
for $a,b\in H^{2}(\widetilde{S},{\Bbb Z})$, 
we get
\begin{equation}
{\rm Vol}(X_{(S,T)},\gamma)
=
\frac{1}{2}\int_{\widetilde{S}\times T}
\frac{(H+{\bf v})^{3}}{(2\pi)^{3}3!}
=
\frac{1}{2^{5}\pi^{3}}\,\langle H,H\rangle.
\end{equation}
By the Ricci-flatness of $\gamma$, Remark 4.2,
and (13.1), we get
\begin{equation}
{\cal A}(X_{(S,T)},\gamma)
=
{\rm Vol}(X_{(S,T)},\gamma)^{\chi(X_{(S,T)})/12}
=1.
\end{equation}

\begin{lem}
The following identity holds:
$$
{\rm Vol}_{L^{2}}
(H^{2}(X_{(S,T)},{\Bbb Z}),\kappa)
=
\frac{\langle H,H\rangle}{2^{35}\pi^{33}}.
$$
\end{lem}

\begin{pf}
Let $H^{2}_{+}(\widetilde{S}\times T,{\Bbb Z})$ 
be the invariant subspace of 
$H^{2}(\widetilde{S}\times T,{\Bbb Z})$
with respect to the $\iota_{S}\times(-1_{T})$-action.
Similarly, let $H^{2}_{+}(T,{\Bbb Z})$ be
the invariant subspace of
$H^{2}(T,{\Bbb Z})$ with respect to 
the $-1_{T}$-action. We have
\begin{equation}
\pi^{*}H^{2}(X_{(S,T)},{\Bbb Z})_{\rm fr}
=
H^{2}_{+}(\widetilde{S}\times T,{\Bbb Z})
=
H^{2}_{+}(\widetilde{S},{\Bbb Z})\oplus
H^{2}_{+}(T,{\Bbb Z})
=
H^{2}_{+}(\widetilde{S},{\Bbb Z})\oplus
{\Bbb Z}\,{\bf v}.
\end{equation}
By \cite[Chap.\,8, Lemma 15.1 (iii)]{BPV84}, 
there exists an integral basis
$\{{\bf e}_{1},\ldots,{\bf e}_{10}\}$ of
$H^{2}_{+}(\widetilde{S},{\Bbb Z})$
such that
\begin{equation}
\det(\langle{\bf e}_{i},{\bf e}_{j}\rangle
)_{1\leq i,j\leq10}
=-2^{10}.
\end{equation}
By (13.4), we fix the basis 
$\{\bar{\bf e}_{1},\ldots,\bar{\bf e}_{10},
\bar{\bf v}\}$
of $H^{2}(X_{(S,T)},{\Bbb Z})_{\rm fr}$ 
such that
$$
\pi^{*}(\bar{\bf e}_{i})={\bf e}_{i}
\quad
(1\leq i\leq 10),
\qquad
\pi^{*}(\bar{\bf v})={\bf v}.
$$
Recall that the cubic form $c=c_{X_{(S,T)}}$ 
on $H^{2}(X_{(S,T)},{\Bbb Z})_{\rm fr}$ 
was defined in Sect.\,4.4. Then we get
$$
c(\bar{\bf e}_{i},\bar{\bf v},\kappa)
=
\frac{1}{2(2\pi)^{3}}\int_{\widetilde{S}\times T}
{\bf e}_{i}\wedge{\bf v}\wedge\pi^{*}\kappa
=
\frac{1}{2(2\pi)^{3}}\int_{\widetilde{S}\times T}
{\bf e}_{i}\wedge{\bf v}\wedge(H+{\bf v})
=
\frac{1}{2(2\pi)^{3}}\langle{\bf e}_{i},H\rangle,
$$
$$
c(\bar{\bf e}_{i},\bar{\bf e}_{j},\kappa)
=
\frac{1}{2(2\pi)^{3}}\int_{\widetilde{S}\times T}
{\bf e}_{i}\wedge{\bf e}_{j}\wedge\pi^{*}\kappa
=
\frac{1}{2(2\pi)^{3}}\int_{\widetilde{S}\times T}
{\bf e}_{i}\wedge{\bf e}_{j}\wedge(H+{\bf v})
=
\frac{1}{2(2\pi)^{3}}
\langle{\bf e}_{i},{\bf e}_{j}\rangle,
$$
$$
c(\bar{\bf e}_{i},\kappa,\kappa)
=
\frac{1}{2(2\pi)^{3}}\int_{\widetilde{S}\times T}
{\bf e}_{i}\wedge(\pi^{*}\kappa)^{2}
=
\frac{1}{2(2\pi)^{3}}\int_{\widetilde{S}\times T}
{\bf e}_{i}\wedge(H+{\bf v})^{2}
=
\frac{1}{(2\pi)^{3}}\langle{\bf e}_{i},H\rangle,
$$
$$
c(\bar{\bf v},\bar{\bf v},\kappa)
=
\frac{1}{2(2\pi)^{3}}\int_{\widetilde{S}\times T}
{\bf v}\wedge{\bf v}\wedge\pi^{*}\kappa
=
\frac{1}{2(2\pi)^{3}}\int_{\widetilde{S}\times T}
{\bf v}\wedge{\bf v}\wedge(H+{\bf v})
=0,
$$
$$
c(\bar{\bf v},\kappa,\kappa)
=
\frac{1}{2(2\pi)^{3}}\int_{\widetilde{S}\times T}
{\bf v}\wedge(\pi^{*}\kappa)^{2}
=
\frac{1}{2(2\pi)^{3}}\int_{\widetilde{S}\times T}
{\bf v}\wedge(H+{\bf v})^{2}
=
\frac{1}{2(2\pi)^{3}}\langle H,H\rangle,
$$
$$
c(\kappa,\kappa,\kappa)
=
\frac{1}{2(2\pi)^{3}}\int_{\widetilde{S}\times T}
(\pi^{*}\kappa)^{3}
=
\frac{1}{2(2\pi)^{3}}\int_{\widetilde{S}\times T}
(H+{\bf v})^{3}
=
\frac{3}{2(2\pi)^{3}}\langle H,H\rangle.
$$
By Lemma 4.12 and these formulae, we get
$$
(2\pi)^{3}\langle
\bar{\bf e}_{i},\bar{\bf e}_{j}
\rangle_{L^{2},\kappa}
=
\frac{3}{2}
\frac{c(\bar{\bf e}_{i},\kappa,\kappa)
c(\bar{\bf e}_{j},\kappa,\kappa)}
{c(\kappa,\kappa,\kappa)}
-
c(\bar{\bf e}_{i},\bar{\bf e}_{j},\kappa)
=
\frac{\langle{\bf e}_{i},H\rangle\,
\langle{\bf e}_{j},H\rangle}{\langle H,H\rangle}
-
\frac{1}{2}\langle{\bf e}_{i},{\bf e}_{j}\rangle,
$$
$$
(2\pi)^{3}\langle
\bar{\bf e}_{i},\bar{\bf v}
\rangle_{L^{2},\kappa}
=
\frac{3}{2}
\frac{c(\bar{\bf e}_{i},\kappa,\kappa)
c(\bar{\bf v},\kappa,\kappa)}
{c(\kappa,\kappa,\kappa)}
-
c(\bar{\bf e}_{i},\bar{\bf v},\kappa)
=
\frac{1}{2}\frac{\langle{\bf e}_{i},H\rangle\,
\langle H,H\rangle}{\langle H,H\rangle}
-
\frac{1}{2}\langle{\bf e}_{i},H\rangle
=0,
$$
$$
(2\pi)^{3}\langle
\bar{\bf v},\bar{\bf v}
\rangle_{L^{2},\kappa}
=
\frac{3}{2}
\frac{c(\bar{\bf v},\kappa,\kappa)
c(\bar{\bf v},\kappa,\kappa)}
{c(\kappa,\kappa,\kappa)}
-
c(\bar{\bf v},\bar{\bf v},\kappa)
=
\frac{1}{4}\frac{\langle H,H\rangle\,
\langle H,H\rangle}{\langle H,H\rangle}
-0
=
\frac{1}{4}\,\langle H,H\rangle,
$$
which yields that
\begin{equation}
\begin{aligned}
\,&
{\rm Vol}_{L^{2}}(H^{2}(X_{(S,T)},{\Bbb Z}),\kappa)
\\
&
=
\det
\begin{pmatrix}
\langle
\bar{\bf e}_{i},\bar{\bf e}_{j}
\rangle_{L^{2},\kappa}
&
\langle
\bar{\bf e}_{i},\bar{\bf v}
\rangle_{L^{2},\kappa}
\\
\langle
\bar{\bf e}_{i},\bar{\bf v}
\rangle_{L^{2},\kappa}
&
\langle
\bar{\bf v},\bar{\bf v}
\rangle_{L^{2},\kappa}
\end{pmatrix}
\\
&
=
(2\pi)^{-33}2^{-10}\,\frac{\langle H,H\rangle}{4}\,
\det
\left(
\langle{\bf e}_{i},{\bf e}_{j}\rangle
-
2\frac{\langle{\bf e}_{i},H\rangle\,
\langle{\bf e}_{j},H\rangle}{\langle H,H\rangle}
\right)_{1\leq i,j\leq10}.
\end{aligned}
\end{equation}
\par
Define a $10\times10$ symmetric matrix $A$ by
$A=(\langle{\bf e}_{i},{\bf e}_{j}\rangle)$.
Write $H=\sum_{i=1}^{10}h_{i}\,{\bf e}_{i}$
and define a column vector 
${\bf h}\in{\Bbb Z}^{10}$ by ${\bf h}=(h_{i})$.
We set
$$
B:=A-2\,\frac{(A{\bf h})\cdot({}^{t}{\bf h}A)}
{{}^{t}{\bf h}A{\bf h}}.
$$
Since $A$ is invertible and since
${}^{t}{\bf h}\,A\,{\bf h}=\langle H,H\rangle>0$,
we get the decomposition
${\Bbb R}^{10}={\Bbb R}\,{\bf h}\oplus 
(A\,{\bf h})^{\perp}$. 
Since $B\,{\bf h}=-A\,{\bf h}$ and 
$B\,{\bf x}=A\,{\bf x}$ for
${\bf x}\in(A\,{\bf h})^{\perp}$, 
we get $\det B=-\det A=2^{10}$ by (13.5), 
which, together with (13.6), yields that
$$
{\rm Vol}_{L^{2}}(H^{2}(X_{(S,T)},{\Bbb Z}),\kappa)
=
(2\pi)^{-33}2^{-10}\,
\frac{\langle H,H\rangle}{4}\,\det B
=
\frac{\langle H,H\rangle}{2^{35}\pi^{33}}.
$$
This completes the proof of Lemma 13.4.
\end{pf}

Let $\square_{H}$ (resp. $\square_{T}$)
be the $\bar{\partial}$-Laplacain of 
$(\widetilde{S},\gamma_{H})$
(resp. $(T,\gamma_{T})$) acting on
$C^{\infty}(\widetilde{S})$
(resp. $C^{\infty}(T)$).
We define
$$
A^{\pm}(\widetilde{S})
:=\{f\in C^{\infty}(\widetilde{S});\,
\iota_{S}^{*}f=\pm f\},
\qquad
A^{\pm}(T):=\{f\in C^{\infty}(T);\,
(-1_{T})^{*}f=\pm f\}.
$$
Since $\iota_{S}$ (resp. $-1_{T}$)
preserves $\gamma_{H}$ (resp. $\gamma_{T}$), 
$\square_{H}$ commutes 
with the $\iota_{S}$-action on 
$C^{\infty}(\widetilde{S})$ and
$\square_{T}$ commutes 
with the $(-1)_{T}$-action on $C^{\infty}(T)$.
Hence $\square_{H}$ preserves 
$A^{\pm}(\widetilde{S})$, and
$\square_{T}$ preserves $A^{\pm}(T)$. 
We set 
$$
\square_{H}^{\pm}:=\square_{H}
|_{A^{\pm}(\widetilde{S})},
\qquad
\square_{T}^{\pm}:=\square_{T}
|_{A^{\pm}(T)}.
$$
\par
Let $\zeta_{H}^{\pm}(s)$ (resp. $\zeta_{T}^{\pm}(s)$)
be the spectral zeta function of 
$\square_{H}^{\pm}$ (resp. $\square_{T}^{\pm}$).
Then $\zeta_{H}^{\pm}(s)$ and
$\zeta_{T}^{\pm}(s)$ converges absolutely
for ${\rm Re}\,s\gg0$, 
they extend meromorphically 
to the complex plane $\Bbb C$, 
and they are holomorphic at $s=0$.

\begin{lem}
The following identity holds
$$
\log{\cal T}_{\rm BCOV}(X_{(S,T)},\gamma)
=
-24\,(\zeta^{+}_{T})'(0)
-8\left\{(\zeta^{+}_{H})'(0)-(\zeta^{-}_{H})'(0)
\right\}.
$$
\end{lem}

\begin{pf}
See
\cite[Sect.\,V]{HarveyMoore98}, in particular 
\cite[Eqs.\,(5.3), (5.9), (5.10)]{HarveyMoore98}.
\end{pf}

\begin{rem}
The signs in \cite[Eqs.\,(5.10), (5.11)]{HarveyMoore98}
are not correct. 
In \cite[Eqs.\,(5.10), (5.11)]{HarveyMoore98},
the formula
$\log\det'\square^{\pm}_{H}=(\zeta^{\pm}_{H})'(0)$
was used, while the correct formula is
$\log\det'\square^{\pm}_{H}=-(\zeta^{\pm}_{H})'(0)$.
\end{rem}

\begin{lem}
There exists a constant $C_{0}$ such that 
for every Enriques surface $S$ and 
for every K\"ahler class $H$ on $\widetilde{S}$,
the following identity holds
$$
8\left\{(\zeta^{+}_{H})'(0)-(\zeta^{-}_{H})'(0)
\right\}
+4\,\log\langle H,H\rangle
=
-\log\|\Phi([S])\|^{2}+C_{0}.
$$
\end{lem}

\begin{pf}
The result follows from
\cite[Eq.\,(8.3)]{Yoshikawa04} and
\cite[Lemma 4.3, Eq.\,(4.4)]{Yoshikawa05b}.
\end{pf}

\begin{lem}
There exists a constant $C_{1}$ such that 
for every elliptic curve $T$,
$$
24\,(\zeta_{T}^{+})'(0)
=
-\log\|\Delta([T])\|^{2}+C_{1}.
$$
\end{lem}

\begin{pf}
Since $\zeta^{+}_{T}(s)=\zeta^{-}_{T}(s)$
by \cite[p.166 l.8 and l.10]{RaySinger73}
and since $\zeta^{+}_{T}(s)+\zeta^{-}_{T}(s)$
is the spectral zeta function of $\square_{T}$,
the result follows from the Kronecker limit
formula. See e.g. \cite[Th.\,4.1]{RaySinger73}
or \cite[Th.\,13.1]{BismutBost90}.
\end{pf}

\subsection
{}
{\bf Proof of Theorem 13.3}
\par
By Lemmas 13.5, 13.7, 13.8, we get
\begin{equation}
\log{\cal T}_{\rm BCOV}(X_{(S,T)},\gamma)
=
\log(\|\Phi([S])\|^{2}\,\|\Delta([T])\|^{2})
+4\,\log\langle H,H\rangle
-C_{0}-C_{1}.
\end{equation}
By (13.2), (13.3), (13.7) and Lemma 13.4,
we get
$$
\begin{aligned}
\,&
\tau_{\rm BCOV}(X_{(S,T)},\gamma)
\\
&
=
{\rm Vol}(X_{(S,T)},\frac{\gamma}{2\pi})^{-3}\,
{\rm Vol}_{L^{2}}
(H^{2}(X_{(S,T)},{\Bbb Z}),[\gamma])^{-1}\,
{\cal A}(X_{(S,T)},\gamma)\,
{\cal T}_{\rm BCOV}(X_{(S,T)},\gamma)
\\
&
=
\left(
\frac{\langle H,H\rangle}{2^{5}\pi^{3}}
\right)^{-3}
\cdot
\left(
\frac{\langle H,H\rangle}{2^{35}\pi^{33}}
\right)^{-1}
\cdot
1\cdot
\frac{\|\Phi([S])\|^{2}\,\|\Delta([T])\|^{2}
\langle H,H\rangle^{4}}{e^{C_{0}+C_{1}}}
\\
&
=
C\,\|\Phi([S])\|^{2}\,\|\Delta([T])\|^{2},
\end{aligned}
$$
where we set $C=2^{50}\pi^{42}\,e^{-C_{0}-C_{1}}$.
This completes the proof of Theorem 13.3.
\qed

%%%%%%%%%%%%%%%%%%%%%%%%%%%%%%%%%%%%%%%%%%%%%%%%%%%%%%%%%%%
%
%                 References
%
%%%%%%%%%%%%%%%%%%%%%%%%%%%%%%%%%%%%%%%%%%%%%%%%%%%%%%%%%%%

\end{document}